\documentclass{amsart}

\usepackage[foot]{amsaddr}
\makeatletter
\def\BState{\State\hskip-\ALG@thiprstlm}
\renewcommand{\email}[2][]{%
	\ifx\emails\@empty\relax\else{\g@addto@macro\emails{,\space}}\fi%
	\@ifnotempty{#1}{\g@addto@macro\emails{\textrm{(#1)}\space}}%
	\g@addto@macro\emails{#2}%
}

\makeatother

\usepackage{color}

\usepackage{tikz}

\usepackage{hhline}
\usepackage{multirow}
\usepackage{multicol}

\usepackage{amsfonts}
\usepackage{amstext}
\usepackage{wrapfig}

\usepackage{stmaryrd}
\usepackage{amsthm}
\usepackage{amsmath}

\usepackage{mathtools}
\usepackage{amssymb}
\usepackage{commath}
\usepackage{bigints}
\usepackage{romanbar}
\usepackage{algorithm2e}
\usepackage{enumerate}

\usepackage{comment}

\usepackage[shortlabels]{enumitem}

\newtheorem{remark}{Remark}

\newcommand{\vy}{\mathbf{y}}
\newcommand{\vu}{\mathbf{u}}
\newcommand{\vv}{\mathbf{v}}
\newcommand{\vw}{\mathbf{w}}
\newcommand{\vx}{\mathbf{x}}
\newcommand{\vz}{\mathbf{z}}
\newcommand{\vp}{\mathbf{p}}
\newcommand{\vq}{\mathbf{q}}
\newcommand{\vf}{\mathbf{f}}
\newcommand{\vn}{\mathbf{n}}

\newcommand{\vnu}{\boldsymbol{\nu}}
\newcommand{\vvarphi}{\boldsymbol{\varphi}}
\newcommand{\vphi}{\boldsymbol{\phi}}

\newcommand{\jump}[1]{[ #1 ]} 
\newcommand{\avrg}[1]{\left\{ #1 \right\}} 
\newcommand{\Th}{\mathcal{T}_h} 
\newcommand{\Eh}{\mathcal{E}_h} 
\newcommand{\Ehg}{\mathcal{E}_h^{\rm grad}}
\newcommand{\Ehv}{\mathcal{E}_h^{\rm val}}
\newcommand{\E}{\mathcal{E}} 
\newcommand{\K}{T} 
\newcommand{\V}{\mathbb{V}} 
\newcommand{\W}{\mathbb{W}} 
\newcommand{\h}{ {\rm h}} 

\newcommand{\Ghg}{\Gamma_h^{\rm grad}} 
\newcommand{\Ghv}{\Gamma_h^{\rm val}} 
 
\newcommand{\I}{\normalfont{\Romanbar{1}}} 
\newcommand{\II}{\normalfont{\Romanbar{2}}} 
 
\newcommand{\tr}{{\rm tr}} 
\newcommand{\di}{\mathop{\rm div}\nolimits} 
 
\newcommand{\A}{\mathbb{A}}

\DeclareMathOperator*{\argmin}{arg\,min}

\def\restriction#1#2{\mathchoice
	{\setbox1\hbox{${\displaystyle #1}_{\scriptstyle #2}$}
		\restrictionaux{#1}{#2}}
	{\setbox1\hbox{${\textstyle #1}_{\scriptstyle #2}$}
		\restrictionaux{#1}{#2}}
	{\setbox1\hbox{${\scriptstyle #1}_{\scriptscriptstyle #2}$}
		\restrictionaux{#1}{#2}}
	{\setbox1\hbox{${\scriptscriptstyle #1}_{\scriptscriptstyle #2}$}
		\restrictionaux{#1}{#2}}}
\def\restrictionaux#1#2{{#1\,\smash{\vrule height .8\ht1 depth .85\dp1}}_{\,#2}}

\title[Large deformations of plates]{Numerical approximations of thin structure deformations}
\author{Andrea Bonito}
\address[Andrea Bonito]{Department of Mathematics, Texas A\&M University, College Station, TX 77845, USA.}
\email{bonito@math.tamu.edu}
\author{Diane Guignard}
\address[Diane Guignard]{Department of Mathematics and Statistics
	\\ University of Ottawa,
	Ottawa, ON K1N 6N5, Canada.}
\email{dguignar@uottawa.ca}
\author{Angelique Morvant}
\address[Angelique Morvant]{Department of Mathematics, Texas A\&M University, College Station, TX 77845, USA.}
\email{mae4102@tamu.edu}

\thanks{AB was partially supported by the NSF Grant DMS-2110811. DG was partially supported by NSERC Grant RGPIN-2021-04311. AM was partially supported by the NSF Grants DMS-1817691 and DMS-2110811.}

\dedicatory{In memory of Roland Glowinski} 
\date{\today}

\begin{document}

\maketitle

\begin{abstract} 
	We review different (reduced) models for thin structures using bending as principal mechanism to undergo large deformations. Each model consists in the minimization of a fourth order energy, potentially subject to a nonconvex constraint. Equilibrium deformations are approximated using local discontinuous Galerkin (LDG) finite elements. The design of the discrete energies relies on a discrete Hessian operator defined on discontinuous functions with better approximation properties than the piecewise Hessian. Discrete gradient flows are put in place to drive the minimization process. They are chosen for their robustness and ability to preserve the nonconvex constraint. Several numerical experiments are presented to showcase the large variety of shapes that can be achieved with these models.
\end{abstract}

\vspace*{0.2cm}

\vspace*{0.2cm}
\subjclass{\textbf{AMS Subject Classification}: 65N12, 65N30, 74K20, 74G65
	
\section{Introduction}

Deformations of thin materials are widely observed in nature, from the snapping of the venus flytrap \cite{FSDM2005, Speck2020} to the natural growth of soft tissues in leaves and flowers \cite{goriely2005differential, yavari2010geometric, BK2014, LM2021}. They also appear in a variety of man-made applications \cite{klein2007shaping, KVS2011}. Some thin structures, called {\em bilayers}, are made of two thin layers of different materials that react differently to external stimuli (changes in temperature or humidity, electrical or chemical stimuli, etc.). Examples include the bimetal strips in thermostats, microactuators  \cite{JSI2000, BALG2010}, and plywood panels in climate-responsive architectures \cite{Menges2015} inspired from the morphology of conifer cones. Bending of thin sheets can also occur when folding is present; examples of this include origami and flexible structures \cite{SLPSK2015, LPFJ2021}. We refer to the review paper \cite{LM2021} for additional examples.

All the phenomena listed above can be modelled as thin limits of hyperelastic materials. A thin structure of thickness $s$ endowed with a hyperelastic energy $E_s$ is approximated by its two dimensional midplane endowed with a limiting energy as the thickness vanishes. A hierarchy of reduced models for elastic deformations is described in \cite{FJM2006}, in which the type of the model depends on the scaling of the elastic energy $E_s$ with the thickness $s$ of the plate. The scaling $E_s \sim s$ corresponds to stretching, the scaling $E_s \sim s^3$ corresponds to bending, and the scaling $E_s \sim s^5$ corresponds to the von K\'arm\'an theory of plate bending. Each of these reduced energies are the $\Gamma$-limit of $E_s$. We refer to \cite{DR1995} for the membrane theory, to \cite{FJM2002a, FJM2002b, FJMM2003}  for the bending theory, and to \cite{FJM2002c, FJM2006} for the von K\'arm\'an theory. In particular, cluster points of sequences of minimizers of $\{ E_s\}_{s>0}$ are minimizers of the $\Gamma$-limit.

In this work, we are mainly interested in the bending regime, namely when the hyperelastic three dimensional energy
\begin{equation} \label{eqn:Es_intro}
	E_s(
	\vu
	) = \int_{\Omega\times\big(-\frac{s}{2},\frac{s}{2}\big)}W\big(\nabla\vu(\vx)G^{-\frac{1}{2}}(\vx)\big){\rm d}\vx
\end{equation}
scales like the cube of the thickness of the plate $\Omega \times (-\frac s2,\frac s2)$. Here, $\Omega$ is the midplane of the thin structure, $\vu: \Omega\times (-\frac s2,\frac s2) \to \mathbb{R}^3$ is a deformation of the plate, $G$ is a prescribed Riemannian metric (deformations satisfying the target metric $G$ are stress-free), and $W$ is some energy density function satisfying \cite{FJM2002b} 
\begin{itemize}
	\item $W\in C^0(\mathbb{R}^{3\times 3})$, $W\in C^2$ in a neighborhood of ${\rm SO}(3)$;
	\item $W$ is frame indifferent: $W(F)=W(RF)$ for all $F\in\mathbb{R}^{3\times 3}$ and all $R\in{\rm SO}(3)$;
	\item $W(F)\ge C{\rm dist}^2(F,{\rm SO}(3))$ for all $F\in\mathbb{R}^{3\times 3}$ and $W(F)=0$ if $F\in{\rm SO}(3)$,
\end{itemize}
where ${\rm SO}(3)$ denotes the special orthogonal group of rotations in $\mathbb{R}^3$.
In this work, we restrict our considerations to the St.~Venant-Kirchhoff stored energy, see \eqref{def:E_density}.

The prestrain metric $G$ characterizes the material. From the third property of the energy density, we see that when $G$ is the identity matrix, the minimum of the hyperelastic energy \eqref{eqn:Es_intro} (without boundary conditions or external forces) is zero and is achieved by rigid motions. That is, at equilibrium, the plate is stress free and flat. When the Riemannian curvature tensor of $G$ is not identically zero, then the minimum of the hyperelastic energy is strictly positive \cite{LP2011}. In this case, there is no stress-free configuration. The plate is then said to be non-Euclidean or prestrained and generally exhibits more complex equilibrium shapes. 

In 1850, Kirchhoff \cite{Kirchhoff1850} obtained a reduced energy that can be formally obtained assuming that the material deformation reads
\begin{equation}\label{KL_intro}
	\vu(x_1,x_2,x_3) = \vy(x_1,x_2) + x_3 \vnu(x_1,x_2),
\end{equation}
where $\vnu$ is the normal to the deformed midplane $\vy(\Omega)$. This assumption, used by Love \cite{Love1906}, is usually referred to as (nonlinear) Kirchhoff-Love ansatz in the literature, although it is not due to Kirchhoff \cite{FJM2002b}. When the thickness $s$ of the plate is small, it is convenient to derive reduced models for the deformation $\vy$ of the midplane $\Omega$. 
In \cite{ESK2009}, see also \cite{BK2014}, the Kirchhoff-Love assumption \eqref{KL_intro} is made for prestrained plates.
The resulting model consists of the sum of a stretching and a bending component multiplied by different powers of the plate thickness $s$. We derive a similar energy referred to as the {\em preasymptotic} energy to express that it is obtained upon assuming that $s$ is small but not vanishing. To do this, we follow  \cite{LP2011} and assume that 
\begin{equation}\label{e:G_intro}
	G(x_1,x_2,x_3)=\left(\begin{array}{cc}g(x_1,x_2) & \mathbf{0} \\ \mathbf{0} &1\end{array}\right),
\end{equation}
i.e. $G$ is uniform within the thickness and no stretching occurs in the direction orthogonal to the plane. The function $g:\Omega\rightarrow\mathbb{R}^{2\times 2}$ is assumed to be symmetric and uniformly  positive definite.

The vanishing thickness limit of prestrained materials has been studied in several works. 
When $G$ is the identity matrix, as mentioned above, a reduced energy was obtained formally by Kirchhoff \cite{Kirchhoff1850}. Much later, an ansatz-free derivation was obtained via $\Gamma$-convergence in the seminal work \cite{FJM2002b}. The limiting energy involves the second fundamental form of the midplane deformation and is finite only for isometric immersions.

For the derivation in the general case, namely when $G$ is not the identity matrix, we refer to  \cite{LP2011, BK2014, BLS2016, MS2019, BNPO2022}. The first $\Gamma$-convergence results for prestrained plates were obtained in \cite{LP2011} for target metrics as in \eqref{e:G_intro}, and later extended to more general metrics in \cite{BLS2016}. In both cases, the metric is assumed to be independent of the thickness and uniform throughout the thickness.
For the particular metric given in \eqref{e:G_intro}, the limiting energy involves the second fundamental form of the midplane deformation and is finite only for isometric immersion of the two dimensional target metric $g$. 
In \cite{BGNY2022a} a formal derivation based on a modified Kirchhoff-Love assumption \eqref{KL_intro} of the reduced model is provided.

Bilayer materials are composed of two thin layers of materials with different properties. External (thermal, electrical or chemical) stimuli correspond to a prestrain metric
\begin{equation} \label{eqn:G_bilayer_intro}
	G=(I_3\pm\zeta s N)^T(I_3\pm s\zeta  N),
\end{equation}
where the product $s\zeta $ describes the lattice mismatch between the layers, and $N$ encodes the inhomogeneity and anisotropy of the bilayer. 
When actuated, the material deforms to relax its internal stress and can typically undergo large deformations even with relatively small stimuli.
The reduced model for bilayer plates was derived via $\Gamma$-convergence in \cite{S2007} and formally explained in \cite{BBN2017} for a metric of the form \eqref{eqn:G_bilayer_intro}. The limiting energy penalizes deviations of the second fundamental form of the midplane from a spontaneous curvature $Z$ which depends on the mismatch between the two layers (i.e. on $\zeta$ and $N$). In this model, as in the others when $G=I_3$, the energy is finite only when the deformations are isometries.

The ability of the plate to fold along creases was incorporated in \cite{BBH2021}. It considers hyperelastic materials with $G=I_3$ assumed to be weakened in a neighborhood of a curve $\Sigma\subset\Omega$ modelling a crease. Assuming an asymptotic behavior of the defect width and strength (see \eqref{Hyp:folding_param}), the limiting model reduces to the original model without the crease, except that jumps on the deformations gradients are allowed across the crease without suffering any energetic penalty.  
Note that in principle, any plate - with an isometry constraint or prestrain, a single layer or a bilayer - can have folding. However, only the case of a single layer plate with isometry has been derived rigorously so far.

\vskip1ex\paragraph{\bf Motivations and Novelties.} One of the goals of this paper is to review different models for the deformation of plates, including single layer and bilayer plates, plates with an isometry constraint as well as general prestrain metric constraints, and plates with folding along some curves inside the domain. Although emphasis is made on the pure bending of plates, the preasymptotic regime, in which both the bending and the stretching/shearing of the plate are considered, is also analyzed. To give an intuitive derivation, the reduced two dimensional energies are obtained formally using a (modified) Kirchhoff-Love assumption for the deformation of the three dimensional plate. However, as mentioned above, ansatz-free limiting energies have been obtained via $\Gamma$-convergence for most cases.

A second goal is to collect in one place the numerical methods proposed recently by the authors and collaborators for the approximation of near minimizers of the discrete energies. These methods are based on a \emph{local discontinuous Galerkin} (LDG) approach for the space discretization and a gradient flow for the minimization of the discrete energy. This methodology has already been successfully applied to a wide range of problems, see for instance \cite{BGNY2022a,BNY2022,BBH2021}. 
Moreover, the $\Gamma$-convergence as the mesh size goes to zero of the discrete energies obtained using LDG is analyzed in \cite{BGMpreprint}, \cite{BGNY2022b}, \cite{BNY2022}, and \cite{BBH2021} for the preasymptotic, prestrain, bilayer, and folding case, respectively. Note that other discretizations are also possible. We refer to \cite{BNN2020,BNN2021,bartels2022error} for approaches using the symmetric interior penalty discontinuous Galerkin method, and to \cite{bartels2013,BBN2017,BP2021,bartels2020finite} for a method based on Kirchhoff elements.

This work also offers several novel contributions. We present a recently developed algorithm for the preasymptotic model, extend the existing algorithms to time-dependent input data, present a uniform treatment of folding effects, and introduce an accelerated algorithm that drastically improves the performance of the unconstrained discrete gradient flows. We perform new and challenging simulations highlighting the various types of deformations that can be achieved with the models, as well as the versatility of the code. 

\vskip1ex\paragraph{\bf Outline.} The rest of this paper is organized as follows. In Section~\ref{sec:math_models}, we formally derive the various two dimensional models for the deformation of thin structures and discuss how to incorporate forcing and boundary conditions into these models. In Section~\ref{sec:method}, we present the spatial discretization of the two dimensional energies based on the LDG method, including the introduction of a discrete Hessian operator - the key component of the method. The discrete gradient flows put forward to minimize the discrete energies are described in Section~\ref{sec:GF}, which also includes strategies to produce good initial deformations and to potentially improve the convergence rate of the iterative processes. Section~\ref{sec:results} is devoted to numerical experiments demonstrating the variety of deformations that can be achieved. 

\vskip1ex\paragraph{\bf Notation.} In what follows, we will denote by $I_n$ the identity matrix in $\mathbb{R}^{n\times n}$. Moreover, for a multivariate function, we write $\partial_i$ the partial derivative with respect to the $i^{th}$ variable.
Uppercase letters and bold lowercase letters will be used for matrix-valued and vector-valued functions, respectively, and subindices will denote their components. For instance, the $(i,j)$ component of $M:\mathbb{R}^m\rightarrow\mathbb{R}^{n_1\times n_2}$ will be denoted $M_{ij}:\mathbb{R}^m\rightarrow\mathbb{R}$, $1\le i\le n_1$ and $1\le j\le n_2$, while the $j^{th}$ component of $\vv:\mathbb{R}^m\rightarrow\mathbb{R}^n$ will be denoted $v_j:\mathbb{R}^m\rightarrow\mathbb{R}$, $1\le j\le n$. Furthermore, for $\Phi:\mathbb{R}^m\rightarrow\mathbb{R}^{n\times m}$, $\Phi_j:\mathbb{R}^m\rightarrow\mathbb{R}^n$ will denote the $j^{th}$ row of $\Phi$, $1\le j\le n$.
For $\Phi:\mathbb{R}^m\rightarrow\mathbb{R}^{m\times m}$, the divergence operator is applied row-wise, namely $[\di(\Phi)]_i=\sum_{j=1}^m\partial_j(\Phi_{ij})$ for $i=1,2,\ldots,m$. 
Finally, we use the notation $A:B$ to denote the Euclidean scalar product between two tensors $A,B$ and $|.|$ for the corresponding Frobenius norm. For the particular case of vectors, the Euclidean scalar product between $\vv, \vw \in \mathbb R^n$ is instead denoted by $\vv \cdot \vw$.

\section{Mathematical Models} \label{sec:math_models}

In this section, we provide a formal but intuitive derivation of two dimensional models for three dimensional thin hyperelastic structures. References to rigorous derivations of the bending model with isometry and prestrain constraint and of the bilayer model via $\Gamma$-convergence can be found in Sections~\ref{sec:Limit2D_Kirchhoff}, \ref{sec:prestrain}, and \ref{sec:bilayer}, respectively.
The thin structure is denoted $\Omega_s:=\Omega\times(-\frac s2,\frac s2)\subset\mathbb{R}^3$, where $s>0$ stands for the thickness of the structure and $\Omega\subset\mathbb{R}^2$ is a open, bounded domain with Lipschitz boundary which represents the midplane (see the diagram on the left of Figure~\ref{fig:plate_diagrams}). A deformation of the plate is denoted by $\vu:\Omega_s\subset\mathbb{R}^3\rightarrow\mathbb{R}^3$ while its restriction to the midplane is denoted $\vy:\Omega\subset\mathbb{R}^2\rightarrow\mathbb{R}^3$ (see the diagram on the right of Figure~\ref{fig:plate_diagrams}).
For later use, we also introduce the unit normal vector to the surface $\vy(\Omega)$ at the point $\vy(\vx')$
\begin{equation*}
	\vnu(\vx'):=\frac{\partial_1\vy(\vx')\times\partial_2\vy(\vx')}{|\partial_1\vy(\vx')\times\partial_2\vy(\vx')|}, \qquad\vx'\in\Omega,
\end{equation*}
and the first and second fundamental forms of the surface $\vy(\Omega)$ 
\begin{equation} \label{def:1st2nd}
	\I(\vy):=(\nabla'\vy)^T\nabla'\vy  \quad \mbox{and} \quad \II(\vy):= -(\nabla'\vy)^T\nabla'\vnu = (\partial_{ij}\vy\cdot\vnu)_{i,j=1}^2.
\end{equation}
Here $\nabla'$ denotes the gradient with respect to the variable $\vx' \in \Omega \subset \mathbb R^2$.
\begin{figure}[htbp]
\hspace*{-0.2cm}\includegraphics[width=13cm]{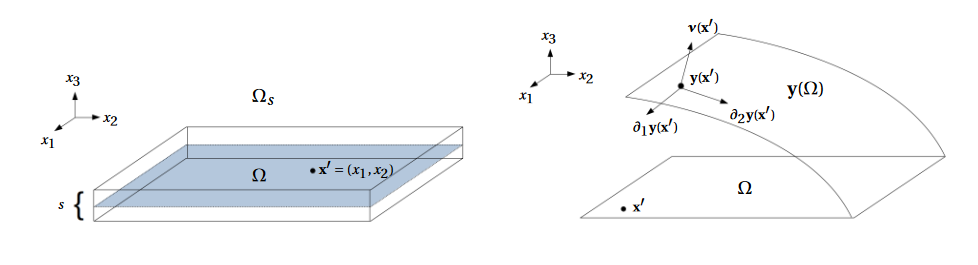}
	\caption{The undeformed reference plate (left) and the deformed midplane (right).} \label{fig:plate_diagrams}
\end{figure}

In all the models presented below, the structure is endowed with an energy favoring deformations with Cauchy-Green strain tensor $(\nabla \vu)^T \nabla \vu$ matching a given target metric $G:\Omega_s\rightarrow\mathbb{R}^{3\times 3}$ (symmetric and positive definite), i.e. for which the strain tensor
\begin{equation}\label{e:strain}
	\epsilon(\nabla\vu):=\frac{(\nabla\vu)^T\nabla\vu-G}{2}
\end{equation}
is as close to zero as possible. Note that $\epsilon(\nabla\vu)\equiv\mathbf{0}$ signifies that the first fundamental form of the deformed surface $\vu(\Omega_s)$ is $G$ and that $\vu$ is an isometric immersion of $G$. Whether a given $n$ dimensional Riemannian manifold has an isometric immersion into $R^N$ is a long standing problem in differential geometry, see for instance \cite{HH2006,CS17}. By the Nash-Kuiper embedding theorem \cite{Nash1954,Kuiper1955a,Kuiper1955b}, there always exists an isometric immersion $\vu:\Omega_s \rightarrow\mathbb{R}^3$ of $G$ in $W^{1,\infty}(\Omega_s)$. However, the deformation $\vu$ cannot be orientation preserving unless the Riemannian curvature tensor of $G$ is identically zero, see \cite{LM2021} and the references therein.

\subsection{Bending of Isotropic Thin Structures}\label{sec:single_iso}

We start by considering thin structures endowed with the isotropic hyperelastic energy
\begin{equation} \label{def:E_3D}
	E_{s}(\vu) := \int_{\Omega_s}W(\nabla \vu(\vx)){\rm d}\vx,
\end{equation}
where $W:\mathbb{R}^{3\times 3}\rightarrow \mathbb{R}$ is the St.~Venant-Kirchhoff stored energy density function defined by
\begin{equation} \label{def:E_density}
	W(F) := \mu|\epsilon(F)|^2+\frac{\lambda}{2}\tr(\epsilon(F))^2, \quad \epsilon(F):= \frac{F^TF-I_3}{2}.
\end{equation}
Here $F:=\nabla\vu$ is the deformation gradient assumed to satisfy the orientation condition $\textrm{det}(\nabla \vu)>0$, $\epsilon(F)$ is the Green-Lagrange strain tensor corresponding to $G=I_3$ in \eqref{e:strain}, and $\lambda$ and $\mu$ are the Lam\'e constants. As we shall see below, deformations for which $E_s(\vu)$ scales like $s$ correspond to a stretching of the midplane (membrane theory), while bending occurs when $E_s(\vu)$ scales like $s^3$ (bending theory). 

Note that in absence of other constraints such as boundary conditions or external forcing, energy \eqref{def:E_3D} is minimized on $\textrm{SO}(3)$, i.e. when the deformations are isometries:
\begin{equation*}
	(\nabla\vu)^T\nabla\vu = I_3 \quad \mbox{in } \Omega_s.
\end{equation*}

\subsection{Limiting Model and Modified Kirchhoff Assumption}\label{sec:Limit2D_Kirchhoff}

Two dimensional models for hyperelastic thin structures are derived in \cite{FJM2002b}. It is shown that when bending is the chief mechanism of deformations, the energy $E_s$ must scale like $s^3$ and
\begin{equation*}
	\mathop{\Gamma\!-\!\lim}_{s \to 0^+} s^{-3} E_s = E
\end{equation*}
where for $\vy: \Omega \rightarrow \mathbb R^3$
\begin{equation}\label{e:2d}
	E(\vy) = \left\lbrace \begin{array}{ll} \displaystyle
		\frac{1}{24}\int_{\Omega}\left(2\mu|\II(\vy(\vx'))|^2+\frac{2\mu\lambda}{2\mu+\lambda}\tr(\II(\vy(\vx')))^2\right){\rm d}\vx', \ & \textrm{when } \I(\vy) = I_2, \\
		+\infty, & \textrm{otherwise}. \end{array}\right.
\end{equation}

Part of the $\Gamma$-convergence argument provided in \cite{FJM2002b} is based on the construction of a recovery sequence given for $\vx'\in\Omega$ and $x_3 \in (-\frac s2,\frac s2)$ by
\begin{equation}\label{Hyp:kirchhoff_2}
	\vu(\vx',x_3)=\vy(\vx')+\left(x_3+ \underbrace{\frac{\lambda}{2\mu+\lambda}\tr\left(\II(\vy(\vx')\right)}_{=:\beta(\vx')} \frac12 x_3^2 \right)\vnu(\vx')
\end{equation}
for a midplane deformation $\vy:\Omega \rightarrow \mathbb R^3$ with finite energy $E(\vy)<\infty$.
This indicates that the two dimensional model \eqref{e:2d} requires that the fibers orthogonal to the midplane are not stretched homogeneously on $\Omega$. This differentiates the theory in \cite{FJM2002b} from the standard Cosserat (or Kirchhoff-Love) theory based on the ansatz
\begin{equation*}
	\label{Hyp:kirchhoff}
	\vu(\vx',x_3)=\vy(\vx')+x_3\vnu(\vx')
\end{equation*}
and leading to the same two dimensional energy improperly scaled with coefficient $\lambda$ instead of $\frac{2\mu\lambda}{2\mu+\lambda}$ in front of the trace term \cite{FJM2002a,FJM2002b}.

\subsubsection{Preasymptotic} \label{sec:preasymptotic}

We follow \cite{BGNY2022a} and provide a formal derivation of the preasymptotic energy for thin structures with non-vanishing but small thickness under assumption \eqref{Hyp:kirchhoff_2} on the deformations. 
We are interested in deformations where bending is the chief mechanism and thus also fix $0< s\ll 1$ and assume that
\begin{equation}\label{e:bending_regim}
	E_s(\vu) \leq \Lambda s^3
\end{equation}
for a constant $\Lambda$ independent of $s$.
We recall that $\nabla'$ stands for the gradient with respect to $\vx' \in \Omega \subset \mathbb R^2$ so that
\begin{equation*}
	\nabla\vu = \left[ \nabla'\vy+x_3\nabla' \vnu + \frac{1}{2}x_3^2\nabla'(\beta\vnu), \vnu+x_3\beta\vnu\right].
\end{equation*}
Because $\partial_i\vy\cdot\vnu=0$ for $i=1,2$ and $\vnu^T\vnu=1$, we have 
\begin{equation*}
	\vnu^T\nabla'\vy=\mathbf{0}, \quad (\nabla'\vy)^T\vnu=\mathbf{0}, \quad  \vnu^T\nabla'\vnu=\mathbf{0}, \quad (\nabla'\vnu)^T\vnu=\mathbf{0},
\end{equation*}
and thus
\begin{align*}
	(\nabla\vu)^T\nabla\vu = &
	\begin{bmatrix}
		(\nabla'\vy)^T\nabla'\vy & \mathbf{0} \\
		\mathbf{0} & 1
	\end{bmatrix}
	+x_3
	\begin{bmatrix}
		(\nabla'\vy)^T\nabla'\vnu+(\nabla'\vnu)^T\nabla'\vy & \mathbf{0} \\
		\mathbf{0} & 2\beta
	\end{bmatrix} \\
	&+x_3^2
	\begin{bmatrix}
		(\nabla'\vnu)^T\nabla'\vnu+\frac{1}{2}(\nabla'\vy)^T\nabla'(\beta\vnu)+\frac{1}{2}(\nabla'(\beta\vnu))^T\nabla'\vy & \frac{1}{2}(\nabla'(\beta\vnu))^T\vnu \\
		\frac{1}{2}\vnu^T\nabla'(\beta\vnu) & \beta^2
	\end{bmatrix} \\
	& + \mathcal{O}(x_3^3). 
\end{align*}
We rewrite the above expression using the first and second fundamental forms \eqref{def:1st2nd} of the deformed midplane $\vy(\Omega)$. In addition, we take advantage of the symmetry of the second fundamental form to obtain
\begin{equation} \label{eqn:Aj}
	(\nabla\vu)^T\nabla\vu-I_3 =  A_1+2x_3A_2+x_3^2A_3 + \mathcal{O}(x_3^3),
\end{equation}
where
\begin{eqnarray*}
	A_1 & := & \begin{bmatrix}
		\I(\vy)-I_2 & \mathbf{0} \\
		\mathbf{0} & 0
	\end{bmatrix}, 
	\qquad A_2  :=   \begin{bmatrix}
		-\II(\vy) & \mathbf{0} \\
		\mathbf{0} & \beta
	\end{bmatrix}, \\
	A_3 & := & \begin{bmatrix}
		(\nabla'\vnu)^T\nabla'\vnu+\frac{1}{2}(\nabla'\vy)^T\nabla'(\beta\vnu)+\frac{1}{2}(\nabla'(\beta\vnu))^T\nabla'\vy & \frac{1}{2}(\nabla'(\beta\vnu))^T\vnu \\
		\frac{1}{2}\vnu^T\nabla'(\beta\vnu) & \beta^2
	\end{bmatrix}.
\end{eqnarray*}

Using this relation in the energy density \eqref{def:E_density} with $F=\nabla\vu$, we obtain
\begin{align*}
	W(\nabla \vu) = &  \frac{\mu}{4}\left(|A_1|^2+4x_3A_1:A_2+4x_3^2|A_2|^2+2x_3^2A_1:A_3+\mathcal{O}(x_3^3)\right) \\ & +\frac{\lambda}{8}\left(\tr(A_1)^2+4x_3\tr(A_1)\tr(A_2)+4x_3^2\tr(A_2)^2+2x_3^2\tr(A_1)\tr(A_3)+\mathcal{O}(x_3^3)\right).
\end{align*}
Note that the terms with odd power of $x_3$ vanish when integrating $x_3$ on $(-\frac s2,\frac s 2)$ to derive an expression of the energy \eqref{def:E_3D}.
Thus, we get
\begin{align*}
	E_s(\vu) =&  \frac{\mu}{4}\int_{\Omega}\left(s|A_1|^2+\frac{s^3}{3}|A_2|^2+\frac{s^3}{6}A_1:A_3\right){\rm d}\vx' \\
	&+\frac{\lambda}{8}\int_{\Omega}\left(s~\tr(A_1)^2+\frac{s^3}{3}\tr(A_2)^2+\frac{s^3}{6}\tr(A_1)\tr(A_3)\right){\rm d}\vx'+\mathcal{O}(s^4).
\end{align*}

The terms with $A_1:A_3$ and $\tr(A_1)\tr(A_3)$ are of higher order. Indeed, thanks to the scaling assumption \eqref{e:bending_regim}, we have
\begin{equation*}
	\int_{\Omega}A_1:A_3{\rm d}\vx'\le \left(\int_{\Omega}|A_1|^2{\rm d}\vx'\right)^{\frac{1}{2}}\left(\int_{\Omega}|A_3|^2{\rm d}\vx'\right)^{\frac{1}{2}}\le \left(\int_{\Omega}|A_3|^2{\rm d}\vx'\right)^{\frac{1}{2}}\left(\frac{4}{\mu}\Lambda\right)^{\frac{1}{2}}s,
\end{equation*}
where $\left(\int_{\Omega}|A_3|^2{\rm d}\vx'\right)^{\frac{1}{2}}$ is independent of $s$.
A similar argument holds for the trace terms. 
As a consequence, in view of the definitions of $A_1$ and $A_2$ in \eqref{eqn:Aj}, the expression \eqref{def:1st2nd} for the first and second fundamental forms, and the definition of $\beta$ in \eqref{Hyp:kirchhoff_2},
we arrive at the final expression for the three dimensional energy per unit volume $s^{-1}E_s(\vu)$:
\begin{equation} \label{eqn:ES_EB}
	s^{-1}E_s(\vu) = \widetilde E^S(\vy) + s^2 \widetilde E^B(\vy) + \mathcal O(s^3),
\end{equation}
where
\begin{eqnarray*}
	\widetilde E^S(\vy) & := & \frac{1}{8}\int_{\Omega}\left(2\mu\left|\I(\vy)-I_2\right|^2+\lambda\tr\left(\I(\vy)-I_2\right)^2\right){\rm d}\vx' \\
	\widetilde E^B(\vy) & := & \frac{1}{24}\int_{\Omega}\left(2\mu\left|\II(\vy)\right|^2+\frac{2\mu\lambda}{2\mu+\lambda}\tr\left(\II(\vy)\right)^2\right){\rm d}\vx' 
\end{eqnarray*}
denote the stretching and bending energies, respectively.

\subsubsection{Limiting Bending Model} \label{sec:reduced2D}

We are now interested in taking the limit when $s \to 0^+$.
We recall that we are assuming \eqref{e:bending_regim}, which, as we shall see, implies that the limiting plate deformation $\vy$ can not stretch nor shear the midplane $\Omega$ but can bend it to reduce its energy.

The starting point to  derive the limiting energy as the thickness $s$ vanishes is the preasymptotic expression \eqref{eqn:ES_EB} for the energy. 
The bending regime condition \eqref{e:bending_regim} implies that $\widetilde E^S(\vy) = 0$ and 
\begin{equation} \label{def:EB_iso}
	E(\vy):= \lim_{s\rightarrow 0^+}s^{-3} E_s(\vu)  =  \frac{1}{24}\int_{\Omega}\left(2\mu\left|\II(\vy)\right|^2+\frac{2\mu\lambda}{2\mu+\lambda}\tr\left(\II(\vy)\right)^2\right){\rm d}\vx'.
\end{equation}
Note that when $\widetilde E^S(\vy) = 0$, 
\begin{equation*}
	\I(\vy) = (\nabla'\vy)^T\nabla'\vy = I_2 \quad \mbox{a.e. in } \Omega,
\end{equation*}
i.e. is an isometry and in particular $\vy \in [W^{1,\infty}(\Omega)]^3$. Furthermore, for isometries the following relations for the second fundamental form hold  (see for instance \cite{bartels2013}):
\begin{equation} \label{def:II_D2y}
	\left|\II(\vy)\right| = \left|\nabla'\nabla'\vy\right| = \tr\left(\II(\vy)\right).
\end{equation}
This allows to rewrite the limiting energy as
\begin{equation} \label{def:EB_iso_simple}
	E(\vy)= \frac{\alpha}{2}\int_{\Omega}|\nabla'\nabla'\vy|^2{\rm d}\vx',
\end{equation}
when $\vy \in [H^2(\Omega) \cap W^{1,\infty}(\Omega)]^3$ is an isometry, i.e. $\I(\vy) = I_2$, and where
$\alpha := \frac{\mu(\mu+\lambda)}{3(2\mu+\lambda)}$.

\subsection{Single Layer with a General Metric Constraint} \label{sec:prestrain}

Without additional constraints such as boundary conditions, the minimum of energy \eqref{def:E_3D} with density \eqref{def:E_density} is achieved by the identity deformation $\vu(\vx) = \vx$, which is an isometry. In particular, the reference flat configuration is stress-free. We now consider prestrained materials characterized by the presence of internal stresses in the flat configuration. Examples are nematic
glasses \cite{modes2010disclination,modes2010gaussian}, natural growth of soft tissues \cite{goriely2005differential,yavari2010geometric}, and manufactured polymer gels \cite{kim2012thermally,klein2007shaping,wu2013three}.

We follow \cite{LP2011,BLS2016} (see also \cite{BGNY2022a}) and modify the energy \eqref{def:E_3D} to read
\begin{equation} \label{def:E_G}
	E_{s}^{\rm pre}(\vu):=\int_{\Omega_s}W\big(\nabla\vu(\vx)G(\vx)^{-\frac12}\big){\rm d}\vx,
\end{equation}
where $G:\Omega_s\rightarrow\mathbb{R}^{3\times 3}$ is a given symmetric positive definite target metric and $G^{-\frac{1}{2}}$ denotes the positive definite symmetric square root of the inverse $G^{-1}$ of $G$. 
Note that the previous case is recovered when $G=I_3$. Furthermore, when $G$ is a metric immersion, i.e. there exists a deformation $\vu:\Omega_s \rightarrow \mathbb R^3$ such that 
\begin{equation} \label{def:G_cstr_3D}
	(\nabla\vu)^T\nabla\vu = G \quad \mbox{a.e. in } \Omega_s,
\end{equation}
we have
$$
\left(\nabla\vu(\vx)G(\vx)^{-\frac12}\right)^T \nabla\vu(\vx)G(\vx)^{-\frac12} = I_3, 
$$
and thus $\vu$ is a minimizer of $E^{\rm pre}_s$ with $E^{\rm pre}_s(\vu)=0$.

Following \cite{ESK2009,LP2011}, we further assume that $G$ has the form
\begin{equation} \label{def:metric_G}
	G(\vx',x_3) = \left(\begin{array}{cc}
		g(\vx') & \mathbf{0} \\
		\mathbf{0} & 1
	\end{array}\right),  \qquad \vx'\in\Omega, ~ x_3 \in (-s/2,s/2),
\end{equation}
where $g:\Omega\rightarrow\mathbb{R}^{2\times 2}$ is symmetric and uniformly  positive definite. In other words, the metric $G$ is independent of the variable $x_3$, is uniform throughout the thickness, and no stretching is allowed in the vertical direction.

The $\Gamma$-convergence of this model towards the two dimensional energy given below is obtained in \cite{LP2011,BLS2016}.
Proceeding as in Section~\ref{sec:preasymptotic}, an intuitive derivation of the two dimensional $\Gamma$-limit can be obtained, using again assumption \eqref{Hyp:kirchhoff_2}.
Indeed, when $0<s\ll 1$, neglecting the $\mathcal O(s^3)$ terms the three dimensional prestrain energy, \eqref{def:E_G} satisfies
\begin{equation} \label{eqn:ES_EB_G}
	s^{-1}E_{s}^{\rm pre}(\vu) = E^S(\vy) + s^2E^B(\vy) =: E_{s}^{\rm pre}(\vy),
\end{equation}
where
\begin{equation}\label{def:ES}
	E^S(\vy)  :=  \frac{1}{8}\int_{\Omega}\left(2\mu\left|g^{-\frac{1}{2}}(\I(\vy)-g)g^{-\frac{1}{2}}\right|^2+\lambda\tr\left(g^{-\frac{1}{2}}(\I(\vy)-g)g^{-\frac{1}{2}}\right)^2\right){\rm d}\vx' 
\end{equation}
and
\begin{equation}\label{def:EB}
	E^B(\vy)  :=  \frac{1}{24}\int_{\Omega}\left(2\mu\left|g^{-\frac{1}{2}}\II(\vy)g^{-\frac{1}{2}}\right|^2+\frac{2\mu\lambda}{2\mu+\lambda}\tr\left(g^{-\frac{1}{2}}\II(\vy)g^{-\frac{1}{2}}\right)^2\right){\rm d}\vx'.
\end{equation}
Note that we slightly abuse the notation by using $E^{\rm pre}_s$ to refer to both the three dimensional and two dimensional prestrain energy. Which energy $E^{\rm pre}_s$ is referring to will be clear from the context. 

For the limit $s^{-3}E_{s}(\vu)$ to be finite as the thickness $s$ goes to zero, the two dimensional deformation $\vy:\Omega \rightarrow \mathbb R^3$ must satisfy the metric constraint
\begin{equation}\label{metric_constraint}
	(\nabla \vy)^T\nabla \vy = g
\end{equation}
and in that case
\begin{equation}\label{E:prestrain-bending}
	\lim_{s \to 0^+} s^{-3}E_{s}(\vu)= \frac{1}{24}\int_{\Omega}\Big(2\mu\left|g^{-\frac{1}{2}}\II(\vy) g^{-\frac{1}{2}}\right|^2+\frac{2\mu\lambda}{2\mu+\lambda}\tr\left(g^{-\frac{1}{2}}\II(\vy) g^{-\frac{1}{2}}\right)^2\Big){\rm d}\vx'.
\end{equation}

\begin{remark}
	Nash's theorem guarantees that there exists an isometric immersion of $g$ into $\mathbb{R}^{10}$ \cite{Gromov1986}, but whether or not there exists an isometric immersion into $\mathbb{R}^3$ depends on $g$. We refer to \cite{HH2006} for a discussion of positive and negative results for metrics with specific properties. Note that the existence of $\vy\in H^2(\Omega)$ satisfying \eqref{metric_constraint} is equivalent to the boundedness condition \eqref{e:bending_regim} on $E_s$. We refer to \cite{LP2011} for the case where $G$ has the form \eqref{def:metric_G} and to \cite{KS2014} for a general Riemannian metric.
\end{remark}

Comparing the limiting energy \eqref{E:prestrain-bending} with \eqref{def:EB_iso}, we find that the isometry case is recovered for $g=I_2$.
However, in the case of isometries, the energy is further reduced to \eqref{def:EB_iso_simple} 
thanks to \eqref{def:II_D2y}. When $g \not = I_2$, the second fundamental form of $\vy$ in \eqref{E:prestrain-bending} cannot be substituted by the Hessian of $\vy$. 
Nevertheless, Proposition~1 in \cite{BGNY2022a} and Proposition~A.1 in \cite{BGNY2022b} guarantee that for $\vy \in [H^2(\Omega)]^3$ and $g\in[H^1(\Omega)\cap L^{\infty}(\Omega)]^{2\times 2}$, up to an additive term only depending on $g$, the limiting energy \eqref{E:prestrain-bending} is equal to 
\begin{equation*}
	\frac{1}{24}\sum_{m=1}^3\int_{\Omega}\left(2\mu\left|g^{-\frac{1}{2}}(\nabla' \nabla' y_m)g^{-\frac{1}{2}}\right|^2+\frac{2\mu\lambda}{2\mu+\lambda}\tr\left(g^{-\frac{1}{2}}(\nabla' \nabla' y_m) g^{-\frac{1}{2}}\right)^2\right){\rm d}\vx',
\end{equation*}
where $y_m$ is the $m^{th}$ component of $\vy$.
In addition, in view of the metric constraint \eqref{metric_constraint} we have $\vy \in [W^{1,\infty}(\Omega)]^3$. As a consequence, the equilibrium deformations of the thin limit of prestrained materials are deformations $\vy \in [H^2(\Omega) \cap W^{1,\infty}(\Omega)]^3$  satisfying $\I(\vy) = g$ and minimizing
\begin{equation}\label{E:prestrain-bending-D2y}
E^{\rm pre}(\vy):= \frac{1}{24}\sum_{m=1}^3\int_{\Omega}\left(2\mu\left|g^{-\frac{1}{2}}(\nabla' \nabla' y_m) g^{-\frac{1}{2}}\right|^2 + \frac{2\mu\lambda}{2\mu+\lambda}\tr\left(g^{-\frac{1}{2}}(\nabla' \nabla' y_m) g^{-\frac{1}{2}}\right)^2\right){\rm d}\vx'.
\end{equation}

\subsection{Bilayer Plates}\label{sec:bilayer}

We now discuss an anisotropic model where the three dimensional hyperelastic structure is a compound of two thin layers with different mechanical properties. Each material exhibits a different response under external stimuli, generating a material mismatch compensated by the bending of the structure. Temperature, humidity, pH, and electric current are typical external stimuli. 
We refer for instance to \cite{SIPL1993,JSI2000,ABS2011,BNS2015} for devices actuated by electric current and to \cite{Menges2015,WVMR2018} for humidity controlled materials.

From a mathematical point of view, this corresponds to a prestrain tensor of the form
\begin{equation} \label{def:Metric_Bilayer}
	G(\vx',x_3) = (I_3 \pm \zeta s  N)^T(I_3 \pm  \zeta sN) = I_3 \pm 2 \zeta s N + \zeta^2s^2 N^2, \qquad  \pm x_3 \geq 0,    
\end{equation}
where the inhomogeneity and anisotropy of the bilayer material is encoded in the tensor 
$$
N := 
\begin{pmatrix}
	N_{\vx'\vx'} & {\mathbf m} \\
	{\mathbf m}^T & n
\end{pmatrix}
$$
with $N_{\vx' \vx'}: \Omega \rightarrow \mathbb{R}^{2\times 2}$ uniformly symmetric, $n: \Omega \rightarrow \mathbb R$, and ${\mathbf m} : \Omega \rightarrow \mathbb R^2$.
Note that the expression \eqref{def:Metric_Bilayer} does not reduce to the previous setting with target metric \eqref{def:G_cstr_3D}.

The limiting model is derived in \cite{S2007} via $\Gamma$-convergence when $N=I_3$ and can again be formally recovered using a modified Kirchhoff-Love assumption, namely 
\begin{equation*}
	\vu(\vx',x_3)=\vy(\vx')+x_3(1\pm c\zeta s) \vnu(\vx')+\frac{1}{2}x_3^2d(\vx')\vnu(\vx')
\end{equation*}
with
$$c:=\frac{2\mu+3\lambda}{2\mu+\lambda} \quad \mbox{and} \quad d(\vx'):=-\frac{\lambda}{2\mu+\lambda}\tr\left(\II(\vx')\right).
$$
Proceeding as in the previous section, the limiting energy as the thickness $s$ goes to zero is
$$\frac{1}{24}\int_{\Omega}\left(2\mu\left|\II(\vy)-3\zeta N_{\vx'\vx'}\right|^2+\frac{2\mu\lambda}{2\mu+\lambda}\tr\left(\II(\vy)-3\zeta N_{\vx'\vx'}\right)^2\right){\rm d}\vx'+C$$
provided that $\vy$ is an isometry, where the constant $C$ depends on $\mu$, $\lambda$, $\zeta$, and $N$, but not on the deformation $\vy$. We refer to \cite{BBN2017} for more details in the case $\mu=6$ and $\lambda=0$.

In what follows, we consider an equivalent limiting two dimensional energy
$$
\frac{\alpha}{2}\int_{\Omega}\left|\II(\vy)-Z\right|^2{\rm d}\vx'
$$
for $\alpha>0$ and with $Z(\vx') := 3 \zeta N_{\vx'\vx'}$. Note that $Z$ acts as an intrinsic spontaneous curvature tensor favoring deformations with principal curvatures matching the eigenvalues of $Z$. Now recall that the entries of the second fundamental form $\II(\vy)$ are given by
$$
\left[\II(\vy)\right]_{ij} = \partial_{ij}\vy\cdot \vnu = \partial_{ij}\vy\cdot\frac{\partial_1\vy\times\partial_2\vy}{|\partial_1\vy\times\partial_2\vy|}, \quad i,j=1,2,
$$
and is thus nonlinear in $\vy$. 
Taking advantage of relation \eqref{def:II_D2y} valid for isometries, we have
$$
\left| \II(\vy) \right| = \left| \nabla' \nabla' \vy\right| \qquad \textrm{and} \qquad \left| \partial_1 \vy \times \partial_2 \vy \right|  = 1,
$$
which allows us to write
\begin{equation*}
\frac{\alpha}{2}\int_{\Omega}\left|\II(\vy)-Z\right|^2{\rm d}\vx' = \frac{\alpha}{2}\int_{\Omega}\left|\nabla'\nabla'\vy\right|^2{\rm d}\vx'-\alpha \sum_{i,j=1}^2 \int_\Omega \partial_{ij}\vy\cdot \left(\partial_1\vy\times\partial_2\vy\right) Z_{ij}{\rm d}\vx'+\frac{\alpha}{2}\int_{\Omega}\left|Z\right|^2{\rm d}\vx'.
\end{equation*}
Deformations of bilayer materials are thus deformations $\vy  \in [H^2(\Omega) \cap W^{1,\infty}(\Omega)]^3$ satisfying $\I(\vy) = I_2$ and minimizing
\begin{equation} \label{def:E_2D_bilayer}
	E^{\rm bil}(\vy):= \frac{\alpha}{2}\int_{\Omega}\left|\nabla'\nabla'\vy\right|^2{\rm d}\vx'-\alpha \sum_{i,j=1}^2 \int_\Omega\partial_{ij}\vy\cdot \left(\partial_1\vy\times\partial_2\vy\right)Z_{ij}{\rm d}\vx'.
\end{equation}

\subsection{Folding} \label{sec:folding}

We incorporate the capability for the plates to fold along a given crease.
In particular, we make structural assumptions leading to a limiting model for which folding does not require any energy. We restrict our description to single layer plates with isometry constraint (i.e. $g=I_2$). The derivation of the folding model described here was originally proposed in \cite{BBH2021}.
Extensions to other models, although feasible, have not been treated in the existing literature. 
However, numerical simulation are available in \cite{BNY2022} for the bending of a bilayer plate with folding. 

Folding models are obtained by assuming that the thin structure has a material defect in a cylindrical neighborhood of a $C^2$ curve $\Sigma \subset \Omega$ (the crease) splitting $\Omega$ in two parts and intersecting the boundary transversely. Let $\Omega_s=\Omega\times(-s/2,s/2)$ be as above and let $\Sigma_{s,r} := B(\Sigma,r) \times (-s/2,s/2)$ be the location of the defect. Here, for $r>0$ to be determined, $B(\Sigma,r) := \cup_{\vx' \in \Sigma} \{ \vz' \in \Omega \ : \ | \vz' - \vx' | < r \}$ (see Figure~\ref{fig:folding_diagram}).

\begin{figure}[htbp]
	\hspace*{-0.2cm}\includegraphics[width=10cm]{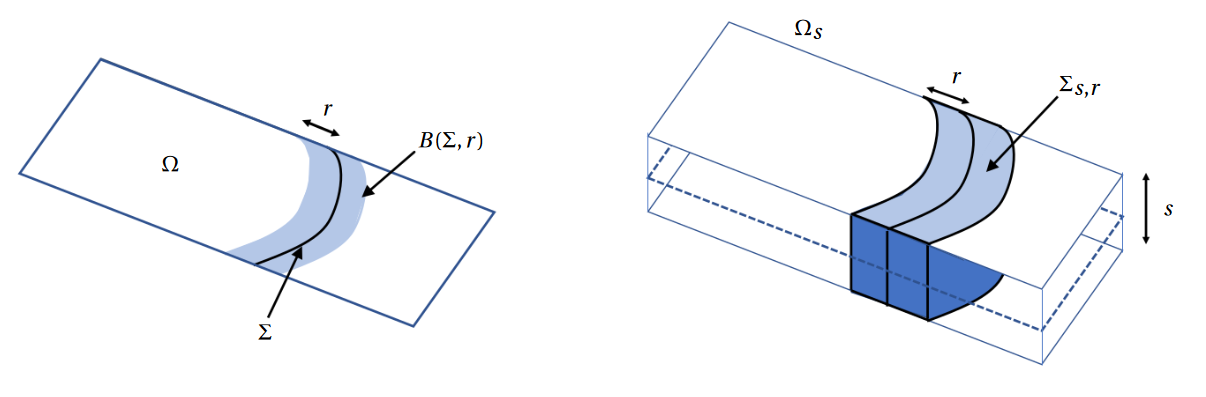}
	\caption{Material with defect. (Left) Midplane domain $\Omega$ with a crease $\Sigma$ and its tubular neighborhood $B(\Sigma,r)$. (Right) Thin domain $\Omega_s$ with inflated crease $\Sigma_{s,r}$.} \label{fig:folding_diagram}
\end{figure}

The hyperelastic energy \eqref{def:E_3D} in presence of a defect reads
\begin{equation*}
	E^{\rm fold}_s(\vu):=\int_{\Omega_s}f_{\varepsilon,r}(\vx')W(\nabla\vu(\vx)){\rm d}\vx,
\end{equation*}
where $f_{\varepsilon,r}:\Omega\rightarrow(0,1]$ defined by
\begin{equation*}
	f_{\varepsilon,r}(\vx'):=\varepsilon\chi_{B(\Sigma,r)}(\vx')+1-\chi_{B(\Sigma,r)}(\vx')
\end{equation*}
is $1$ away from $B(\Sigma,r)$ and $\varepsilon$ in $B(\Sigma,r)$ for $0<\varepsilon\ll 1$. Here $\chi_{B(\Sigma,r)}$ is the characteristic function of $B(\Sigma,r)$.

We assume the following asymptotic behaviors of $r$ and $\varepsilon$:
\begin{equation} \label{Hyp:folding_param}
	\limsup_{s\rightarrow 0^+}\frac{s^2}{\varepsilon}<\infty, \quad \limsup_{s\rightarrow 0^+}\frac{s}{r}<\infty, \quad \mbox{and} \quad \limsup_{s\rightarrow 0^+}\frac{\varepsilon r}{s^2}=0.
\end{equation}
They guarantee that the limiting deformations are globally Lipschitz and that folding is performed at no energy cost.
The $\Gamma$-limit of $s^{-3} E^{\rm fold}_s$ under assumptions \eqref{Hyp:folding_param} and in the bending regime $s^{-3} E^{\rm fold}_s \leq \Lambda < \infty$, is
\begin{equation*}
	E^{\rm fold}(\vy) := \frac{1}{24}\int_{\Omega\setminus\Sigma}\left(2\mu\left|\II(\vy)\right|^2+\frac{2\mu\lambda}{2\mu+\lambda}\tr\left(\II(\vy)\right)^2\right){\rm d}\vx'= \frac{\alpha}{2}\int_{\Omega\setminus\Sigma}\left|\nabla' \nabla' \vy\right|^2{\rm d}\vx',
\end{equation*}
provided $\vy\in[H^2(\Omega\setminus\Sigma)\cap W^{1,\infty}(\Omega)]^3$ is an isometry and where
$\alpha=\alpha(\mu,\lambda)$ is defined as in \eqref{def:EB_iso_simple}; see \cite{BBH2021}.
In particular, $E^{\rm fold}$ reduces to $E$ in \eqref{def:EB_iso_simple} when no creases are present. 

Although not rigorously justified, in the numerical section below we provide experiments with multiple folding arcs not necessarily intersecting the boundary.
Also, we assume that the prestrain and bilayer energies can be modified similarly to account for possible folding along a crease $\Sigma$. Without modifying the notation but with the understanding that $\Sigma = \emptyset$ in the original cases \eqref{E:prestrain-bending-D2y} and \eqref{def:E_2D_bilayer}, we redefine for $\vy \in [H^2(\Omega \setminus \Sigma) \cap W^{1,\infty}(\Omega)]^3$ 
\begin{equation}\label{E:prestrain-bending-D2y-fold}
E^{\rm pre}(\vy) := \frac{1}{24}\sum_{m=1}^3\int_{\Omega \setminus \Sigma}\Big(2\mu\left|g^{-\frac{1}{2}}(\nabla' \nabla' y_m) g^{-\frac{1}{2}}\right|^2 + \frac{2\mu\lambda}{2\mu+\lambda}\tr\left(g^{-\frac{1}{2}}(\nabla' \nabla' y_m) g^{-\frac{1}{2}}\right)^2\Big){\rm d}\vx'
\end{equation}
and
\begin{equation} \label{def:E_2D_bilayer-fold}
	E^{\rm bil}(\vy):= \frac{\alpha}{2}\int_{\Omega\setminus \Sigma}\left|\nabla'\nabla'\vy\right|^2{\rm d}\vx'-\alpha \sum_{i,j=1}^2 \int_{\Omega \setminus \Sigma}\partial_{ij}\vy\cdot \left(\partial_1\vy\times\partial_2\vy\right)Z_{ij}{\rm d}\vx'.
\end{equation}

\subsection{Forcing Term and Boundary Conditions} \label{sec:force_BC}

To complete the model, it remains to discuss how to incorporate external forces and boundary conditions. 

In the presence of an external force $\vf_s:\Omega_s \rightarrow \mathbb R^3$, the hyperelastic energy \eqref{def:E_3D} must be modified to read
$$
E_{s}(\vu) = \int_{\Omega_s}W(\nabla \vu(\vx)){\rm d}\vx - \int_{\Omega_s} \vf_s (\vx)\cdot\vu(\vx){\rm d}\vx.
$$
In view of \eqref{e:bending_regim}, the force applied to the three dimensional thin structure $\Omega_s = \Omega \times (-s/2,s/2)$ must satisfy
\begin{equation*}
	\left| \lim_{s\rightarrow 0^+}s^{-3}\int_{-s/2}^{s/2}\vf_s(\vx',x_3){\rm d}x_3 \right| < \infty, \quad \vx'\in\Omega.
\end{equation*}
This allows us to define $\vf(\vx'):=\lim_{s\rightarrow 0^+}s^{-3}\int_{-s/2}^{s/2}\vf_s(\vx',x_3){\rm d}x_3$ which yields
\begin{equation*}
	\lim_{s\rightarrow 0^+}s^{-3}\int_{\Omega_s}\vf_s(\vx)\cdot\vu(\vx){\rm d}\vx = \int_{\Omega}\vf(\vx')\cdot\vy(\vx'){\rm d}\vx'
\end{equation*}
when using the modified Kirchhoff-Love ansatz \eqref{Hyp:kirchhoff_2} on the three dimensional deformation $\vu$.

Therefore, when external forces are present, the term
\begin{equation}\label{e:rhs}
	F(\vy):=\int_{\Omega}\vf(\vx')\cdot\vy(\vx'){\rm d}\vx'
\end{equation}
must be subtracted from the energies $E^{\rm pre}_s(\vy)$, $E^{\rm pre}(\vy)$, and $E^{\rm bil}(\vy)$ (scaled by $s^2$ for the preasymptotic energy).

Different types of boundary conditions can be incorporated into the system. 
We say that \emph{Dirichlet boundary conditions} are imposed on $\Gamma^D \subset \partial \Omega$ when 
\begin{equation} \label{eqn:Dirichlet_BC}
	\vy = \vvarphi \quad \mbox{and} \quad \nabla'\vy = \Phi \quad \mbox{on } \Gamma^D,
\end{equation}
where $\vvarphi:\Omega\rightarrow\mathbb{R}^3$ and $\Phi:\Omega\rightarrow\mathbb{R}^{3\times 2}$ are sufficiently smooth, and with $\Phi$ satisfying the compatibility condition $\Phi^T\Phi=g$ a.e. in $\Omega$. In the case where $\vvarphi(\vx)=(\vx,0)^T$ and $\Phi=\nabla\vvarphi$, the boundary conditions \eqref{eqn:Dirichlet_BC} are referred to as \emph{clamped boundary conditions}. \emph{Mixed boundary condition} occurs when only the value of the deformation on $\Gamma^M \subset \partial \Omega$ is specified, that is
\begin{equation} \label{eqn:Mixed_BC}
	\vy = \vvarphi \quad \mbox{on } \Gamma^M
\end{equation}
for some smooth function $\vvarphi:\Omega\rightarrow\mathbb{R}^3$. 

All the energy models are defined on $[W^{1,\infty}(\Omega)]^3$ and thus also on continuous functions $[C^0(\overline{\Omega})]^3$. This means that \emph{pointwise boundary conditions} can be enforced on the deformations as well. Given a collection of points $\Gamma^P:=\{\vx_1,\vx_2,\ldots,\vx_P\}$ lying on $\partial\Omega$ and  corresponding values $\{\vvarphi_1,\vvarphi_2,\ldots,\vvarphi_P\}\subset\mathbb{R}^3$, we set
\begin{equation*}
	\vy(\vx_i) = \vvarphi_i \quad \forall\, \vx_i\in\Gamma^P. 
\end{equation*}

The \emph{free boundary case} refers to the case where $\Gamma^D\cup\Gamma^M\cup\Gamma^P=\emptyset$.

\section{Space Discretization} \label{sec:method}

We now describe numerical methods based on finite element methods to approximate minimizers of the preasymptotic energy $E_s^{\rm pre}$ in \eqref{eqn:ES_EB_G} and the limiting energies $E^{\rm pre}$ in \eqref{E:prestrain-bending-D2y-fold} and $E^{\rm bil}$ in \eqref{def:E_2D_bilayer-fold}. We point out two difficulties addressed below: (i) all the energies are defined on a subspace of $[H^2(\Omega \setminus \Sigma)]^3$ (recall that $\Sigma=\emptyset$ when no creases are present) and thus conforming discretizations can only be achieved by globally $C^0$ elements which are also $C^1$ on $\Omega \setminus \Sigma$; (ii) all but the preasymptotic energy require the deformations to satisfy the metric constraint \eqref{metric_constraint}. The latter is too rigid to be satisfied a.e. in $\Omega$ when using piecewise polynomial approximations. 

Discontinuous Galerkin finite elements are put forward to circumvent the conformity requirement. 
In essence, discrete energies are designed by substituting the Hessian of the deformation by a suitable discrete Hessian. It is worth pointing out that discontinuous Galerkin elements are also preferred for their flexibility; in particular, they are more efficient for the approximation of the metric constraint \cite{BNN2020,BNN2021}. 

We first set up the notation related to discontinuous Galerkin finite element methods in Section~\ref{sec:DG}.
Then, in Section~\ref{sec:H_h}, we introduce a discrete Hessian operator, the crucial ingredient in all the numerical schemes discussed in Section~\ref{sec:E_h}.

To simplify the notation, from now on we will write $\vx$ instead of $\vx'$, and accordingly for the differential operators (e.g. $\nabla'$ will be denoted $\nabla$). Moreover, the variable of integration will no longer be indicated. 

\subsection{Discontinuous Galerkin Finite Elements} \label{sec:DG}

For the description of the method we assume that the midplane $\Omega\subset\mathbb{R}^2$ is a polygonal domain, but extensions to general domains can be performed using standard techniques. Let $\{\Th\}_{h>0}$ be a sequence of shape-regular but possibly graded partitions of $\Omega$ made of elements $\K$, either triangles or quadrilaterals, of diameter at most $h$. Hanging nodes are thus allowed, and we assume that all the elements within each domain of influence have comparable diameters \cite{bonito2010quasi}. 
Let $\Eh:=\Eh^0\cup\Eh^b$ denote the set of edges, where $\Eh^0$ stands for the set of interior edges and $\Eh^{b}$ for the set of boundary edges. We also define $\Eh^D$ to be the set of Dirichlet boundary edges and $\Eh^M$ to be the set of mixed boundary edges.

In the presence of a crease $\Sigma$, we assume that the crease is well approximated by the subdivision $\Th$ in the spirit of \cite{BBH2021}. That is, the folding curve $\Sigma$ is approximated by the piecewise linear curve
\begin{equation} \label{def:sigma_h}
	\Sigma_h := \cup_{l=1}^Le_l,
\end{equation}
where $\mathcal E_h^\Sigma:= \{ e_l\}_{l=1}^L \subset \mathcal E_h^0$ is a collection of interior edges such that the endpoints of each edge $e_l$ belong to $\Sigma$.

It will be convenient later to group all the \emph{active edges} into two sets: $\Ehv:=\Eh^0  \cup\Eh^D\cup\Eh^M$ for the deformations  and $\Ehg:=(\Eh^0\setminus \E_h^\Sigma)\cup\Eh^D$ for the gradient of the deformations. These sets contain all the edges over which averages and jumps will be computed. Notice that in the presence of a crease, the crease edges are not included in the active edge set for the gradients to allow for discontinuous gradients. As we shall see, this is the only difference between the schemes with and without folding.

We define the diameter function $\h$ on $\Th \cup \E_h \cup \Gamma^P$ (recall that $\Gamma^P=\{\vx_i\}_{i=1}^P$ are the locations where pointwise boundary conditions are prescribed)
by
\begin{equation*}
	\begin{split}
		\restriction{\h}{\K}&:= h_{\K} := \textrm{diam}(\K), \quad \forall\, \K\in\Th, \quad
		\restriction{\h}{e}:= h_e:=\textrm{diam}(e), \quad \forall\, e\in\Eh, \\
		\restriction{\h}{\vx_i}&:= h_i := \frac{1}{\#(\omega_i)}\sum_{\K\in\omega_i}h_{\K}, \quad i=1,2,\ldots,P, 
	\end{split}
\end{equation*}
where $\omega_i:=\left\{\K\in\Th:\,\, \vx_i\in\partial\K\right\}$. The minimum mesh size is $h_{\min}:=\min_{\K\in\Th}h_T$.

We assume that the boundary regions $\Gamma^D$ and $\Gamma^M$, where Dirichlet and mixed boundary conditions are enforced, are exactly captured by all the subdivisions.
This means that for all $h>0$ we have
\begin{equation*}
	\Gamma^D=\{e: \,\, e\in\Eh^D\} \quad \mbox{and} \quad \Gamma^M=\{e: \,\, e\in\Eh^M\}    
\end{equation*}
for some $\Eh^D,\Eh^M\subset\Eh^b$. 
We also assume that the locations $\Gamma^P= \{\vx_i\}_{i=1}^P$, where pointwise boundary conditions are prescribed, are vertices of all the subdivisions.
To keep a uniform notation, we introduce various skeletons of the mesh $\Th$
$$
\Gamma_h^{a}=\cup\{e: e\in\Eh^a\} \quad \mbox{for } a\in\{0,D,M,\rm val,\rm grad\}.
$$

For an integer $r\ge 0$, let $\mathbb{P}_r$ and $\mathbb{Q}_r$ denote polynomials of total degree at most $r$ and degree at most $r$ in each variable, respectively. We then define the finite element space
\begin{equation} \label{def:Vhr}
	\V_h^r:=\left\{v_h\in L^2(\Omega): \,\, \restriction{v_h}{\K}\circ F_{\K}\in\mathbb{P}_r \quad \forall\, \K \in\Th \right\},
\end{equation}
where $\mathbb{P}_r$ is replaced by $\mathbb{Q}_r$ if the subdivision is made of quadrilaterals. Here $F_{\K}:\widehat \K\rightarrow \K$ is the map from the reference element $\widehat \K$ (unit simplex for triangulation and unit square for subdivisions made of quadrilaterals) to the physical element $T\in\Th$. 
In what follows, the deformation of the plate $\vy$ will be approximated by a discrete function $\vy_h\in[\V_h^k]^3$ with $k\ge 2$.

The broken (elementwise) gradient of a scalar function $v_h\in\V_h^k$ is denoted $\nabla_h v_h$. We use a similar notation for other  differential operators, for instance $D_h^2 v_h=\nabla_h\nabla_h v_h$ denotes the broken Hessian. Moreover, these operators are applied componentwise for vector-valued functions.

Let us now introduce the jump and average operators instrumental for discontinuous Galerkin methods to relax conformity conditions. For each $e\in\Eh^0$ let $\vn_e$ be a unit normal to $e$ (the orientation is arbitrary but fixed once for all). With this notation, we set for $v_h \in \V_h^r$ and $e \in \E_h^0$
\begin{equation*}
	\restriction{\jump{v_h}}{e} := v_h^{-}-v_h^+, \qquad \restriction{\avrg{v_h}}{e} := \frac{v_h^{+}+v_h^{-}}{2},
\end{equation*}
where $v_h^{\pm}(\vx):=\lim_{\varrho\rightarrow 0^+}v_h(\vx\pm \varrho\vn_e)$ for any $\vx \in e$. Note that if $e=\partial\K_1\cap\partial\K_2$ with $\K_1,\K_2\in\Th$, then
\begin{equation*}
	\restriction{\jump{v_h}}{e}\vn_e = \restriction{v_h}{\K_1}\vn_{\partial\K_1}+\restriction{v_h}{\K_2}\vn_{\partial\K_2}
\end{equation*}
with $\vn_{\partial\K_i}$ the outward unit normal to $\partial\K_i$, $i=1,2$.
For edges on the boundary $e\in\Eh^b$, $\vn_e$ is the outward unit normal vector to $\partial\Omega$ and $\restriction{\jump{v_h}}{e}=\restriction{\avrg{v_h}}{e}:=v_h$. 
The jumps and averages of non-scalar functions are computed componentwise.

We introduce the mesh-dependent bilinear form $\langle\cdot,\cdot\rangle_{H_h^2(\Omega)}$ defined for any $v_h, w_h\in \V_h^k$ by
\begin{eqnarray} \label{def:H2_inner_product_A}
	\langle v_h,w_h\rangle_{H_h^2(\Omega)} & := &  (D_h^2 v_h,D_h^2 w_h)_{L^2(\Omega)} \nonumber \\
	& & + (\h^{-1}\jump{\nabla_h v_h},\jump{\nabla_h w_h} )_{L^2(\Ghg)} + (\h^{-3}\jump{ v_h},\jump{ w_h} )_{L^2(\Ghv)} \nonumber \\
	& & +\sum_{\vx_i \in \Gamma^P} h_i^{-2}v_h(\vx_i)w_h(\vx_i)
\end{eqnarray}
and we set 
$$
| v_h |_{H^2_h(\Omega)}:= \langle v_h,v_h\rangle_{H^2_h(\Omega)}^{1/2}.
$$
These definitions are extended componentwise to vector-valued functions. 

\subsection{Discrete Hessian}\label{sec:H_h}

All the energies considered in this work involve the Hessian of the deformation. 
Conforming approximations come at the price of using costly and rigid $C^1$ finite element methods. 
Instead, we describe an LDG approach, which retains the simplicity of the discrete problem while not requiring $C^1$ elements.  

The proposed LDG methods consist of replacing the continuous Hessian $D^2\vy$ by a suitable discrete (aka reconstructed) Hessian $H_h(\vy_h)$ that we introduce now. The discrete Hessian described here was originally proposed in \cite{BGNY2022a} following ideas from, for example, \cite{DPE2011,DPE2010} but modified to possibly account for subdivision containing non-affine elements.

For $\vv_h \in [\V_h^k]^3$, we recall that $D^2_h \vv_h: \Omega \rightarrow \mathbb{R}^{3\times 2\times 2}$ is the piecewise Hessian of $\vv_h$. Because $\vv_h$ does not belong to $[H^2(\Omega)]^3$, a better approximation of the Hessian must account for the jumps of $\vv_h$ and $\nabla_h \vv_h$ across elements. To achieve this, we use (local) lifting operators \cite{bassi1997high,arnold2002unified,brezzi2000discontinuous} $r_e:[L^2(e)]^2\rightarrow[\V_h^{l_1}]^{2\times 2}$ and $b_e:L^2(e)\rightarrow[\V_h^{l_2}]^{2\times 2}$ to construct representations of these jumps in $L^2(\Omega)$.
Given $l_1$, $l_2 \geq 0$, they are defined for $e\in\Eh$ by
\begin{equation*}
	r_e(\vphi) \in [\V_h^{l_1}]^{2\times 2}: \,
	\int_{\Omega}r_e(\vphi):\tau_h = \int_e\avrg{\tau_h}\vn_e\cdot\vphi \quad \forall\, \tau_h\in [\V_h^{l_1}]^{2\times 2},
\end{equation*}
and
\begin{equation*}
	b_e(\phi) \in [\V_h^{l_2}]^{2\times 2}: \,
	\int_{\Omega} b_e(\phi):\tau_h = \int_e\avrg{\di \tau_h}\cdot\vn_e\phi \quad \forall\, \tau_h\in [\V_h^{l_2}]^{2\times 2}.
\end{equation*} 
Note that the support of the liftings $r_e$ and $b_e$ is the union of the two elements sharing $e$ as an edge (which reduces to a single element if $e$ is a boundary edge).

With these lifting operators, we can now define the discrete Hessian of a function $v_h\in\V_h^k$ as
\begin{eqnarray} \label{def:H_h}
	H_h(v_h) & := & D_h^2(v_h)-R_h(\jump{\nabla_h v_h})+B_h(\jump{v_h}) \nonumber \\
	& := & D_h^2(v_h)-\sum_{e\in\Ehg}r_e(\jump{\nabla_h v_h})+\sum_{e\in\Ehv}b_e(\jump{v_h}).
\end{eqnarray} 
This definition is naturally extended to vector-valued functions in $[\V_h^k]^3$ by componentwise application.
For the remainder of this work, we set $l_1=l_2=k$.

To motivate definition \eqref{def:H_h} of the discrete Hessian, we briefly sketch an argument showing a (weak) convergence property of $H^2_h(v_h)$. 
Let $\tau\in [C_0^{\infty}(\Omega)]^{2\times 2}$ and let $\{v_h\}_{h>0}\subset \V_h^k$ be a sequence such that $| v_h |_{H^2_h(\Omega)} \leq C$ for some constant $C$ independent of $h$ and such that $v_h\rightarrow v$ in $L^2(\Omega)$ as $h\rightarrow 0^+$ for some $v\in H^2(\Omega)$.
On the one hand, notice that thanks to the regularity of $v$ and the convergence of $v_h$ to $v$, using two integrations by parts we have
\begin{equation*}
	\int_{\Omega}D^2v:\tau  = \int_{\Omega}v\di(\di(\tau)) \leftarrow \int_{\Omega}v_h\di(\di(\tau))  \quad \mbox{as } h\rightarrow 0^+.
\end{equation*}
On the other hand,  using again two integrations by parts, we have for any $\K\in\Th$
\begin{eqnarray*}
	\int_{\K}v_h\di(\di(\tau)) & = &  \sum_{i,j=1}^2\int_{\K}v_h\partial_i(\partial_j \tau_{ij}) \\
	& = &  \sum_{i,j=1}^2\left[\int_{\K}\partial_j(\partial_iv_h)\tau_{ij}-\int_{\partial\K}\partial_iv_h\tau_{ij}(\vn_{\partial\K})_j+\int_{\partial\K}v_h\partial_j\tau_{ij}(\vn_{\partial\K})_i\right] \\
	& = & \int_{\K}D^2v_h:\tau-\int_{\partial\K}\nabla v_h\cdot(\tau\vn_{\partial\K})+\int_{\partial\K}v_h\di(\tau)\cdot\vn_{\partial\K}. 
\end{eqnarray*}
Therefore, summing over the elements $\K\in\Th$ and using the fact that $\tau$ is smooth, namely that $\jump{\tau}=\avrg{\tau}=\tau$ and similarly for $\di(\tau)$, we get
\begin{equation*} \label{eqn:motivation_Hh}
	\int_{\Omega}v_h\di(\di(\tau))=\int_{\Omega}D_h^2v_h:\tau-\sum_{e\in\Eh^{\rm grad}}\int_{e}\jump{\nabla_h v_h}\cdot\avrg{\tau}\vn_e+\sum_{e\in\Eh^{\rm val}}\int_{e}\jump{v_h}\avrg{\di(\tau)}\cdot\vn_e.
\end{equation*}

If $\tau$ were an admissible test function for the lifting operators, then the right hand side of the last relation would be $\int_{\Omega}H_h(v_h):\tau$ and we would conclude that for all $\tau \in [C_0^{\infty}(\Omega)]^{2\times 2}$
$$
\int_\Omega D^2 v : \tau \leftarrow \int_\Omega H_h(v_h): \tau \qquad \textrm{as }h \to 0^+,
$$
i.e. $H_h(v_h)$ converges weakly to $D^2 v$ in $L^2(\Omega)^{2\times 2}$. Because $\tau$ does not belong to $[\V_h^k]^{2\times 2}$ in general, it remains to show that for the Lagrange interpolant $I_h:[C_0^{\infty}(\Omega)]^{2\times 2} \rightarrow [\V_h^k \cap H^1_0(\Omega)]^{2\times 2}$ we have
\begin{equation*}
	\int_{\Omega}D_h^2v_h:(\tau-I_h\tau)-\sum_{e\in\Eh^{\rm grad}}\int_{e}\jump{\nabla_h v_h}\cdot\avrg{\tau-I_h \tau}\vn_e +\sum_{e\in\Eh^{\rm val}}\int_{e}\jump{v_h}\avrg{\di(\tau - I_h \tau)}\cdot\vn_e \rightarrow 0
\end{equation*}
and
$$
\int_\Omega H_h(v_h): (\tau - I_h \tau) \rightarrow 0
$$
as $h \to 0^+$.
This fact follows from the uniform boundedness assumptions $|v_h|_{H^2_h(\Omega)}\leq C$ and standard interpolation estimates. We refer to Lemma~2.4 in \cite{BGNY2022b} for additional details when $\Sigma = \emptyset$ (see Appendix B for a proof in the \emph{free boundary case} and Appendix C for the Dirichlet boundary conditions case). The case $\Sigma \not = \emptyset$ is derived similarly. 
Strong convergence properties are also available, for example, in \cite{BGNY2022b}.

We now turn our attention to the boundary conditions and start by noting that the Dirichlet/mixed boundary edges $e\in\Eh^D\cup\Eh^M$ are incorporated in the definition of the discrete Hessian.
They do not influence the weak convergence property sketched above (because $\tau$ is compactly supported).
However, the value of the deformation and its gradient imposed on $\Gamma^D \cup \Gamma^M$ must be added separately for strong consistency \cite{nitsche1971variationsprinzip}. 
To this end, we introduce the lifting of the boundary conditions \eqref{eqn:Dirichlet_BC} and \eqref{eqn:Mixed_BC}
\begin{equation*}
	L_h(\Phi_m,\varphi_m) := -\sum_{e\in\Eh^D}r_e(\Phi_m)+\sum_{e\in\Eh^D\cup\Eh^M}b_e(\varphi_m),
\end{equation*}
where we recall that $\Phi_m$ stands for the $m^{th}$ row of $\Phi$, $m=1,2,3$.
The discrete Hessian operator accounting for the boundary condition is then chosen to be
\begin{equation} \label{def:H_h_with_BC}
	H_h(\varphi_m,\Phi_m; y_{h,m}) := H_h(y_{h,m}) - L_h(\Phi_m,\varphi_m)
\end{equation}
with the usual extension for vector-valued functions.

\subsection{Discrete Energies}\label{sec:E_h}

We are now in a position to describe the approximation of the different energies $E_s^{\rm pre}(\vy)$ (preasymptotic), $E^{\rm pre}(\vy)$ (prestrain), and $E^{\rm bil}(\vy)$ (bilayer). 

\subsubsection{Preasymptotic Energy}
We start with the preasymptotic energy 
\begin{equation*}
	E_s^{\rm pre}(\vy) = E^S(\vy)+s^2E^B(\vy)
\end{equation*}
and approximate the stretching and bending components separately. 
For the stretching part, we simply replace the derivatives appearing in \eqref{def:ES} by their piecewiese counterparts and define for $\vy_h \in [\V_h^k]^3$
\begin{equation} \label{def:E_h_stretching}
E_h^S(\vy_h) := \frac{\mu}{4}\int_{\Omega}\left|g^{-\frac{1}{2}}\big((\nabla_h\vy_h)^T\nabla_h\vy_h-g\big)g^{-\frac{1}{2}}\right|^2 +  \frac{\lambda}{8}\int_{\Omega}\tr\Big(g^{-\frac{1}{2}}\big((\nabla_h\vy_h)^T\nabla_h\vy_h-g\big)g^{-\frac{1}{2}}\Big)^2.
\end{equation}
For the bending part, we first replace as in \cite{BK2014} the second fundamental form in $E^B$ by the Hessian. This step is not justified but taken for computational convenience. We then take advantage of the discrete Hessian \eqref{def:H_h_with_BC} and set for $\vy_h \in [\V_h^k]^3$
\begin{equation}\label{def:E_h_bending}
	\begin{split}
		E_{h}^B(\vy_h)  := &\frac{\mu}{12}\sum_{m=1}^3\int_{\Omega}\left|g^{-\frac{1}{2}}H_h(\varphi_m,\Phi_m;y_{h,m})g^{-\frac{1}{2}}\right|^2 \\
		&+\frac{\mu\lambda}{12(2\mu+\lambda)}\sum_{m=1}^3\int_{\Omega}\tr\left(g^{-\frac{1}{2}}H_h(\varphi_m,\Phi_m;y_{h,m})g^{-\frac{1}{2}}\right)^2;
	\end{split}
\end{equation}
compare with \eqref{def:EB}. In the preasymptotic case $\Sigma = \emptyset$ but for later use, we still define the approximation of the bending energy and subsequent quantities using $\Sigma$.

At this point neither $E_h^S$ nor $E_h^B$ enforce weak continuity condition nor are the boundary conditions (when $\Gamma^D \cup \Gamma^M \cup \Gamma^P \not = \emptyset$) imposed properly. This is the purpose of the stabilization term defined for $\vy_h \in [\V_h^k]^3$ by
\begin{equation}\label{def:STAB}
	\begin{split}
		S_h(\vy_h) := &\frac{\gamma_0}{2}\|\h^{-\frac{3}{2}}\jump{\vy_h}\|_{L^2(\Gamma_h^0)}^2 + \frac{\gamma_0}{2}\|\h^{-\frac{3}{2}}(\vy_h-\vvarphi)\|_{L^2(\Gamma_h^D\cup\Gamma_h^M)}^2  \\
		& +\frac{\gamma_1}{2}\|\h^{-\frac{1}{2}}\jump{\nabla_h\vy_h}\|_{L^2(\Gamma_h^0 \setminus \Sigma_h)}^2+ \frac{\gamma_1}{2}\|\h^{-\frac{1}{2}}(\nabla_h\vy_h-\Phi)\|_{L^2(\Gamma_h^D)}^2  \\
		& +\frac{\gamma_2}{2}\sum_{\vx_i \in \Gamma^P} h_i^{-2}|\vy_h(\vx_i)-\vvarphi_i|^2,
	\end{split}
\end{equation}
where $\gamma_0, \gamma_1, \gamma_2 >0$ are stabilization parameters. Note that a particularity of the LDG method is that these parameters are not required to be sufficiently large.   

It remains to incorporate possible external forces - see \eqref{e:rhs} - to arrive at the following approximation of $E^{\rm pre}_s$
\begin{equation*}
	E^{\rm pre}_{s,h}(\vy_h):= E^S_h(\vy_h)+s^2\left( E^B_h(\vy_h) + S_h(\vy_h) - F(\vy_h)\right).
\end{equation*}
The discrete problem is thus to seek $\vy_h^{\star} \in [\V_h^k]^3$ such that 
\begin{equation} \label{def:preasymptotic_min}
	\vy_h^{\star} \in \argmin_{\vy_h \in [\V_h^k]^3 } E^{\rm pre}_{s,h}(\vy_h).
\end{equation}
The $\Gamma$-convergence of $E^{\rm pre}_{s,h}$ to $E^{\rm pre}_{s}$ as $h\rightarrow 0^+$ is proved in \cite{BGMpreprint} when using continuous (rather than discontinuous) piecewise polynomial approximations of the deformations and when $\Sigma=\emptyset$ (i.e. without folding). 

\subsubsection{Prestrain Energy}
The approximation of the prestrain energy $E^{\rm pre}$ in \eqref{E:prestrain-bending-D2y} directly follows from the above discussion since $E^{\rm pre}=E^B$. Accounting for the potential external forces and boundary conditions, we thus set
\begin{equation}\label{E:final-prestrain-bending-discrete}
	E^{\rm pre}_h(\vy_h) := E_h^B(\vy_h) + S_h(\vy_h) - F(\vy_h).
\end{equation}
However, unlike the preasymptotic model, the deformations are required to satisfy the metric constraint \eqref{metric_constraint}.
We have already pointed out that when using polynomial approximations, the metric cannot be satisfied everywhere in $\Omega$. We thus consider the following relaxation of the constraint
\begin{equation} \label{def:Defect}
	\mathcal D_h^{\rm aver}(\vy_h):=\sum_{\K\in\Th}\left|\int_{\K}\left((\nabla\vy_h)^T\nabla\vy_h-g\right)\right| \leq \varepsilon
\end{equation}
for some $\varepsilon>0$ sufficiently small.
This leads to the definition of the discrete admissible sets
\begin{equation} \label{def:Discrete_Admissible_g_1}
	\A_{h,\varepsilon}^{\rm aver} := \left\{\vy_h\in [\V_h^k]^3: \mathcal D_h^{\rm aver}(\vy_h)\le\varepsilon\right\}.
\end{equation}

With these notation, the discrete preastrained minimization problem originally proposed in \cite{BGNY2022a,BGNY2022b} is to find $\vy_h^{\star} \in [\V_h^k]^3$ such that
\begin{equation} \label{def:prestrain_constr_min}
	\vy_h^{\star} \in \argmin_{\vy_h \in \A_{h,\varepsilon}^{\rm aver}} E^{\rm pre}_h(\vy_h).
\end{equation}
The $\Gamma$-convergence of $E^{\rm pre}_h$ to $E^{\rm pre}$ as $h\rightarrow 0^+$ is proved in \cite{BGNY2022b} for the case $\Sigma=\emptyset$ (i.e. without folding). Moreover, the case $\Sigma\ne\emptyset$ is considered in \cite{BBH2021} for the bending problem with isometry constraint (i.e. $g=I_2$) assuming that $\Sigma_h$ in \eqref{def:sigma_h} is an exact representation of the folding curve $\Sigma$. Finally, it is worth mentioning that extension to approximate curves obtained using isoparametric finite elements of degree $k\ge 2$ is investigated \cite{bartels2022error} for the linear bending problem, namely for the bending problem without the isometry constraint which is suitable for small displacements.

\subsubsection{Bilayer Energy}

The bilayer energy \eqref{def:E_2D_bilayer} is composed of two terms: the bending term $E^B$ (with $g=I_2$) approximated by $E^B_h$ and the nonlinear term
$$
N(\vy):= \alpha \sum_{i,j=1}^2 \int_\Omega \partial_{ij}\vy \cdot \left(\partial_1\vy\times\partial_2\vy\right)Z_{ij}.
$$
Note that both terms are fourth order terms, requiring the nonlinearity to be discretized with care. To motivate the proposed scheme, we first note that
\begin{equation}\label{e:nonlinear}
	N(\vy) \lesssim \| \vy \|_{H^2(\Omega)} \| \nabla \vy \|_{L^\infty(\Omega)}^2
\end{equation}
which is bounded for deformations $\vy \in [H^2(\Omega)]^3$ satisfying the isometry constraint $\nabla \vy^T \nabla \vy = I_2$ (and thus $|\nabla \vy| \in L^\infty(\Omega)$).
In particular, the control on the nonlinear term is only possible thanks to the isometry constraint. 

The relaxation \eqref{def:Defect} of the metric constraint does not provide enough control.
Instead, we introduce an alternate discrete admissible set which guarantees a pointwise control at all the barycenter $x_\K$ of the elements $\K$ in the subdivision $\Th$:
\begin{equation} \label{def:Defect2}
	\mathcal D_h^{\rm bary}(\vy_h):=\max_{\K \in \Th} | [(\nabla\vy_h)^T\nabla\vy_h-I_2](x_\K) | \leq \varepsilon
\end{equation}
for $\varepsilon>0$ sufficiently small; compare with \eqref{def:Defect}.
This leads to the definition of the discrete admissible set
\begin{equation} \label{def:Discrete_Admissible_g}
	\A_{h,\varepsilon}^{\rm bary} := \left\{\vy_h\in [\V_h^k]^3: \mathcal D_h^{\rm bary}(\vy_h)\le\varepsilon\right\}.
\end{equation}

The approximate constraint \eqref{def:Defect2} implies (see \cite{BNY2022})
\begin{equation*}
	1-\varepsilon \leq | \partial_i \vy_h(x_\K)| \leq 1+\varepsilon, \quad i=1,2,
\end{equation*}
so that $|\nabla_h \vy_h (x_\K)|$ is uniformly bounded for all $\K \in \Th$.
Whence, in order to obtain a discrete version of \eqref{e:nonlinear}, we approximate the nonlinear term $N(\vy_h)$ using a $1$-point quadrature with quadrature points localized at the barycenter of each element $\K \in \Th$
$$
N(\vy) \approx \alpha \sum_{i,j=1}^2 \sum_{\K \in \Th} |\K| \left[\partial_{ij}\vy \cdot \left(\partial_1\vy\times\partial_2\vy\right)Z_{ij}\right](x_T).
$$

In addition, it would be tempting to replace the second derivatives $\partial_{ij}\vy (x_T)$ by the discrete hessian $H_h(\vy_h) (x_T)$. However, estimates on the quadrature error would require regularity of $\vy$ beyond $[H^2(\Omega)]^3$. To circumvent this issue, we introduce the $L^2(\Omega)$ projection $\overline{H}_{h}(\vy_h)$ of $H_h(\vy_h)$ onto the piecewise constant tensors and set 
$$
N_h(\vy_h) := \alpha \sum_{i,j=1}^2 \sum_{\K \in \Th} |\K| \left[\big(\overline{H}_{h}(\vy_h)\big)_{ij} \cdot \left(\partial_1\vy_h\times\partial_2\vy_h\right) Z_{ij}\right](x_T).
$$

The discrete bilayer minimization problem is to seek $\vy_h^{\star} \in [\V_h^k]^3$ such that
\begin{equation} \label{def:bilayer_constr_min}
	\vy_h^{\star} \in \argmin_{\vy_h \in \A_{h,\varepsilon}^{\rm bary} } E_h^{\rm bil}(\vy_h),
\end{equation}
where for $\vy_h \in [\V_h^k]^3$,
\begin{equation} \label{e:bilayer_discrete}
	E_h^{\rm bil}(\vy_h):= E_h^B(\vy_h) + S_h(\vy_h) - F(\vy_h) - N_h(\vy_h).
\end{equation}
Note that the discrete minimization problem \eqref{def:bilayer_constr_min} was introduced in \cite{BNY2022} following ideas from \cite{BBN2017,BP2021}. Moreover, the proof that $E_h^{\rm bil}$ $\Gamma$-converges to $E^{\rm bil}$ as $h\rightarrow 0^+$ can be found in \cite{BNY2022}.

\section{Energy Minimization: Discrete Gradient Flows} \label{sec:GF}

Now that we have established the discrete minimization problems, it remains to discuss  numerical procedures to construct approximate minimizers to each problem.

We recall that the considered energies are nonconvex and/or the minimization problem is subject to a nonconvex constraint. 
Hence, we use a gradient flow for its robustness and ability to approximately satisfy the type of nonconvex constraints encountered in our context. These benefits come at the expense of slow convergence towards minimizers. 

Gradient flows require an initial condition, and their efficiency depends on both the energy of the initial condition and, for constrained problems, how well the initial condition satisfies the constraint. After introducing the (main) discrete gradient flow, we discuss preprocessing algorithms to construct suitable initial deformations. A bi-harmonic problem is advocated to enforce the boundary conditions (recall that the boundary conditions are part of the energies) while an unconstrained gradient flow is put forward for reducing the metric defect of the initial deformation.

For the unconstrained preasymptotic minimization problem and in the metric preprocessing algorithm, an accelerated (Nesterov-type) algorithm considerably improves the number of iterations required to achieve a near minimum. This is the topic of the last part of this section.

\subsection{Main Discrete Gradient Flows}

Discrete gradient flows are used to minimize the discrete energies introduced in Section~\ref{sec:E_h}. Note that we have two types of minimization problems, namely an unconstrained minimization problem of the form
\begin{equation*}
	\min_{\vy_h\in[\V_h^k]^3} E_h(\vy_h)
\end{equation*}
for the preasymptotic case and a constrained minimization problem of the form
\begin{equation*}
	\min_{\vy_h\in\A_{h,\varepsilon}^i} E_h(\vy_h), \qquad i\in\{{\rm aver},{\rm bary}\},
\end{equation*}
for all the other cases.

Gradient flows are chosen to minimize the different energies because of their robustness (energy decreasing property).
In all the cases, we advocate an $H_h^2(\Omega)$ gradient flow based on the mesh-dependent inner product $(\cdot,\cdot)_{H_h^2(\Omega)}$ defined for any $\vv_h,\vw_h\in[\V_h^k]^3$ by
\begin{equation} \label{def:H2_inner_product}
	(\vv_h,\vw_h)_{H_h^2(\Omega)}  :=  \sigma (\vv_h,\vw_h)_{L^2(\Omega)} + \langle\vv_h,\vw_h\rangle_{H_h^2(\Omega)},
\end{equation}
where $\langle\vv_h,\vw_h\rangle_{H_h^2(\Omega)}$ is given by \eqref{def:H2_inner_product_A}.
Here $\sigma=0$ if $\Gamma^D\ne\emptyset$ and $\sigma=1$ otherwise to ensure that \eqref{def:H2_inner_product} is indeed an inner product on $[\V_h^k]^3$ when no Dirichlet boundary conditions are imposed.
The gradient flow metric is thus
\begin{equation*}
	\|\vv_h\|_{H_h^2(\Omega)} := (\vv_h,\vv_h)_{H_h^2(\Omega)}^{1/2}, \quad \vv_h\in[\V_h^k]^3.
\end{equation*}
We detail below the discrete gradient flows for the different energies.
In order to simplify the discussion, we do not include the possible contribution from external forces (in other words we assume that $\vf=\mathbf{0}$).

\subsubsection{Preasymptotic Energy}

We start with the unconstrained minimization problem \eqref{def:preasymptotic_min} and recall that without external forces
$$
E_{s,h}^{\rm pre}(\vy_h) = E_h^S(\vy_h) + s^2 \left( E_h^B(\vy_h) + S_h(\vy_h) \right).
$$
A minimizing movements procedure is adopted to determine successive approximations reducing $E_{s,h}^{\rm pre}$.
The minimization process requires a pseudo time-step $\tau>0$ and an initial guess $\vy_h^0\in[\V_h^k]^3$.
Ideally, each iteration of the algorithm would read: given $\vy_h^n \in [\V_h^k]^3$, find a deformation minimizing 
\begin{equation*}
	\vv_h \mapsto \left( \frac{1}{2\tau}\|\vv_h - \vy_h^n\|_{H_h^2(\Omega)}^2+E_{s,h}^{\rm pre}(\vv_h)\right)
\end{equation*}
over $[\V_h^k]^3$. Note that the Euler-Lagrange equation associated with the above minimization problem is nonlinear, and solving it would entail using an iterative method.
Instead, we consider the following linearization
$$
\delta E_{s,h}^{\rm pre}(\vy_h^{n+1},\vw_h) \approx \delta E^S_h(\vy_h^n; \vy_h^{n+1},\vw_h) + s^2 \left( \delta E^B_h (\vy_h^{n+1},\vw_h) + \delta S_h(\vy_h^{n+1},\vw_h) \right).
$$
Here $\delta E_h^B(\vv_h,\vw_h)$ is the G\^ateau derivative of $E_h^B$ (see \eqref{def:E_h_bending}) at $\vv_h$ in the direction $\vw_h$, $\delta E_h^S(\vy_h^n;\vv_h,\vw_h)$ is the following linearization at $\vy_h^n$ of the G\^ateau derivative of $E_h^S$ (see \eqref{def:E_h_stretching}) at $\vv_h$ in the direction $\vw_h$
\begin{eqnarray*}
	\delta E_h^S(\vy_h^n;\vv_h,\vw_h) & := & \frac{\mu}{2}\int_{\Omega}g^{-\frac{1}{2}}l(\vv_h,\vw_h)g^{-\frac{1}{2}}:g^{-\frac{1}{2}}m(\vy_h^n)g^{-\frac{1}{2}} \\
	& & + \frac{\lambda}{4}\int_{\Omega}\tr\left(g^{-\frac{1}{2}}l(\vv_h,\vw_h)g^{-\frac{1}{2}}\right)\tr\left(g^{-\frac{1}{2}}m(\vy_h^n)g^{-\frac{1}{2}}\right),
\end{eqnarray*}
where
$$l(\vv_h,\vw_h):=(\nabla_h\vv_h)^T\nabla_h\vw_h+(\nabla_h\vv_h)^T\nabla_h\vw_h \quad \mbox{and} \quad m(\vy_h^n):=(\nabla_h\vy_h^n)^T\nabla_h\vy_h^n-g,$$
and
$\delta S_h(\vv_h,\vw_h)$ is the G\^ateau derivative of $S_h$ (see \eqref{def:STAB}) at $\vv_h$ in the direction $\vw_h$.

With these notations, we define $\vy_h^{n+1} \in [\V_h^k]^3$ to be solution of
\begin{equation} \label{def:incrt_variational_BS}
	\tau^{-1}(\vy_h^{n+1}-\vy_h^n,\vw_h)_{H_h^2(\Omega)} + \delta E_h^S(\vy_h^n;\vy_h^{n+1},\vw_h) +s^2 \left(\delta E_h^B(\vy_h^{n+1},\vw_h)+\delta S_h(\vy_h^{n+1},\vw_h) \right) = 0
\end{equation}
for all $\vw_h\in[\V_h^k]^3$. 
Note that \eqref{def:incrt_variational_BS} has a unique solution and the sequence of energies $\{E_{s,h}^{\rm pre}(\vy_h^n)\}_{n\geq 0}$ is decreasing, i.e.
\begin{equation*}
	\frac{1}{\tau}\|\vy_h^{n+1}-\vy_h^n\|_{H_h^2(\Omega)}^2+E_{s,h}^{\rm pre}(\vy_h^{n+1})\le E_{s,h}^{\rm pre}(\vy_h^n)
\end{equation*}
provided $\tau$ is sufficiently small (proportional to $h_{\min}$ and $E_{s,h}^{\rm pre}(\vy_h^0)^{-1}$).
As a consequence, we obtain $E_{s,h}^{\rm pre}(\vy_h^{n+1})<E_{s,h}^{\rm pre}(\vy_h^n)$ if $\vy_h^{n+1}\neq\vy_h^n$. We refer to \cite{BGMpreprint} for more details, see also \cite{BGNY2022b}. 

\subsubsection{Prestrain Energy}

The prestrain energy  $E_h^{\rm pre} = E_h^B + S_h$ appears simpler to reduce because it does not include the problematic non-quadratic stretching energy $E_h^S$. However, the metric condition appears as a constraint on the deformations, which have to belong to the admissible set $\mathbb A_{h,\varepsilon}^{\rm aver}$ (see \eqref{def:Discrete_Admissible_g_1}) for a prescribed $\varepsilon>0$. 

Note that the metric constraint $\nabla \vy^T \nabla \vy = g$ is relaxed in $\mathbb A_{h,\varepsilon}^{\rm aver}$ and only needs to be satisfied approximately. To achieve this, we compute increments in a (pseudo-)tangent space of the constraint and rely on the $H^2_h(\Omega)$ metric of the gradient flow to guarantee that all the deformations generated by the gradient flow indeed belong to $\mathbb A_{h,\varepsilon}^{\rm aver}$ for a specific choice of $\varepsilon$. 

Given a deformation $\vy_h^n \in [\V_h^k]^3$, we define the (pseudo-)tangent space at $\vy_h^n$ as follows:
\begin{equation*} 
	\W_{h}^{\rm aver}(\vy_h^n) := \left\{\vw_h\in[\V_h^k]^3:\,\, \int_{\K} \left((\nabla_h\vw_h)^T\nabla_h\vy_h^n+(\nabla_h\vy_h^n)^T\nabla_h\vw_h \right)=0 \quad \forall\, \K\in\Th\right\}.
\end{equation*}
Now, given a pseudo time-step $\tau>0$ and an initial guess $\vy_h^0\in[\V_h^k]^3$, the discrete gradient flow algorithm for the constrained minimization problem \eqref{def:prestrain_constr_min}  consists in computing successively $\vy_h^{n+1} \in [\V_h^k]^3$ such that $\vy_h^{n+1}-\vy_h^n \in \W_{h}^{\rm aver}(\vy_h^n)$ and
\begin{equation*}
	\tau^{-1}(\vy_h^{n+1}-\vy_h^n,\vw_h)_{H_h^2(\Omega)}+\delta E_h^{\rm pre}(\vy_h^{n+1},\vw_h) = 0 \quad \forall\,\vw_h\in\W_{h}^{\rm aver}(\vy_h^n),
\end{equation*}
where
$$
\delta E_h^{\rm pre}(\vy_h^{n+1},\vw_h) := \delta E_h^B(\vy_h^{n+1},\vw_h) + \delta  S_h(\vy_h^{n+1},\vw_h).
$$
In practice, the linear constraint encoded in $\mathbb W_h^{\rm aver}(\vy_h^n)$ is enforced using piecewise constant Lagrange multipliers $\Lambda_h: \Omega \rightarrow [\V_h^0]^{2\times 2}$ with $\Lambda_h^T = \Lambda_h$.

The resulting sequence of deformations is again energy decreasing, i.e.
\begin{equation} \label{e:energy_decay_pre}
	\frac{1}{\tau}\|\vy_h^{n+1}-\vy_h^n\|_{H_h^2(\Omega)}^2+E_{h}^{\rm pre}(\vy_h^{n+1})\le E_{h}^{\rm pre}(\vy_h^n)
\end{equation}
irrespective of the choice of the pseudo timestep $\tau>0$. 
Furthermore, because $\vy_h^{n+1} - \vy_h^n \in \W_{h}^{\rm aver}(\vy_h^n)$, we can show that the metric defect $\mathcal D_h^{\rm aver}$ defined in \eqref{def:Defect} is uniformly controlled, namely
\begin{equation} \label{rel:control_defect1}
	\mathcal{D}_h^{\rm aver}(\vy_h^n)\le \mathcal{D}_h^{\rm aver}(\vy_h^0)+c_1\tau\left(E_h^{\rm pre}(\vy_h^0)+c_2\right), \qquad n=1,2,\ldots,
\end{equation}
where the constant $c_1>0$ depends only on the constants in Poincar\'e-Friedrichs-type inequalities and $c_2\ge 0$ depends on the input data (namely $\mu$, $g$, $\vvarphi$, and $\Phi$) as well as the stabilization parameters $\gamma_0,\gamma_1$. We refer to \cite{BGNY2022b} for more details when $\Gamma^P = \emptyset$.

The control on the metric defect \eqref{rel:control_defect1} implies that upon setting $\varepsilon = \mathcal{D}_h^{\rm aver}(\vy_h^0)+c_1\tau\left(E_h^{\rm pre}(\vy_h^0)+c_2\right)$, $\vy_h^n \in \mathbb A_{h,\varepsilon}^{\rm aver}$ for all $n=1,2,...$. In particular, the recursive algorithm produces a sequence of deformation in $\mathbb A_{h,\varepsilon}^{\rm aver}$ with decreasing energy as desired.

\subsubsection{Bilayer Energy}

It remains to discuss the constrained minimization problem \eqref{def:bilayer_constr_min} associated with the bilayer energy 
$$
E_h^{\rm bil}=E_h^B+S_h-N_h,
$$ 
see \eqref{e:bilayer_discrete} (recall that $\vf=\mathbf{0}$). 
This case contains the two difficulties encountered above: the G\^ateau derivative $\delta E_h^{\rm bil}$ is not bilinear due to the cubic term $N_h$ and the minimization is subject to an approximate isometry constraint. Note that the latter is enforced at the cells' barycenters rather than in average, see \eqref{def:Discrete_Admissible_g}. 

Nevertheless, the algorithm considered in \cite{BNY2022}, which in turn is inspired from \cite{BP2021}, combines the ideas of the previous two cases. We define for $\vy_h^n \in [\V_h^k]^3$ the set
\begin{equation*}
	\W_{h}^{\rm bary}(\vy_h^n) := \left\{\vw_h\in[\V_h^k]^3:\,\, \left[(\nabla_h\vw_h)^T\nabla_h\vy_h^n+(\nabla_h\vy_h^n)^T\nabla_h\vw_h\right](x_T)=0 \quad \forall\, \K\in\Th\right\}
\end{equation*}
and seek $\vy_h^{n+1} \in [\V_h^k]^3$ such that $\vy_h^{n+1}-\vy_h^n \in \W_{h}^{\rm bary}(\vy_h^n)$
and
\begin{equation} \label{def:incrt_variational_bilayer}
	\tau^{-1}(\vy_h^{n+1}-\vy_h^n,\vw_h)_{H_h^2(\Omega)} + \delta E_h^B(\vy_h^{n+1},\vw_h) + \delta  S_h(\vy_h^{n+1},\vw_h)= N_h(\vy_h^n;\vw_h)
\end{equation}
for all $\vw_h\in\W_{h}^{\rm bary}(\vy_h^n)$.
Here 
\begin{equation*}
	\begin{split}
		N_h(\vy_h^n;\vw_h)  & :=   \alpha \sum_{i,j=1}^2 \sum_{\K \in \Th} |\K| \left[\big(\overline{H}_{h}(\vw_h)\big)_{ij} \cdot \left(\partial_1\vy_h^n\times\partial_2\vy_h^n\right) Z_{ij}\right](x_T) \\
		& +\alpha \sum_{i,j=1}^2 \sum_{\K \in \Th} |\K| \left[\big(\overline{H}_{h}(\vy_h^n)\big)_{ij} \cdot \left(\partial_1\vw_h\times\partial_2\vy_h^n\right) Z_{ij}\right](x_T) \\
		& +\alpha \sum_{i,j=1}^2 \sum_{\K \in \Th} |\K| \left[\big(\overline{H}_{h}(\vy_h^n)\big)_{ij} \cdot \left(\partial_1\vy_h^n\times\partial_2\vw_h\right) Z_{ij}\right](x_T)
	\end{split}
\end{equation*}
so that \eqref{def:incrt_variational_bilayer} corresponds to an explicit treatment of the cubic term $N_h$. 
As in the prestrain case, the linear constraint encoded in $\mathbb W_h^{\rm aver}(\vy_h^n)$ is enforced using piecewise constant Lagrange multipliers $\Lambda_h: \Omega \rightarrow [\V_h^0]^{2\times 2}$ with $\Lambda_h^T = \Lambda_h$.

As it turns out, except for the mild requirement  $\tau \lesssim |\ln(h_{\min})|^{-1}$, the explicit treatment of $N_h$ does not affect the convergence property  (in \cite{BBN2017} its treatment was implicit thereby requiring a fixed-point iteration at each step). In particular, the energy decay property \eqref{e:energy_decay_pre} holds for $E_h^{\rm bil}$ with a factor $\frac{1}{2\tau}$ instead of $\frac 1 \tau$ in front of the metric term and the isometry defect satisfies
\begin{equation} \label{rel:control_defect2}
	\mathcal{D}_h^{\rm bary}(\vy_h^n)\le \mathcal{D}_h^{\rm bary}(\vy_h^0)+c_1\tau|\log(h_{\min})|\left(E_h(\vy_h^0)+c_2\right), \qquad n=1,2,\ldots,
\end{equation}
where similarly to \eqref{rel:control_defect1}, the constant $c_1>0$ depends only on the constants in Poincar\'e-Friedrichs-type inequalities while $c_2\ge 0$ depends on the input data (in this case $Z$, $\vvarphi$, and $\Phi$) as well as the stabilization parameters $\gamma_0,\gamma_1$. Refer to \cite{BNY2022} for additional details in the case $\Gamma^P = \emptyset$. 
Again, this implies that all the iterations of the algorithms belong to $\mathbb A_{\varepsilon,h}^{\rm bary}$ provided that $\varepsilon \geq \mathcal{D}_h^{\rm bary}(\vy_h^0)+c_1\tau|\log(h_{\min})|\left(E_h(\vy_h^0)+c_2\right)$.

\subsection{Preprocessing Step} \label{sec:preprocessing}

In this section we introduce a preprocessing strategy which can be used to construct the initial deformation $\vy_h^0$ required by the gradient flows detailed in Section~\ref{sec:GF}. 
The estimates \eqref{rel:control_defect1} and \eqref{rel:control_defect2} for the metric defect satisfied by the constraint minimization problem suggest that the initial condition should satisfy approximately the boundary conditions (to reduce the energy) and have a small metric defect.
Therefore, the preprocessing strategy consists of two main modules dedicated to each aspect.
Note that in some cases, a suitable initial deformation is known a priori, as in the case of clamped boundary conditions and isometry constraint $(g=I_2)$, which are satisfied by the \emph{identity map} $\vy_h^0(x_1,x_2)=(x_1,x_2,0)^T$. However, in general, the initial deformation is not accessible.

\vskip1ex\paragraph{\bf Boundary Conditions Preprocessing.}
The \emph{BC preprocessing} consists in solving a bi-Lapacian problem to obtain a deformation satisfying (approximately) the prescribed boundary conditions. 
We seek an approximation of the solution to the bi-Laplacian problem
\begin{equation} \label{def:pb_BC}
	\left\{\begin{array}{rcll}
		\Delta^2\hat\vy & = & \hat\vf & \mbox{in } \Omega \\
		\nabla\hat\vy & = & \Phi & \mbox{on } \Gamma^D \\
		\hat\vy & = & \vvarphi & \mbox{on } \Gamma^D\cup\Gamma^M \\
		\hat\vy(\vx_i) & = & \vvarphi_i & \mbox{for all } \vx_i\in\Gamma^P,
 	\end{array}\right.
\end{equation}
supplemented with the following \emph{natural} boundary conditions
\begin{equation} \label{def:pb_BC2}
\left\{
\begin{array}{rcll}
\quad (D^2\hat y_m)\vn & = & \mathbf{0} & \mbox{on } \partial\Omega\setminus\Gamma^D \\
\nabla(\Delta \hat y_m)\cdot\vn & = & 0 & \mbox{on } \partial\Omega\setminus(\Gamma^D\cup\Gamma^M)
 \end{array}\right.
\end{equation}
for $m=1,2,3$.
Note that $\hat\vf$ is a \emph{fictitious} forcing term. Its value is irrelevant but can be employed to obtain a non-planar deformation, namely a deformation $\hat\vy=(y_1,y_2,y_3)$ with $y_3\not\equiv 0$, in particular when $\Gamma^D=\emptyset$.
This happens to be critical because the main discrete gradient flow preserves planar configurations (when $\vf =\mathbf{0}$) and can therefore not reach non-planar minimizers. 

The approximation $\hat \vy_h$ of $\hat \vy$ is obtained using an LDG approach, that is
$$
\hat \vy_h := \argmin_{\vy_h \in [\V_h^k]^3} E_h^{\rm pre}(\vy_h),
$$
where $E_h^{\rm pre}$ is given by \eqref{E:final-prestrain-bending-discrete} with $\mu=6$, $\lambda=0$, and $g=I_2$. In particular, the solution $\hat \vy_h$ is obtained upon solving the associated Euler-Lagrange equation. 

\vskip1ex\paragraph{\bf Metric Constraint Preprocessing.} 

The \emph{metric preprocessing} aims at minimizing the metric defect $\mathcal{D}_h^i$, $i\in \{\rm aver, bary\}$.
It consists of minimizing the \emph{simplified stretching energy}
\begin{equation} \label{def:E_stretching_PP}
	E_h^{\rm str}(\vy_h) := \frac{1}{2}\int_{\Omega}\left|(\nabla_h \vy_h)^T\nabla_h\vy_h-g\right|^2
\end{equation}
over $[\V_h^k]^3$.
The stretching energy $E_h^{\rm str}$ is not quadratic, but its reduction can be performed using an $H_h^2(\Omega)$ discrete gradient flow coupled with a linearization similar to the one for $E_{s,h}^{\rm pre}$. Setting $\tilde\vy_h^0=\hat\vy_h$, the deformation obtained by the BC preprocessing algorithm, recursively computes $\tilde\vy_h^{n+1}\in[\V_h^k]^3$ as the solution to
\begin{equation}\label{e:incr_pp}
	\tilde\tau^{-1}(\tilde\vy_h^{n+1}-\tilde\vy_h^n,\vw_h)_{H_h^2(\Omega)}+\delta E_h^{\rm str}(\tilde\vy_h^n;\tilde\vy_h^{n+1},\vw_h) = 0 \quad \forall\, \vw_h\in[\V_h^k]^3,
\end{equation}
where $\tilde\tau>0$ is a given (sufficiently small) pseudo time-step and where for $\vv_h,\vw_h\in[\V_h^k]^3$
\begin{equation*}
	\delta E_h^{\rm str}(\tilde\vy_h^n; \vv_h,\vw_h) := \int_{\Omega}\left((\nabla_h\vv_h)^T\nabla_h\vw_h+(\nabla_h\vw_h)^T\nabla_h\vv_h\right):\left((\nabla_h\tilde\vy_h^n)^T\nabla_h\tilde\vy_h^n-g\right).
\end{equation*}
The gradient flow produces iterates with decreasing energy provided that $\tilde\tau$ is sufficiently small compared to $E_h^{\rm str}(\tilde\vy_h^0)^{-1}$ and $h_{\rm min}$, see \cite{BGNY2022b}. Similar to the main discrete gradient flow, one drawback of the metric constraint preprocessing is that planar configurations are local minimizers of \eqref{def:E_stretching_PP} regardless of the target metric $g$, see \cite{BGNY2022a}. A non-planar deformation can nevertheless be obtained using a non-planar initial deformation $\tilde\vy_h^0$. Such an initial deformation can be generated, for instance, by solving \eqref{def:pb_BC}-\eqref{def:pb_BC2} with $\hat \vf\not\equiv\mathbf{0}$. Moreover, the metric constraint preprocessing does not necessarily construct a deformation with finite energy $E_{h}^{\rm pre}$. If needed, a bending term could be added to the energy \eqref{def:E_stretching_PP} to generate iterates with a uniformly bounded $H_h^2(\Omega)$ semi-norm, see \cite{BGNY2022b}.

\subsection{Accelerated Algorithm} \label{sec:acceleration}

It is well-documented that gradient flows exhibit slow convergence towards minimizers. Several accelerated algorithms have been introduced to improve the convergence rate of the standard gradient descent algorithm, but they are mainly available for convex problems. 

Below we design an accelerated algorithm suited for the minimization of the preasymptotic energy and the stretching energy used in the preprocessing process. The proposed algorithm is inspired by Nesterov's original accelerated gradient method \cite{Nesterov1983}, see also \cite{Nesterov2004,nesterov2018lectures}.

Typical results for the minimization of functions $f:\mathbb{R}^d\rightarrow\mathbb{R}$ is an $\mathcal{O}(n^{-2})$ convergence in $f$ when $f$ is a $C^1$ convex function with Lipschitz gradient. Here $n$ denotes the number of iterations. This is a significant improvement over the standard gradient method for which the convergence rate is only $\mathcal{O}(n^{-1})$.

The Fast Iterative Shrinkage/Thresholding Algorithm (FISTA) is an extension to when $f$ is the sum of two convex and lower-semicontinuous functions, one of class $C^1$ with Lipschitz gradient and the other potentially non-smooth. The convergence of the FISTA algorithm is also quadratic in the number of iterations, see \cite{BT09} and \cite{CD15}.

Inspired by the FISTA algorithm (case $f_1\equiv 0$), but adapted to the functional setting, our \emph{accelerated algorithm} for the minimization of the preasymptotic energy $E_{s,h}^{\rm pre}$ and the stretching energy $E_h^{\rm str}$ used in the metric constraint preprocessing reads:
Given $\vy_h^0\in[\V_h^k]^3$, set $\vv_h^0=\vy_h^0$. Then for $n=0,1,\ldots$ do 
\begin{enumerate}[a)]
	\item Find $\vy_h^{n+1} \in [\V_h^k]^3$ satisfying \eqref{def:incrt_variational_BS} or \eqref{e:incr_pp} with $\vy_h^n$ replaced by $\vv_h^n$; 
	\item Set $\vv_h^{n+1}=\vy_h^{n+1}+\eta_{n+1}(\vy_h^{n+1}-\vy_h^{n})$, where $\eta_{n+1}=\frac{t_{n+1}-1}{t_{n+2}}$ and $\{t_n\}_{n\ge 1}$ satisfies
	\begin{equation*}
		t_1=1 \quad \mbox{and} \quad  \quad t_{n+1} = \sqrt{t_n^2+\frac{1}{4}}+\frac{1}{2} \quad \mbox{for } n=0,1,\ldots. 
	\end{equation*}
\end{enumerate}
As in the case of the main gradient flow, a linearization is performed at $\vv_h^{n}$ to avoid a nonlinear system at each iteration.

Note that the accelerating algorithm does not have the energy decreasing property that the gradient descent algorithm does. This can be observed experimentally, see Figure~\ref{fig:energydecay_bubble} below. However, the algorithm does appear to converge when the appropriate step size is chosen.

Finding analytic estimates for the convergence rate of the proposed algorithm remains an open problem. In Section~\ref{sec:num_res_acc} we observe at least an $\mathcal{O}(n^{-2})$ convergence rate and illustrate the significant gain when the accelerated version of the algorithm is used.

\subsection{Dynamics: Time Discretization} \label{subsec:dynamics}

Several numerical experiments presented in Section~\ref{sec:results} are dynamical in the sense that the data (boundary conditions and external forces) vary over the time interval $[0,T]$. The minimizers of the different energies $\vy$ are thus time dependent as well.

For a positive integer $\mathtt{M}$, we let $\Delta t:=T/\mathtt{M}$ be the physical time-step (as opposed to the pseudo-time step used in the gradient flows) and consider the points $t_m:=m\Delta t$ for $m=0,1,\ldots,\mathtt{M}$. We assume that the elastic relaxation is faster than the time relaxation, i.e. at each time step, the thin elastic structures have minimal energies. In particular, the deformations $\vy_h(t_m) \in [\V_h^k]^3$ at each time $t_m$ are discrete minimizers of the discrete energy defined using the input data evaluated at $t=t_m$. 

The dynamic algorithm reads:

\vspace*{0.2cm}
\noindent \textbf{Initialization} ($t=0$):
Obtain $\hat \vy_h^{(0)}$ using the BC preprocessing algorithm with the prescribed conditions evaluated at $t=0$; Set $\vy_h^{(0)} = \hat \vy_h^{(0)}$;

\vspace*{0.2cm}
\noindent \textbf{Dynamics} ($t\in(0,T]$): for $m=1,2,\ldots,\mathtt{M}$, do
\begin{enumerate}[i)]
	\item Obtain $\delta \hat \vy_h^{(m)}:= \hat \vy_h^{(m)}-\vy_h^{(m-1)}$ using the BC preprocessing algorithm with the increment boundary conditions
	\begin{equation*} 
		\begin{split}
			\nabla \delta \hat \vy_h^{(m)} &= \Phi(t_m)-\Phi(t_{m-1}), \quad \mbox{on } \Gamma^D \\
			\delta \hat \vy_h^{(m)} &= \vvarphi(t_m)-\vvarphi(t_{m-1}), \quad \mbox{on } \Gamma^D\cup\Gamma^M  \\
			\delta \hat \vy_h^{(m)}(\vx_i) &= \vvarphi_i(t_m)-\vvarphi_i(t_{m-1}), \quad \mbox{for } i=1,2,\ldots,P;
		\end{split}
	\end{equation*}
	
	\item Starting from $\hat \vy_h^{(m)}$, obtain $\tilde \vy_h^{(m)}$ using the metric preprocessing algorithm to satisfy approximately the metric constraint;
	\item Starting from $\tilde \vy_h^{(m)}$, obtain $\vy_h^{(m)}$ using the (main) gradient flow to minimize the discrete energy with the data evaluated at $t=t_m$.
\end{enumerate}

It is important to point out that the BC preprocessing procedure embedded in the dynamic algorithm above determines the increment $\delta \hat \vy_h^{(m)}:= \hat \vy_h^{(m)}-\vy_h^{(m-1)}$ (rather than the preprocessed deformation $\hat \vy_h^{(m)}$ directly) by solving \eqref{def:pb_BC}-\eqref{def:pb_BC2} with incremental boundary conditions. 
This is done to take advantage of the previous step and, in particular, to avoid using a costly metric preprocessing algorithm starting from scratch at each time step.

\section{Numerical Experiments} \label{sec:results}

In this section, we illustrate the performance of the different algorithms. 
All the simulations are performed with $\mu=6$ and $\lambda=8$ (note that $\mu=6$ and $\lambda=0$ yield the factor 1/2 for the bilayer case) and without any external force (i.e. $\vf=\mathbf{0}$) except for the experiment in Section~\ref{sec:sphere}. The polynomial degree used for the approximated deformations $\vy_h$ is chosen to be $k=2$, and, unless specified otherwise, we take $\gamma_0=\gamma_1=1$ and $\gamma_2=10$ for the stabilization parameters. The main gradient flow ends when 
\begin{equation*}
	\tau^{-1}\left|E_h(\vy_h^{n+1})-E_h(\vy_h^{n})\right|\le tol
\end{equation*}
for a prescribed tolerance $tol$ specified for each numerical experiment below. For the metric preprocessing we stop when
\begin{equation*}
	\widetilde\tau^{-1}\left|\widetilde E_h(\tilde\vy_h^{n+1})-\widetilde E_h(\tilde\vy_h^{n})\right|\le \widetilde{tol} \quad \mbox{or} \quad \mathcal{D}_h^i(\tilde\vy_h^{n+1})\le \widetilde\varepsilon_0, \quad i\in\{\rm aver,bary\},
\end{equation*}
for some tolerances $\widetilde{tol}$ and $\widetilde\varepsilon_0$ specified below. Finally, for the experiments in Sections~\ref{sec:num_preasymptotic} and \ref{sec:sphere}, where the domain $\Omega$ consists of the unit disc, we use a quadratic mapping $F_T$ in \eqref{def:Vhr} to better approximate the domain.

\subsection{Preasymptotic} \label{sec:num_preasymptotic}
We first consider two experiments using the preasymptotic model. The first experiment illustrates the effect of the thickness of the plate on the final configuration, while the second demonstrates the advantages of the Nesterov-type acceleration discussed in Section~\ref{sec:acceleration}.
In both experiments the computational domain $\Omega$ is a disc of radius one, and we use $\tau=0.01$ for the gradient flow pseudo time-step.

\subsubsection{Disc with Oscillating Boundary}

Our first experiment is inspired by \cite{KVS2011,LM2021}, in which a hydrogel disc of negative Gaussian curvature is observed to develop more oscillations along the boundary as its thickness is reduced. 
The material is prestrained according to the metric
$$g=J^T(\widetilde g\circ \boldsymbol{\zeta})J,$$
where $J$ is the Jacobian matrix for the change of variables $(r,\theta)=\boldsymbol{\zeta}(x_1,x_2)$ from polar to Cartesian coordinates and $\widetilde g(r, \theta)$ is the first fundamental form of the following deformation
\begin{equation} \label{eqn:target_oscillating}
	\widetilde {\bf y}(r,\theta) = (r\cos(\theta), r\sin(\theta), 0.2r^4\sin(6\theta)).
\end{equation}
This deformation corresponds to a disk with six wrinkles.

For the discretization, our subdivision consists of 1280 quadrilaterals and a total of $34560 = 1280 \cdot 9 \cdot 3$ degrees of freedom. The initial deformation ${\bf y}_h^0$ is taken to be the continuous Lagrange interpolant of (\ref{eqn:target_oscillating}). We set $tol=10^{-8}$ for the stopping criterion and run the experiment for $s^2 = 10^{-1}$, $10^{-2}$, $10^{-3}$, and $0$.

The computed deformations minimizing the preasymptotic energy $E^{\rm pre}_s$ are provided in Figure~\ref{fig:results_oscillating}. In accordance with the laboratory experiments \cite{KVS2011,LM2021}, the number of wrinkles increases from 0 to 6 as $s$ decreases.

\begin{figure}[htbp]
	\begin{center}
		\includegraphics[width=6.27cm]{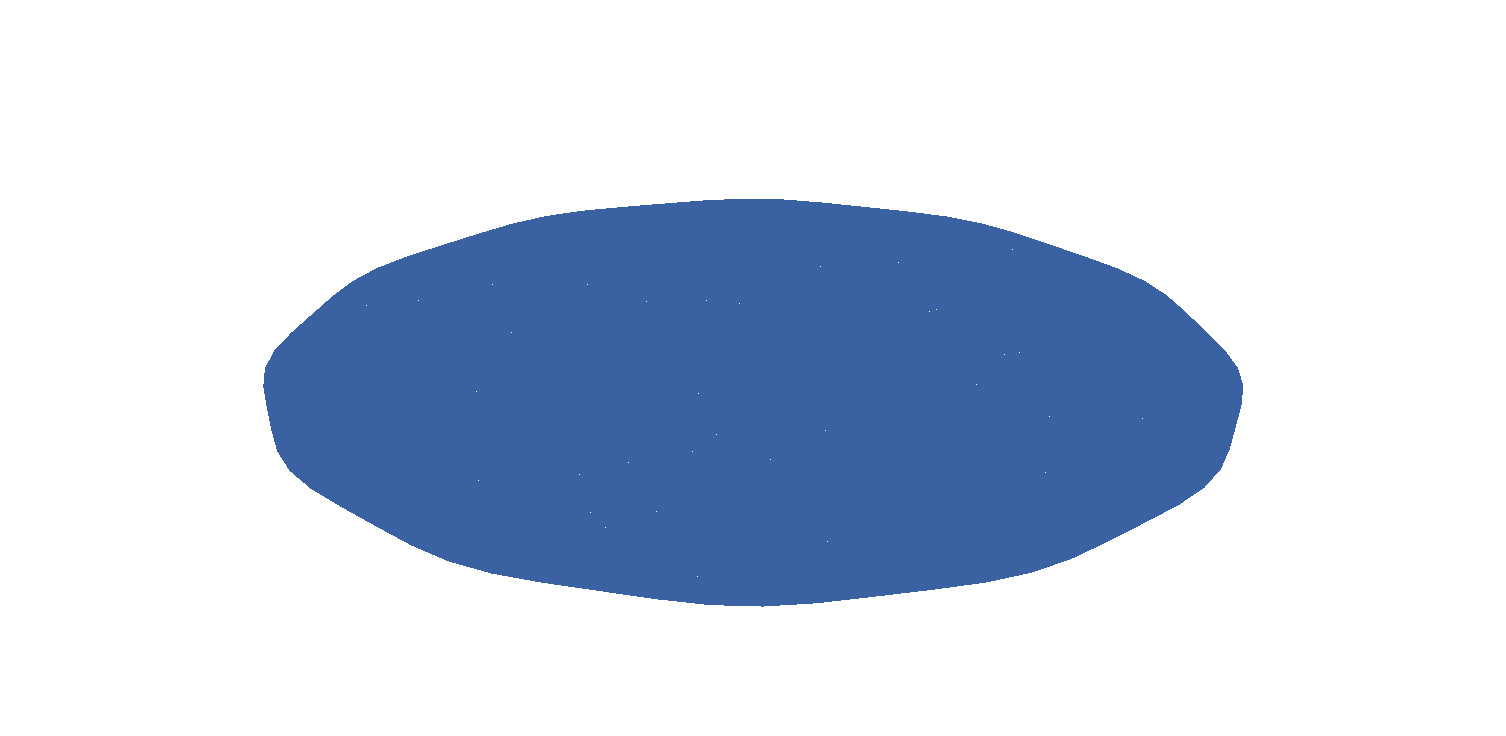}
		\includegraphics[width=6.27cm]{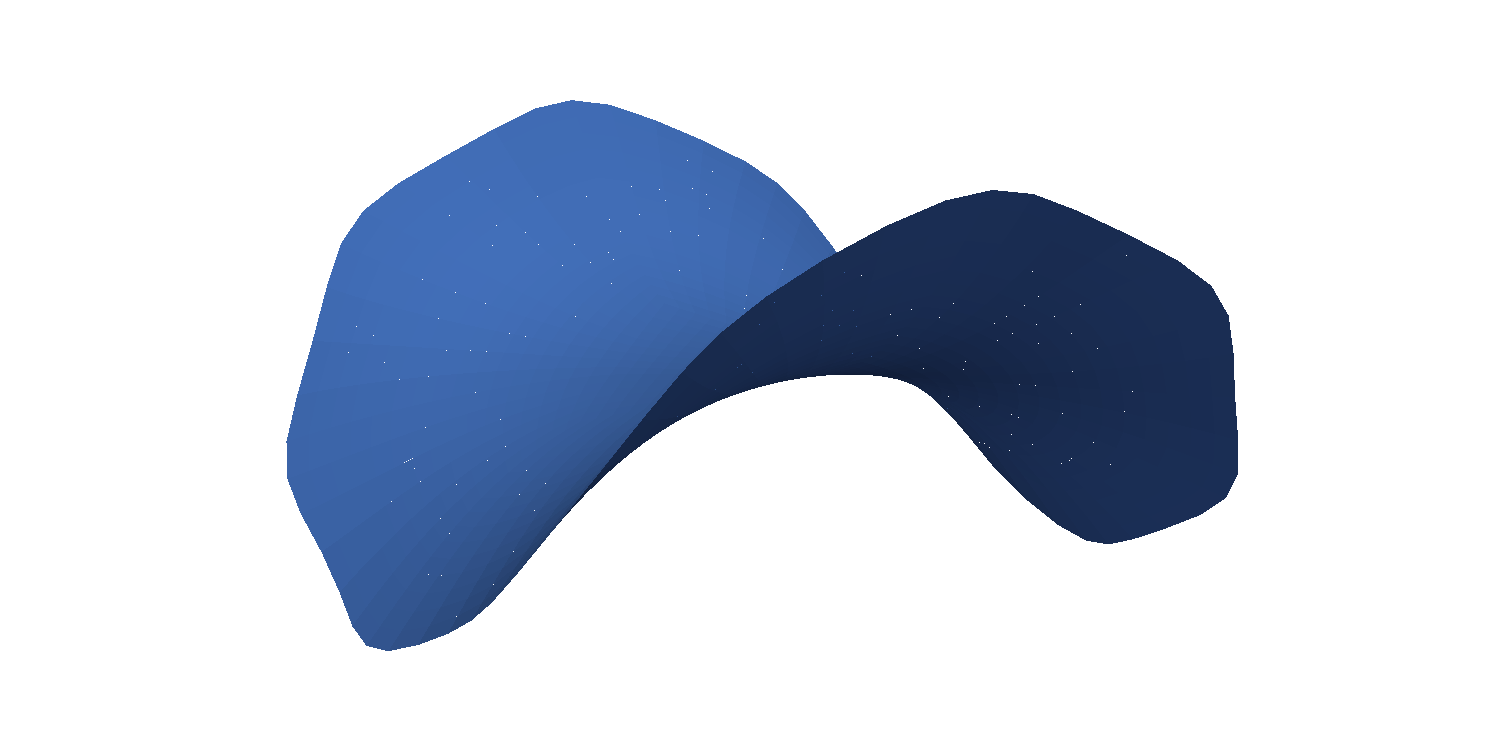} \\
		
		\includegraphics[width=6.27cm]{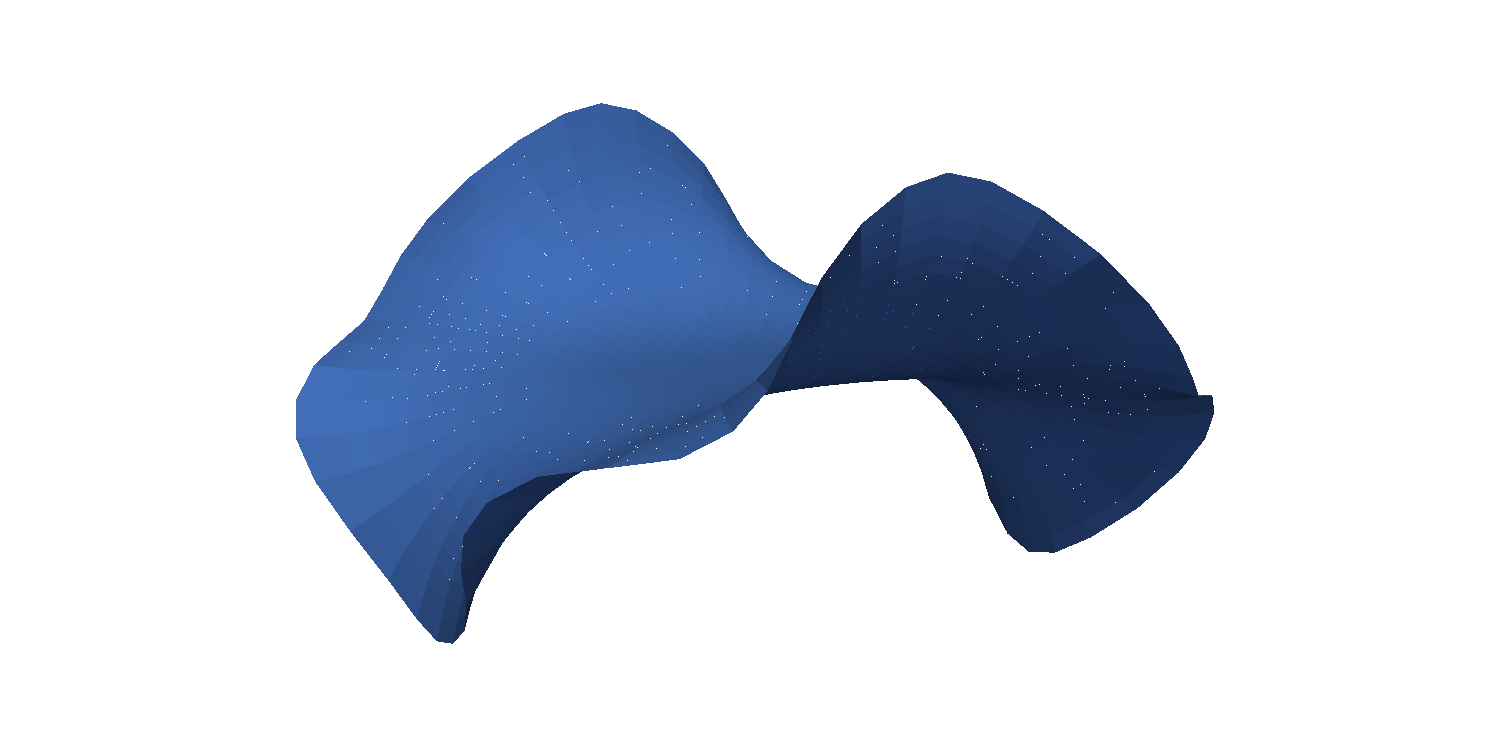}
		\includegraphics[width=6.27cm]{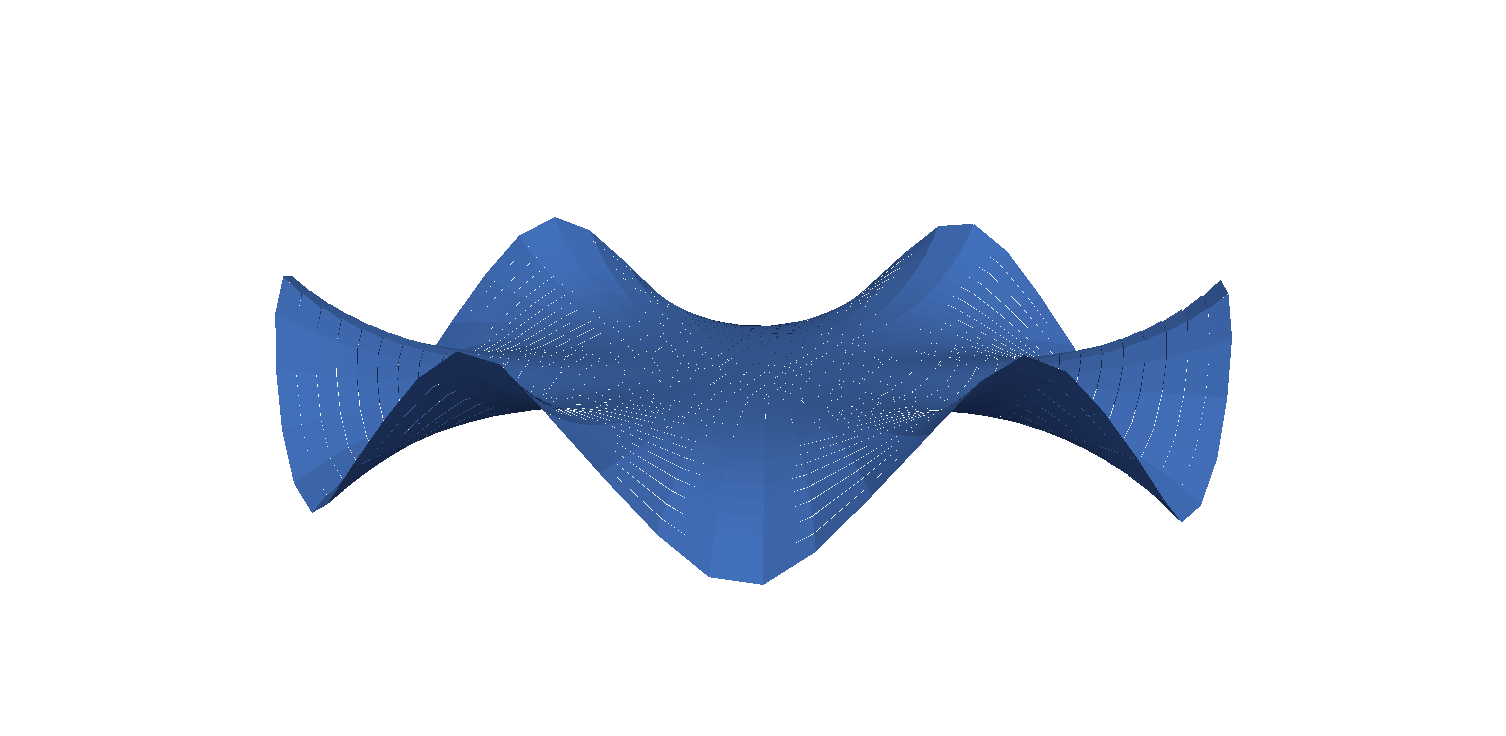}
	\end{center}
	\caption{Final configurations for the oscillating boundary experiment. Left to right and top to bottom: $s^2 = 10^{-1}$, $10^{-2}$, $10^{-3}$, and $0$.}
	\label{fig:results_oscillating}
\end{figure}

\subsubsection{Accelerated Algorithm} \label{sec:num_res_acc}

This example illustrates the benefits of using the accelerated algorithm introduced in Section~\ref{sec:acceleration} over our original gradient flow. This time the material is characterized by a ``bubble'' metric with positive Gaussian curvature, namely
\begin{equation*}
	g(x_1,x_2) = \left(\begin{array}{cc}
		1+\alpha \frac{\pi^2}{4}\cos(\frac{\pi}{2}(1-r))^2 \frac{x_1^2}{r^2} & \alpha \frac{\pi^2}{4}\cos(\frac{\pi}{2}(1-r))^2 \frac{x_1 x_2}{r^2} \\ \alpha\frac{\pi^2}{4}\cos(\frac{\pi}{2}(1-r))^2 \frac{x_1 x_2}{r^2} & 1+\alpha \frac{\pi^2}{4}\cos(\frac{\pi}{2}(1-r))^2 \frac{x_2^2}{r^2}
	\end{array}\right)
\end{equation*}
with $r = \sqrt{x_1^2 + x_2^2}$ and $\alpha = 0.2$. 

The mesh for this experiment consists of 320 quadrilaterals and a total of $8640 = 320\cdot 9 \cdot 3$ degrees of freedom. The initial deformation is taken to be the continuous Lagrange interpolant of the shallow paraboloid $-\frac{x_1^2}{16}+\frac{x_2^2}{16}$. This gives a slight initial bending that allows the simulation to find a non-flat minimizer. (Recall from Section \ref{sec:preprocessing} that the flat configuration is a local minimizer of the preasymptotic energy.) For the thickness of the material, we use the values $s^2 = 10^{-3}$, $10^{-4}$, $10^{-5}$, and $0$ and we set $tol=10^{-6}$ as stopping criterion.

The discrete stretching energy $E_h^S$ and bending energy $s^2(E_h^B+S_h)$ along with the number of iterations needed to reach the tolerance $tol$ are reported in Table~\ref{table:bubble_iterations} when using the standard gradient flow and the accelerated algorithm. 
We observe that the accelerated algorithm finishes in much fewer iterations for all values of the thickness $s$. 
The approximate preasymptotic energy $E^{\rm pre}_{s,h}=E_h^S+s^2 (E_h^B +S_h)$ versus the algorithm iteration number $n$ is depicted in Figure~\ref{fig:energydecay_bubble} for the different values of $s$ considered. 
When using a standard gradient flow, we observe numerically that the energy decay is like $\mathcal{O}(n^{-1})$, while the accelerated algorithm exhibits a decay similar to $\mathcal{O}(n^{-2})$. The latter is the expected convergence rate when using Nesterov-type accelerations. 

\begin{center}
	\begin{table}[htbp!]
		\begin{tabular}{ |c||c|c|c||c|c|c| } 
			\hline
			&\multicolumn{3}{|c||}{Discrete Gradient Flow} & \multicolumn{3}{|c|}{Accelerated Gradient Flow} \\
			$s^2$ & $E_h^S$ & $E_h^B$ & Iterations & $E_{h}^S$  & $E_h^B$ & Iterations \\ 
			\hhline{|=||=|=|=||=|=|=|}
			0 & 2.390e-4 & 0 & 17746 & 8.098e-6 & 0 & 1007 \\ 
			\hline
			$10^{-5}$ & 1.525e-4 & 2.685e-4 & 26784 & 2.694e-4 & 3.990e-4 & 359 \\ 
			\hline
			$10^{-4}$ & 4.479e-5 & 5.951e-4 & 37799 & 6.350e-6 & 4.095e-4 & 896 \\
			\hline
			$10^{-3}$ & 5.595e-5 &2.843e-3 & 39533 & 3.920e-5 & 2.741e-3 & 664 \\
			\hline
		\end{tabular}
		\caption{Final energy and number of iterations for different values of $s^2$, with and without acceleration.} \label{table:bubble_iterations}
	\end{table}
\end{center}

\begin{figure}[htbp!]
	\begin{center}
		\includegraphics[width=6.27cm]{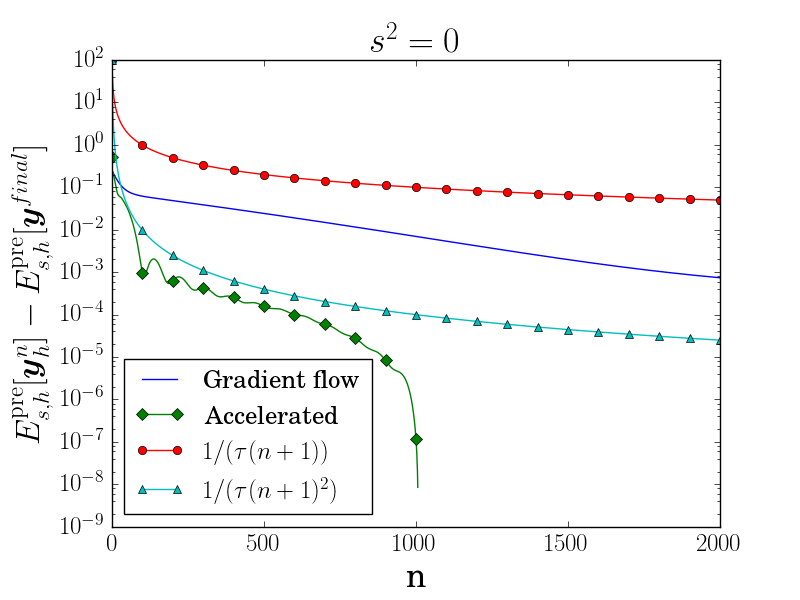}
		\includegraphics[width=6.27cm]{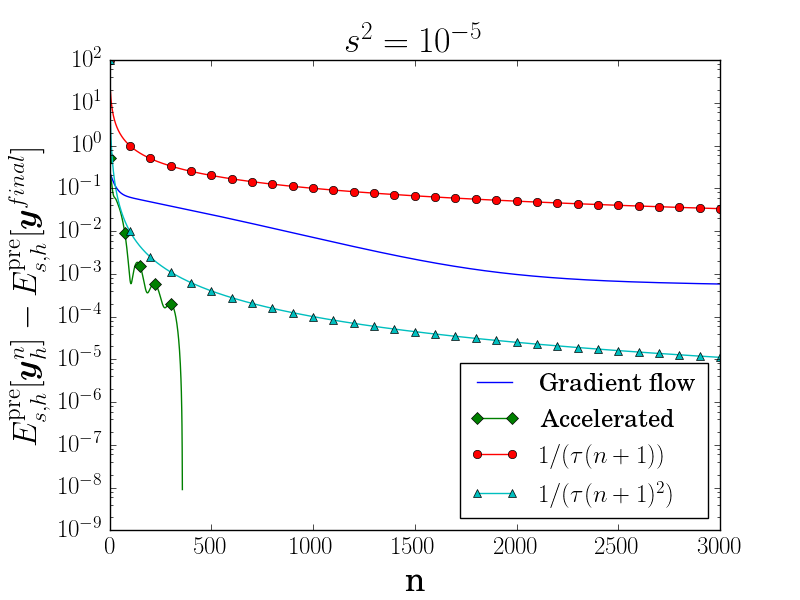} \\
		
		\includegraphics[width=6.27cm]{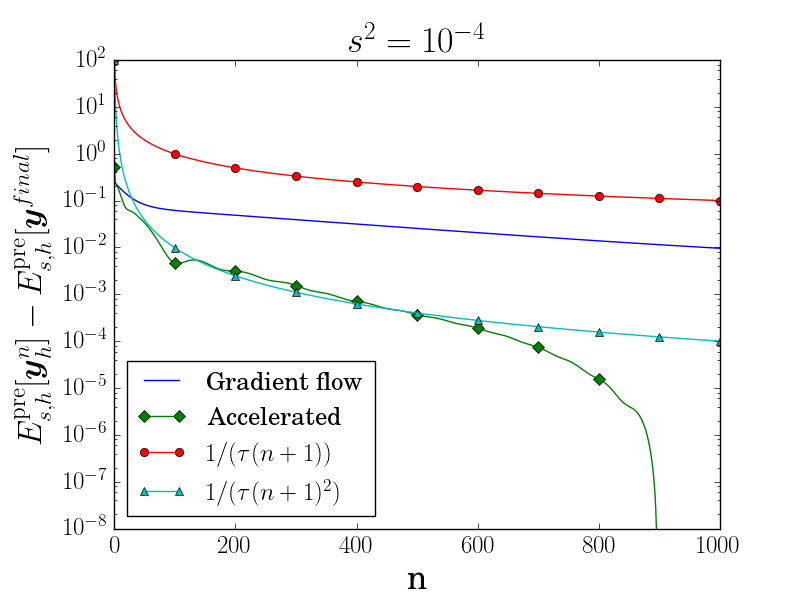}
		\includegraphics[width=6.27cm]{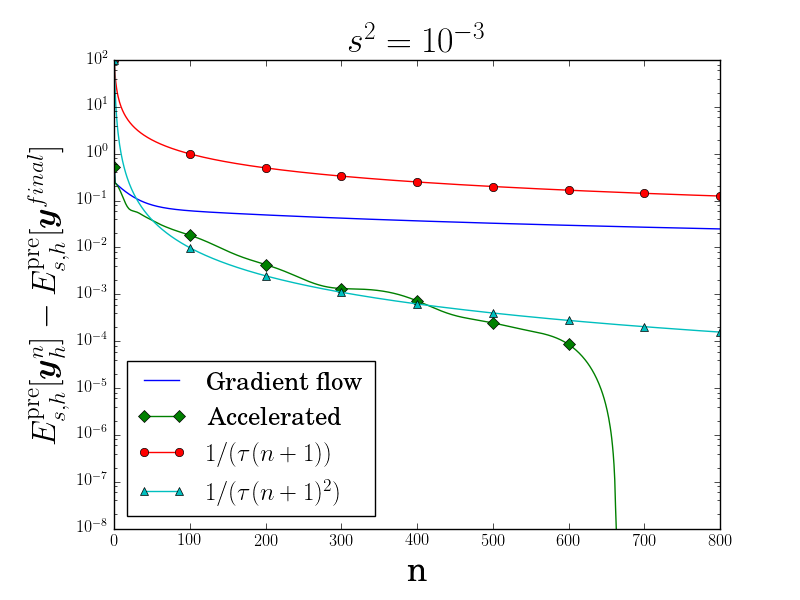}
	\end{center}
	\caption{Evolution of $E_{s,h}^{\rm pre}$ during the minimization when using the standard gradient flow and the accelerated algorithm. The energy are reported for the first few iterations to illustrate better the convergence of the acceleration.} \label{fig:energydecay_bubble}
\end{figure}

\subsection{Bilayer and Folding}

We consider here two experiments taken from \cite{TCWSTBM2020}, referred to as \emph{diamond} and \emph{bird} in the following.
They illustrate the rigidity, robustness, and great variety of shapes achievable with this technology.
In these examples, the computational domain consists of a collection of subdomains delimited by creases where folding is allowed. The bilayer material is designed to yield a piecewise constant spontaneous curvature tensor $Z$ on each subdomain with $\alpha = 1$. 

In both experiments we set $\tau=0.1$ for the gradient flow pseudo time-step and $tol=10^{-3}$ for the stopping criterion. 

\subsubsection{Diamond}

The \emph{diamond} computational domain is provided in Figure~\ref{fig:geo_diamond} (left). It consists of the square $(-1.5,1.5)\times (-1.5,1.5)$ rotated counter-clockwise by an angle of $\pi/4$. The two dashed red curves represent the creases. They are quadratic B\'ezier curves with the origin as control point. The point $\vz_{1}$ is at distance $1/3$ from $\vx_1$ on the segment from $\vx_1$ to $\vx_2$, and similarly for the three other points. We denote by $Z_i \in \mathbb R^{2\times 2}$ the spontaneous curvature associated with each subdomain $i\in\{1,2,3\}$ and set
$$Z_1=Z_3= 0.6 I_2 \quad \mbox{and} \quad Z_2=- 0.6 I_2.$$

\begin{center}
	\begin{figure}[htbp!]
		
		\begin{multicols}{2}	
			
			\begin{tikzpicture}
				
				\draw (0.0,-3.0/1.4142) -- (3.0/1.4142,0) -- (0.0,3.0/1.4142) -- (-3.0/1.4142,0) -- (0.0,-3.0/1.4142);
				
				\draw [red,dashed] (-1.0/1.4142,-2.0/1.4142) .. controls (0.0,0.0) .. (-1.0/1.4142,2.0/1.4142);
				
				\draw [red,dashed] (1.0/1.4142,-2.0/1.4142) .. controls (0.0,0.0) .. (1.0/1.4142,2.0/1.4142);
				
				\draw (-1.1,0.0) node[right]{$1$};
				\draw (-0.9,0.0) circle (2ex);
				
				\draw (0.0,1.0) node[above]{$2$};
				\draw (0.0,1.25) circle (2ex);
				
				\draw (1.1,0.0) node[left]{$3$};
				\draw (0.9,0.0) circle (2ex);
				
				\draw (0.0,-3.0/1.4142) node[below]{$\vx_1$};
				\draw (3.0/1.4142,0) node[right]{$\vx_2$};
				\draw (0.0,3.0/1.4142) node[above]{$\vx_3$};
				\draw (-3.0/1.4142,0) node[left]{$\vx_4$};	
				
				\draw (1.0/1.4142,-2.0/1.4142-0.1) node[right]{$\vz_1$};
				\draw (1.0/1.4142,2.0/1.4142+0.1) node[right]{$\vz_2$};
				\draw (-1.0/1.4142,2.0/1.4142+0.1) node[left]{$\vz_3$};
				\draw (-1.0/1.4142,-2.0/1.4142-0.1) node[left]{$\vz_4$};
				
			\end{tikzpicture}
			
			\includegraphics[width=5cm]{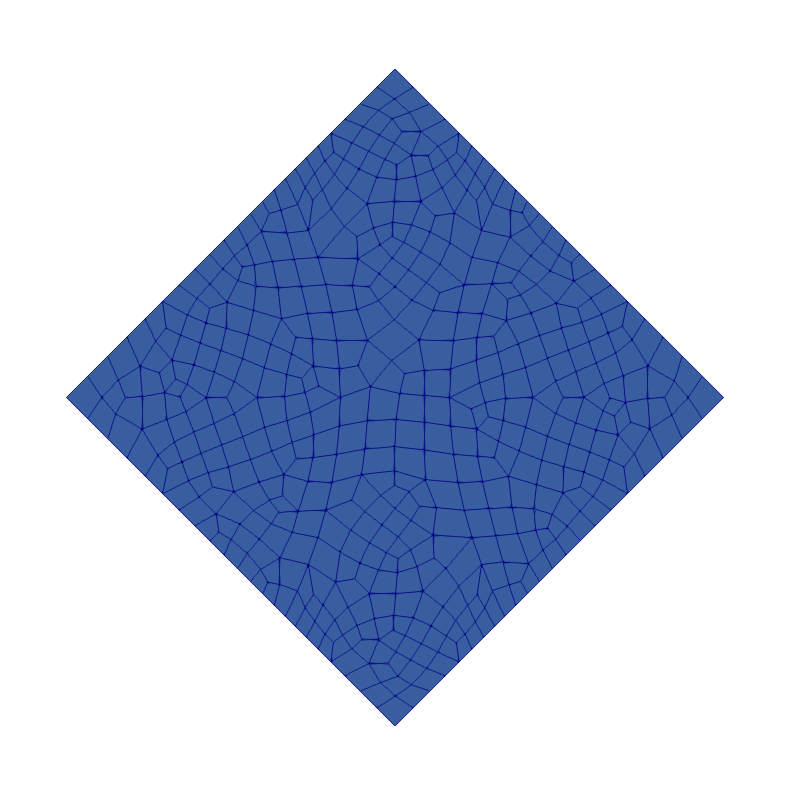}
			
		\end{multicols}
		\caption{Computational domain (left) and subdivision (right) for the \emph{diamond} setting. The dashed red curves are creases across which folding is possible.}\label{fig:geo_diamond}
	\end{figure}
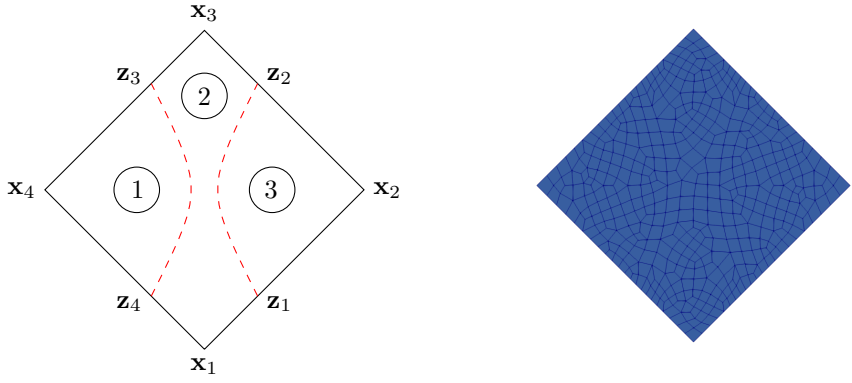
\end{center}

The subdivision  $\Th$ is depicted in Figure~\ref{fig:geo_diamond} (right).
It consists of $449$ quadrilaterals with maximal mesh size $0.331678$ (total of $13470$ degrees of freedom, $12123=449\cdot9\cdot3$ for the deformation $\vy_h$ and $1347=449\cdot 3$ for the Lagrange multipliers used to enforce the linearized isometry constraint). Moreover, the initial deformation $\vy_h^0$ is taken to be the identity map $\vy_h^0(\Omega)=\Omega\times\{0\}$ for which $E_h^{\rm bil}=3.24$ and $\mathcal D_h^{\rm bary}=8.72497\cdot 10^{-15}$. The deformations obtained at several steps of the gradient flow are provided in Figure~\ref{fig:iter_diamond}, including the final (equilibrium) deformation reached in $305$ iterations. For the final deformation, the energy and isometry defect are $E_h^{\rm bil}=0.679318$ and $\mathcal D^{\rm bary}_h=0.0334151$, respectively.

\begin{figure}[htbp]
	\begin{center}
		\includegraphics[width=5cm]{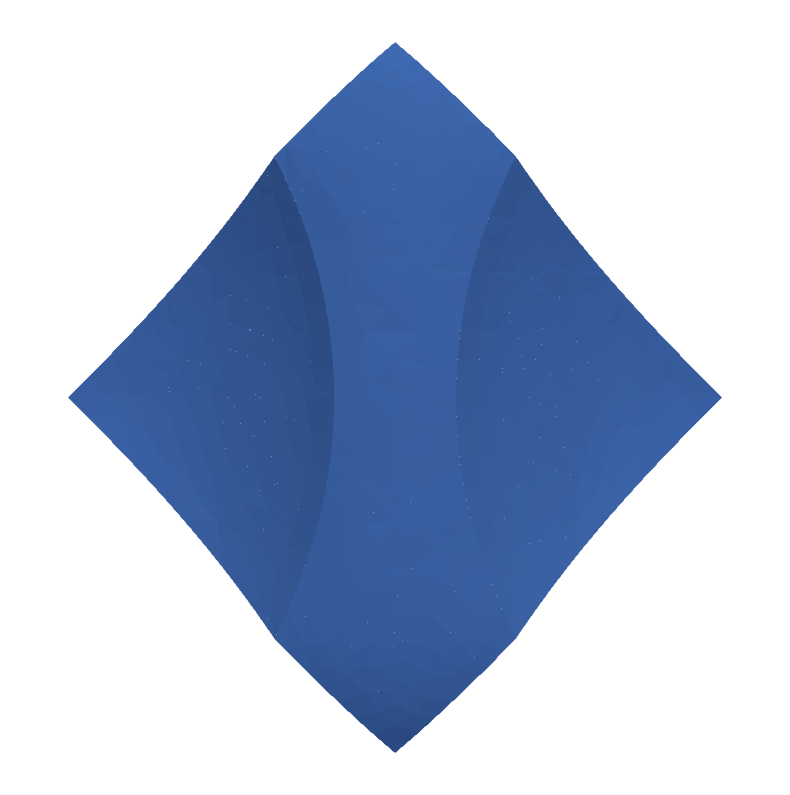}
		\includegraphics[width=5cm]{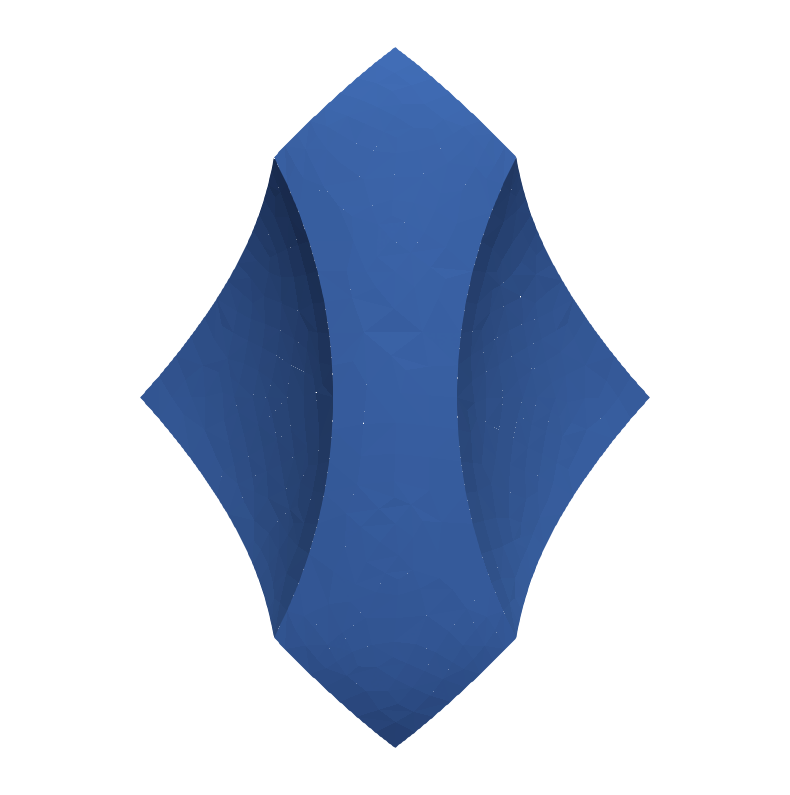} \\
		
		\hspace*{-0.5cm}\includegraphics[width=5cm]{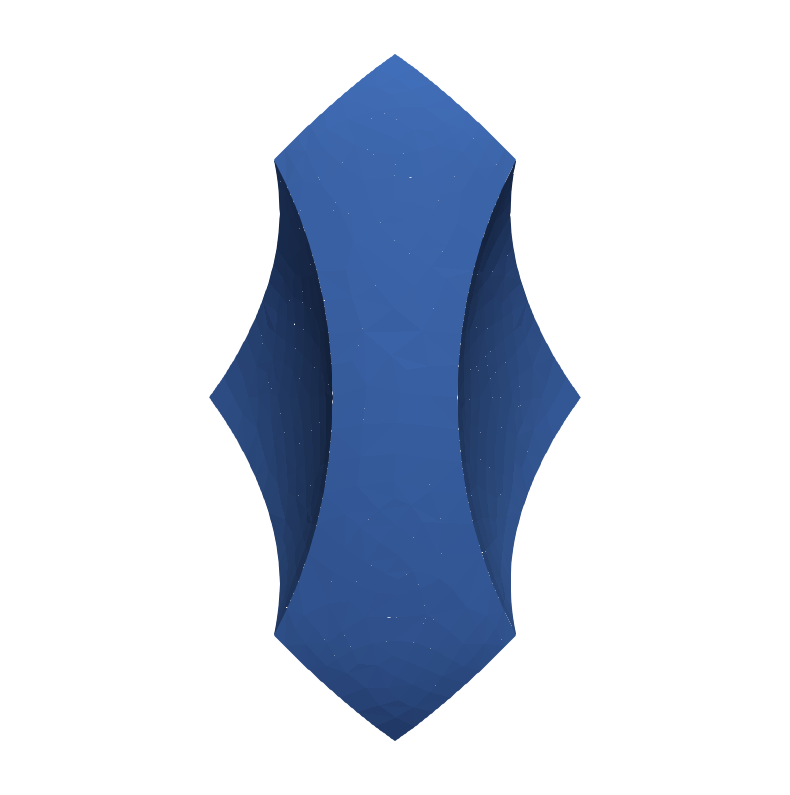}
		\hspace*{-1.2cm}\includegraphics[width=5cm]{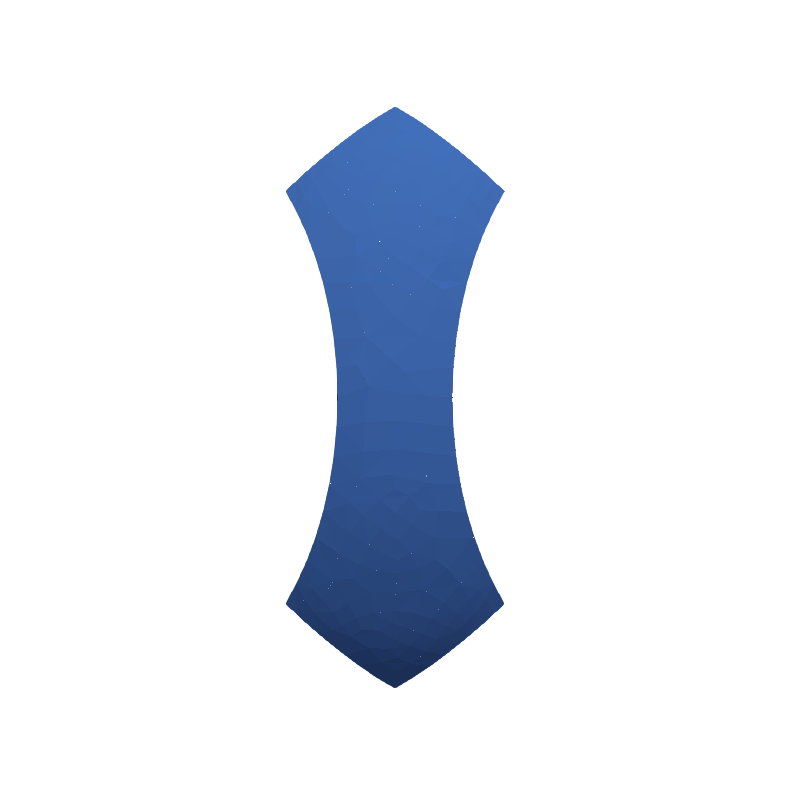}
		\hspace*{-1.2cm}\includegraphics[width=5cm]{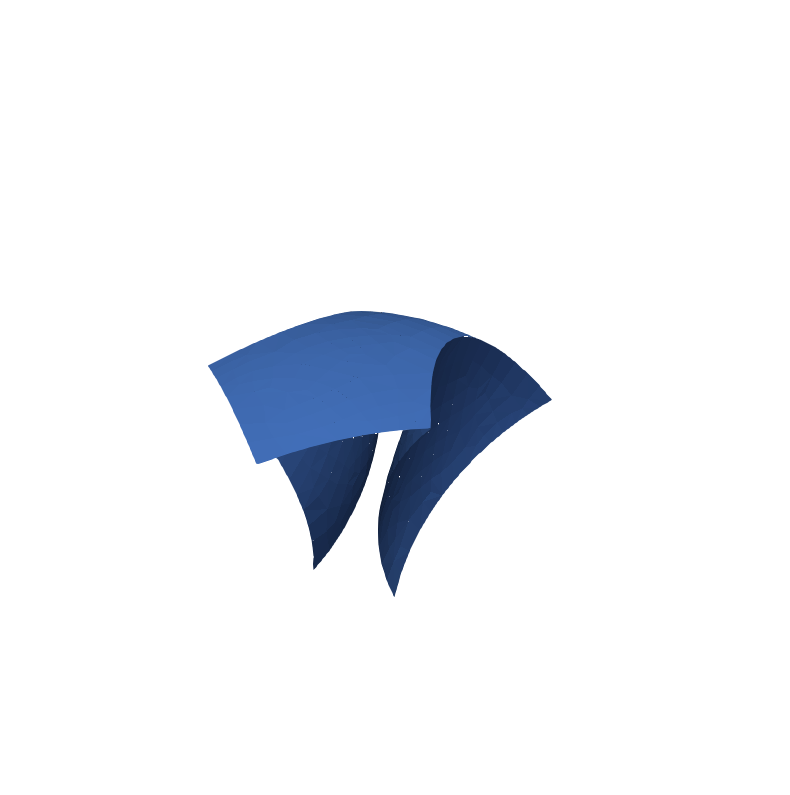}
	\end{center}
	\caption{Deformations for the \emph{diamond} experiment. Left to right and top to bottom: deformation obtained after $20$, $70$, $120$, $305$ (top view), and $305$ (side view) iterations.} \label{fig:iter_diamond}
\end{figure}

\subsubsection{Bird}

The \emph{bird} geometry is depicted in Figure~\ref{fig:geo_bird}. It consists of 26 curves: 17 for the boundary of $\Omega$ (black plain curves), 7 creases (red dashed curves), and 1 extra curve on the tail (black dotted curve) which is only used for the construction of the subdivision (i.e. no folding is possible across this curve). All the curves are quadratic or cubic B\'ezier curves, i.e. obtained using one or two control points. The latter are chosen so that we obtain a geometry similar to \cite{TCWSTBM2020}. In this example, the piecewise constant spontaneous curvature tensors $Z_1,...,Z_9$ are given on each subdomain by

\begin{equation*}
	Z_1=Z_5=\left(\begin{array}{cc}
		0 & 0 \\ 0 & 0.4
	\end{array}\right), \quad Z_2=Z_4=\left(\begin{array}{cc}
		0 & 0 \\ 0 & -0.3
	\end{array}\right), \quad Z_3=0.7I_2,
\end{equation*}
\begin{equation*}
	Z_6=\left(\begin{array}{cc}
		-0.2 & 0 \\ 0 & 0
	\end{array}\right), \quad Z_7=\left(\begin{array}{cc}
		0.7 & 0 \\ 0 & 0
	\end{array}\right) \quad \mbox{and} \quad Z_8=Z_9=-0.7 I_2.
\end{equation*}

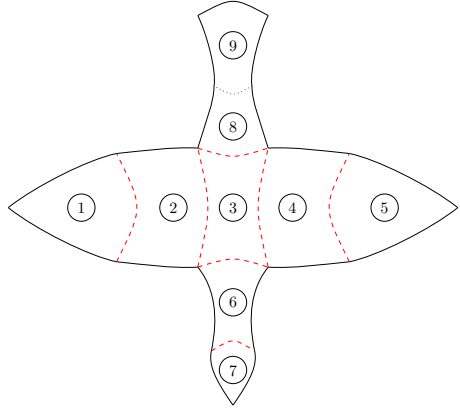
\begin{figure}[htbp]
	\begin{center}
		\scalebox{0.6}{
			\begin{tikzpicture}[x=0.6cm,y=0.6cm]
				
				\draw (-8.3,0.0) .. controls (-7.3,-0.8) and (-5.3,-1.8) .. (-4.3,-2.0);
				\draw [dashed,red] (-4.3,-2.0) .. controls (-3.3,0.0) .. (-4.3,2.0);
				\draw (-4.3,2.0) .. controls (-5.3,1.8) and (-7.3,0.8) .. (-8.3,0.0);
				
				\draw (-4.3,-2.0) .. controls (-3.3,-2.1) and (-2.3,-2.25) .. (-1.3,-2.2);
				\draw [dashed,red] (-1.3,-2.2) .. controls (-0.8,0.0) .. (-1.3,2.2);
				\draw (-1.3,2.2) .. controls (-2.3,2.25) and (-3.3,2.1) .. (-4.3,2.0);
				
				\draw (8.3,0.0) .. controls (7.3,-0.8) and (5.3,-1.8) .. (4.3,-2.0);
				\draw [dashed,red] (4.3,-2.0) .. controls (3.3,0.0) .. (4.3,2.0);
				\draw (4.3,2.0) .. controls (5.3,1.8) and (7.3,0.8) .. (8.3,0.0);
				
				\draw (4.3,-2.0) .. controls (3.3,-2.1) and (2.3,-2.25) .. (1.3,-2.2);
				\draw [dashed,red] (1.3,-2.2) .. controls (0.8,0.0) .. (1.3,2.2);
				\draw (1.3,2.2) .. controls (2.3,2.25) and (3.3,2.1) .. (4.3,2.0);
				
				\draw [dashed,red] (-1.3,-2.2) .. controls (0.0,-1.8) .. (1.3,-2.2);
				\draw [dashed,red] (-1.3,2.2) .. controls (0.0,1.8) .. (1.3,2.2);
				
				\draw (-1.3,-2.2) .. controls (-0.4,-3.3) and (-0.7,-4.8) .. (-0.8,-5.3);
				\draw (-0.8,-5.3) .. controls (-0.9,-5.8) and (-0.7,-6.3) .. (0.0,-7.3);
				\draw (0.0,-7.3) .. controls (0.7,-6.3) and (0.9,-5.8) .. (0.8,-5.3);
				\draw (0.8,-5.3) .. controls (0.7,-4.8) and (0.4,-3.3) .. (1.3,-2.2);
				\draw [dashed,red] (-0.8,-5.3) .. controls (0.0,-4.8) .. (0.8,-5.3);		
				
				\draw (1.3,2.2) .. controls (1.0,3.2) and (0.7,3.9) .. (0.7,4.5);		
				\draw (0.7,4.5) .. controls (0.6,5.1) and (1.0,6.5) .. (1.3,7.1);
				\draw (1.3,7.1) .. controls (0.0,7.8) .. (-1.3,7.1);
				\draw (-1.3,7.1) .. controls (-1.0,6.5) and (-0.6,5.1) .. (-0.7,4.5);
				\draw (-0.7,4.5) .. controls (-0.7,3.9) and (-1.0,3.2) .. (-1.3,2.2);
				
				\draw [dotted] (0.7,4.5) .. controls (0.0,4.1) .. (-0.7,4.5);
				
				\draw (-5.95,0.0) node[right]{$1$};
				\draw (-5.6,0.0) circle (2ex);
				\draw (-2.55,0.0) node[right]{$2$};
				\draw (-2.2,0.0) circle (2ex);
				\draw (-0.35,0.0) node[right]{$3$};
				\draw (0.0,0.0) circle (2ex);
				\draw (2.55,0.0) node[left]{$4$};
				\draw (2.2,0.0) circle (2ex);
				\draw (5.95,0.0) node[left]{$5$};
				\draw (5.6,0.0) circle (2ex);
				
				\draw (0.0,-3.1) node[below]{$6$};
				\draw (0.0,-3.5) circle (2ex);
				\draw (0.0,-5.65) node[below]{$7$};
				\draw (0.0,-6.0) circle (2ex);
				
				\draw (0.0,2.6) node[above]{$8$};
				\draw (0.0,3.0) circle (2ex);
				\draw (0.0,5.6) node[above]{$9$};
				\draw (0.0,6.0) circle (2ex);
				
			\end{tikzpicture}
		}
	\end{center}
	\caption{Computational domain for the \emph{bird} experiment.} \label{fig:geo_bird}
\end{figure}

The mesh is constituted of $1468$ quadrilaterals with maximal mesh size $0.690837$ (total of $44040$ degrees of freedom, $39636=1468\cdot9\cdot3$ for the deformation $\vy_h$ and $4404=1468\cdot 3$ for the Lagrange multipliers used to enforce the linearized isometry constraint), see Figure~\ref{fig:mesh_bird}. 
\begin{figure}[htbp]
	\begin{center}
		\includegraphics[width=9cm]{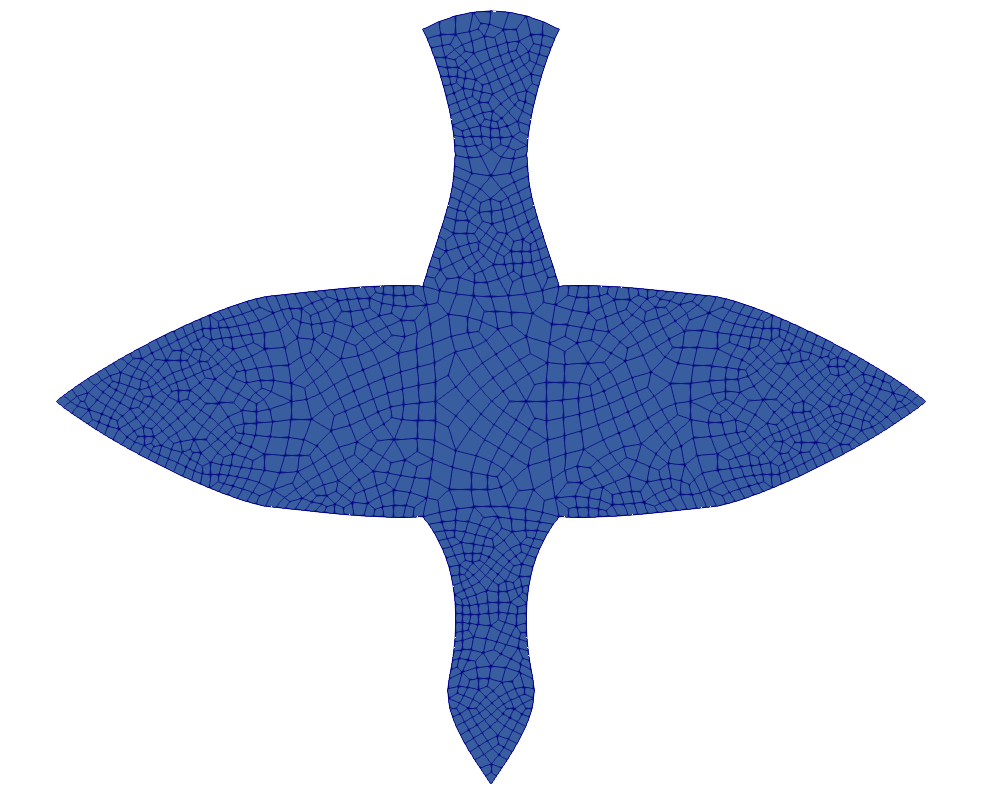}
		\caption{Mesh for the \emph{bird} example.} \label{fig:mesh_bird}
	\end{center}
\end{figure}
The initial deformation $\vy_h^0$ is again taken to be the identity map $\vy_h^0(\Omega)=\Omega\times\{0\}$ for which $E_h^{\rm bil}=12.9388$ and $\mathcal D_h^{\rm bary}=3.0149\cdot 10^{-14}$. The deformation obtained after $100$, $200$, $500$ and $2383$ (final) iterations are given in Figure~\ref{fig:iter_bird}. For the latter the energy is $E_h^{\rm bil}=4.2239$ while the isometry defect is $\mathcal D_h^{\rm bary}=0.0529014$.

\begin{figure}[htbp]
	\begin{center}
		\includegraphics[width=6cm]{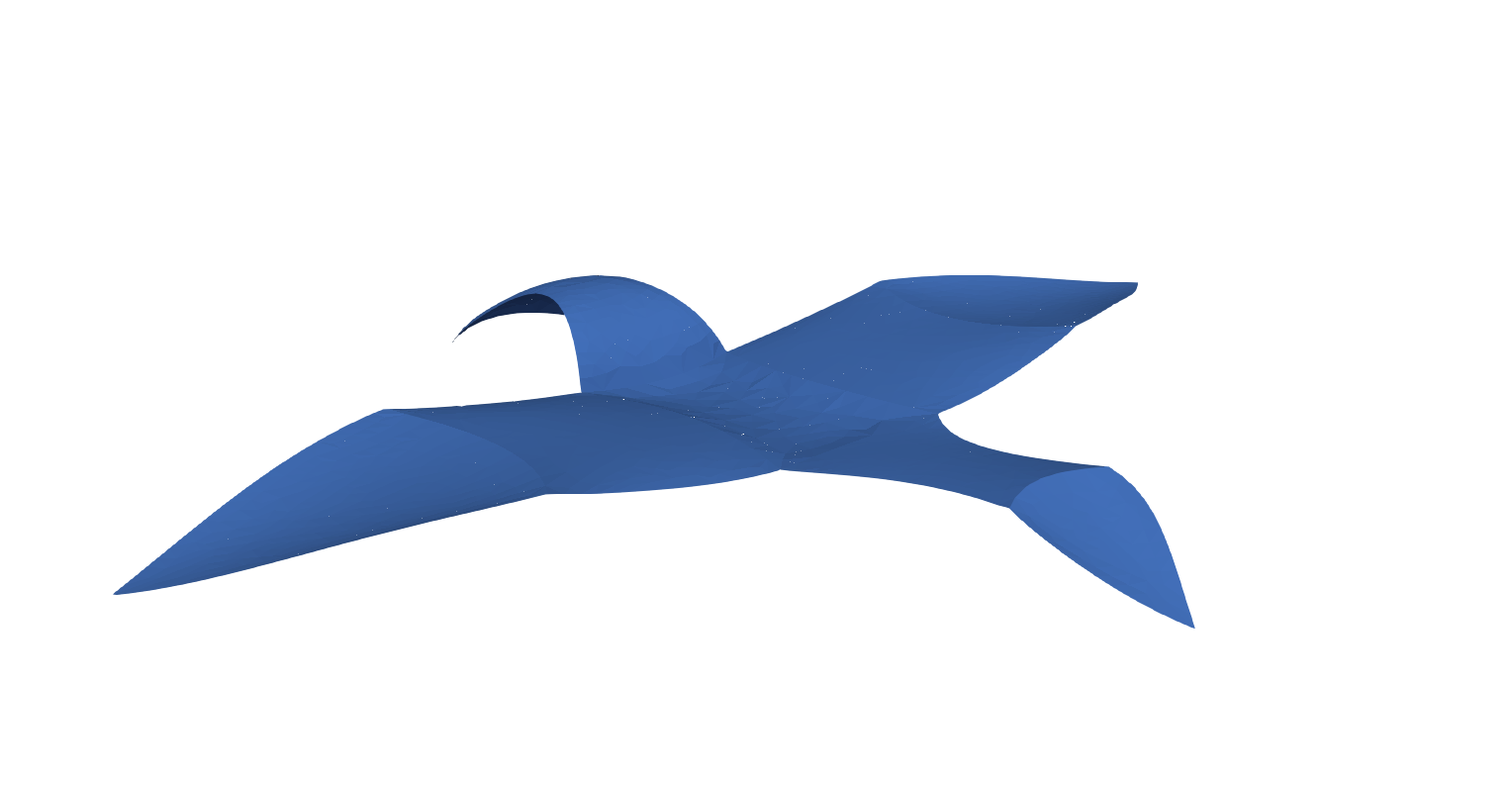}
		\includegraphics[width=6cm]{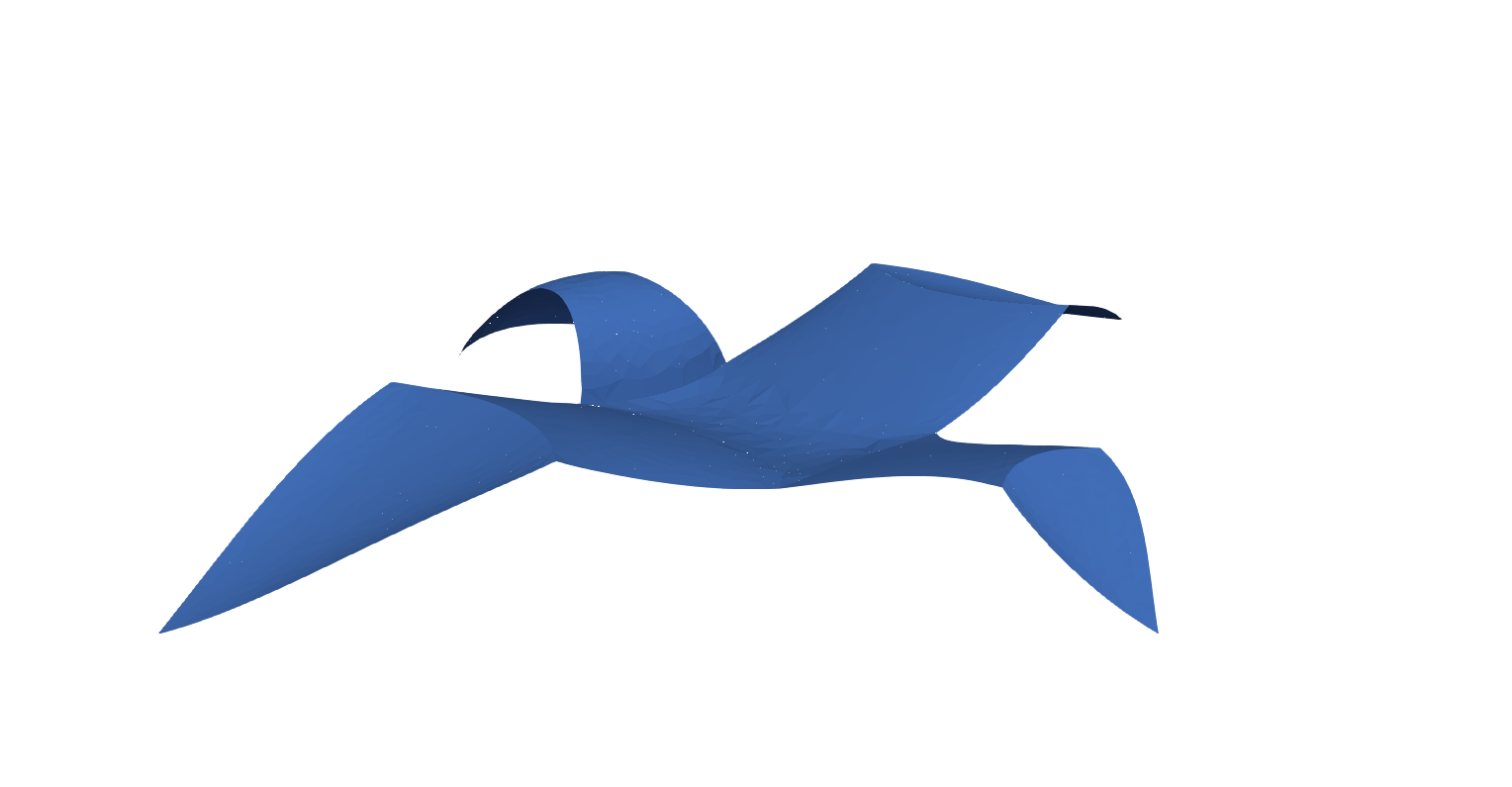} \\
		
		\includegraphics[width=6cm]{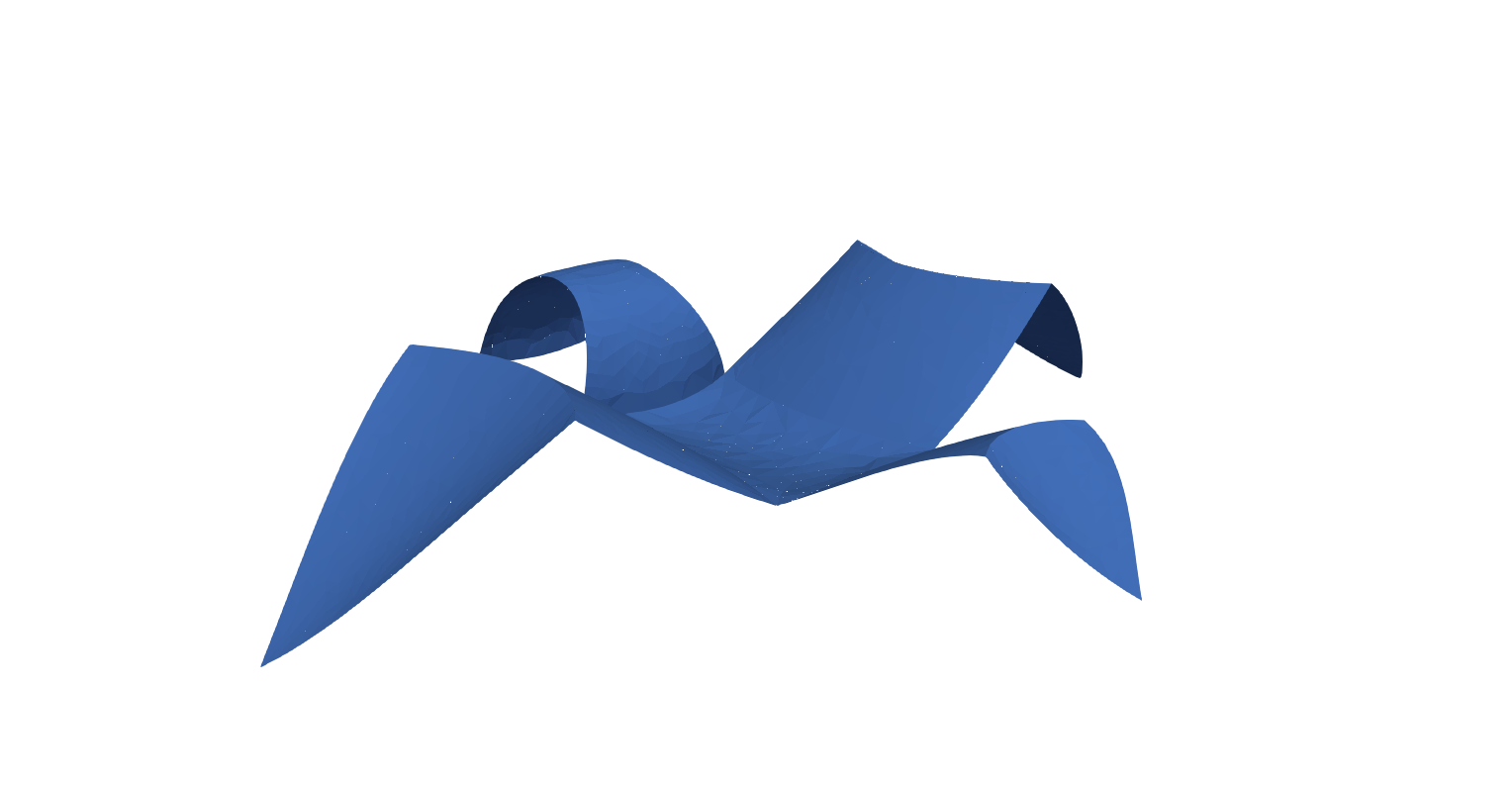}
		\includegraphics[width=6cm]{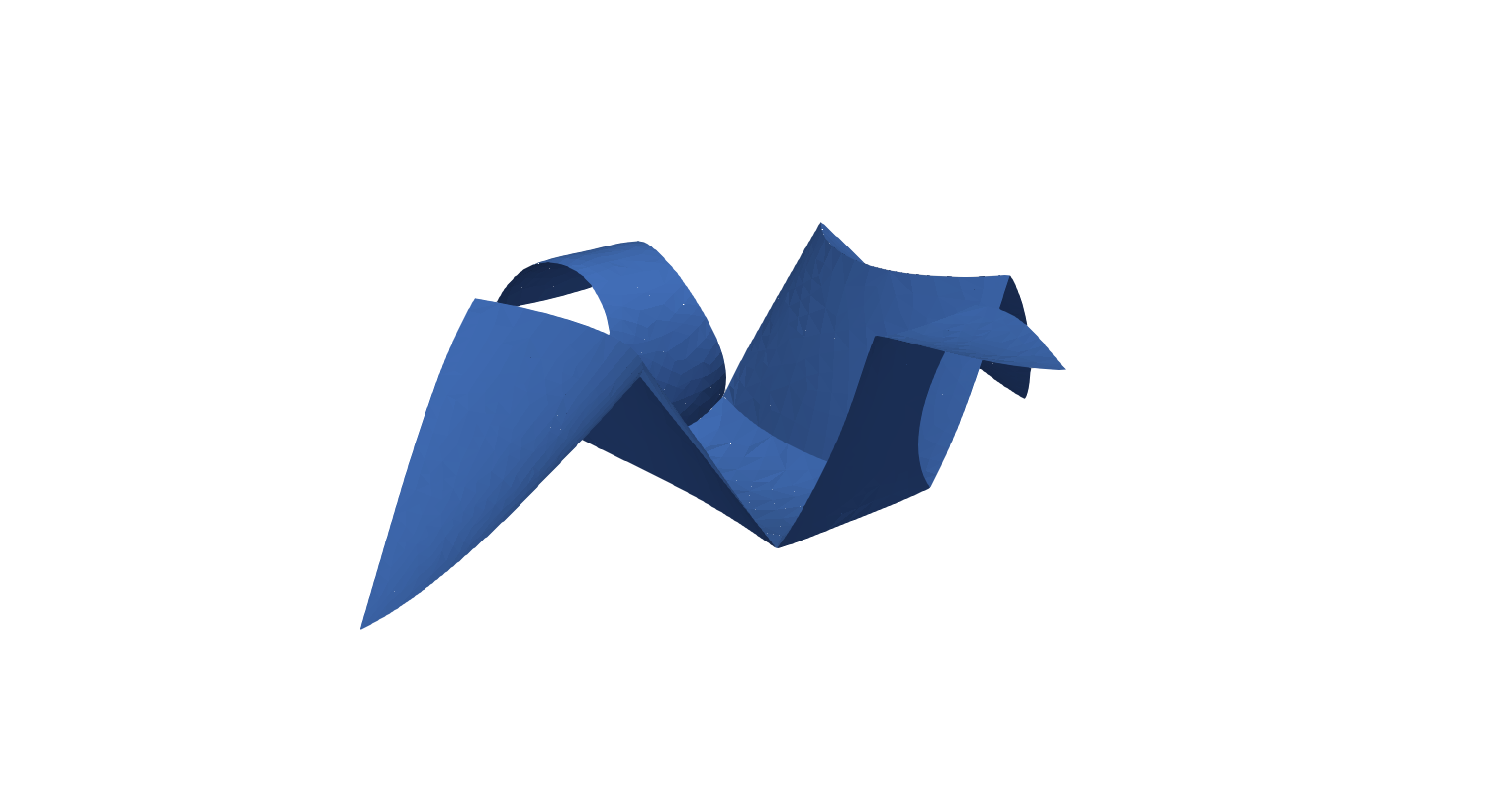}
	\end{center}
	\caption{Deformations for the \emph{bird} experiment. Left to right and top to bottom: deformation obtained after $100$, $200$, $500$, and $2383$ iterations.} \label{fig:iter_bird}
\end{figure}

It is worth pointing out that the equilibrium shape obtained is rather sensitive to the value of the spontaneous curvature tensors. To illustrate this, a coefficient $-0.4$ instead of $-0.3$ in $Z_2$ and $Z_4$ leads to an equilibrium deformation reported in Figure~\ref{fig:iter_bird_04}. After $2383$ iteration of the gradient flow, the deformation is similar to the one obtained in Figure~\ref{fig:iter_bird} (bottom-right), but the final deformation reached in $8085$ iterations is quite different. 
\begin{figure}[htbp]
	\begin{center}
		\includegraphics[width=6cm]{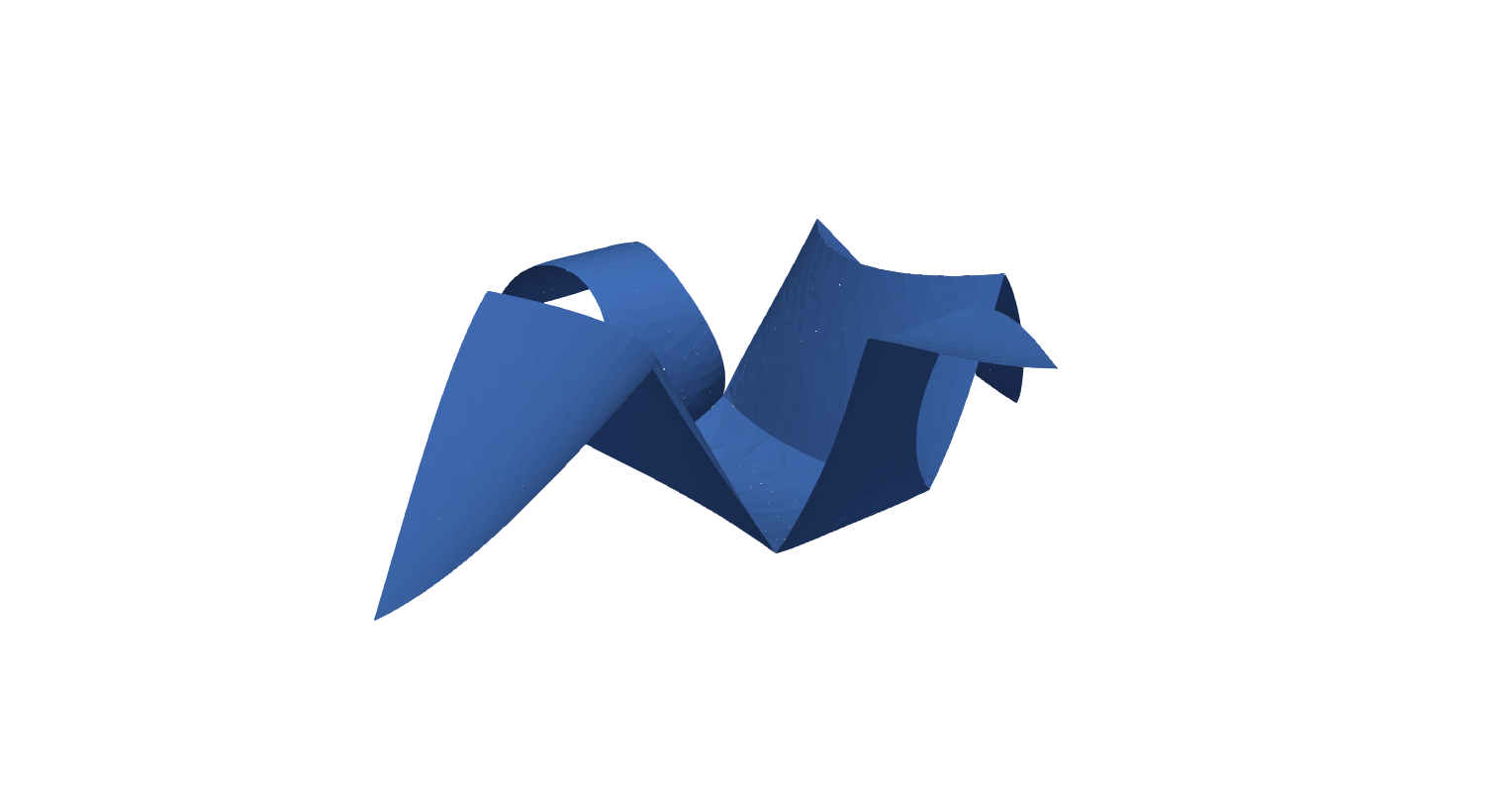}		
		\hspace*{-1cm}
		\includegraphics[width=6cm]{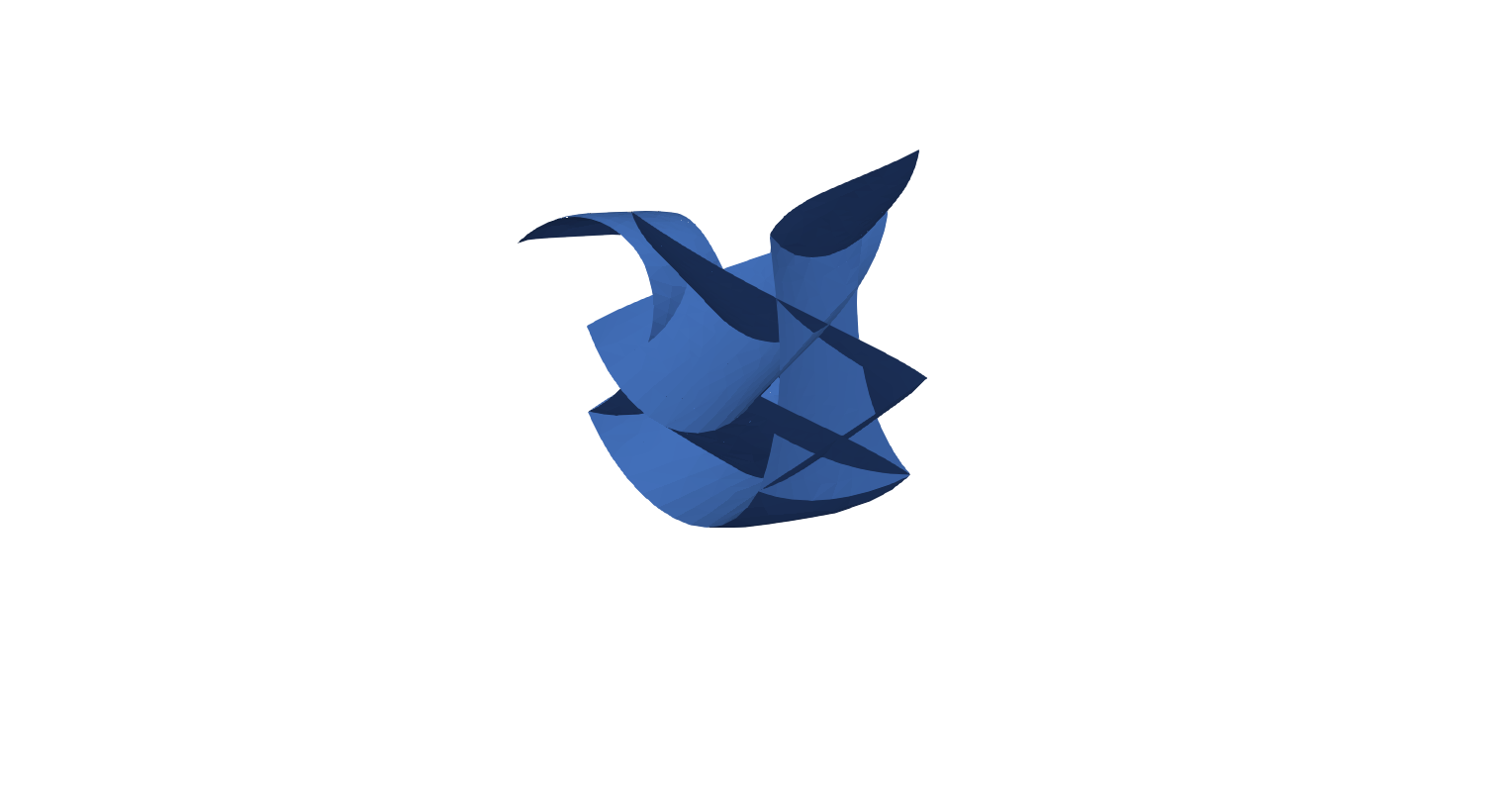}
	\end{center}
	\caption{Deformations for the \emph{bird} experiment when $Z_2=Z_4=[0, 0; 0, -0.4]$. (Left) Deformation obtained after $2383$ iterations; (right) final deformation reached in 8085 iterations.} \label{fig:iter_bird_04}
\end{figure}

\subsection{Single Layer with Time-Dependent Force} \label{sec:sphere} 

In this experiment, the domain $\Omega$ is the unit disc endowed with a prestrain metric corresponding to a \emph{half-sphere}, namely $g=\I(\vy)$ with
\begin{equation} \label{eqn:target_sphere}
	\vy(x_1,x_2)=\left(x_1,x_2,\sqrt{1+\epsilon-x_1^2-x_2^2}\right)^T, \quad \epsilon=10^{-3},    
\end{equation}
where the small parameter $\epsilon$ is introduced to avoid singularities at the boundary $\partial\Omega$. 
We study the effect of a uniform axial force directed towards the center of the sphere.
When the external force is sufficiently strong, we expect the sphere to crush like when the pressure inside a thin spherical reservoir is significantly smaller than the atmospheric pressure.

In this experiment, the external force is given by
\begin{equation} \label{def:force_sphere}
	\vf(x_1,x_2,t) = -t\vnu(x_1,x_2) = -\frac{t}{\sqrt{1+\epsilon}}\vy(x_1,x_2),
\end{equation}
where $\vnu$ is the structure outward pointing normal and $t$ denote the time.
We prescribe a mixed boundary condition with deformation given by $\vvarphi(x_1,x_2)=\left(x_1,x_2,\sqrt{\epsilon}\right)^T$ on $\Gamma^M=\partial \Omega$, which is compatible with \eqref{eqn:target_sphere}. 

The subdivision of $\Omega$ consists of 992 quadrilaterals with hanging nodes and of maximal mesh size $0.194084$ (total of 29760 degrees of freedom, $26874=992\cdot 9 \cdot 3$ for the deformation $\vy_h$ and $2976=992\cdot 3$ for the Lagrange multipliers used to enforce the linearized metric constraint). It is finer near the boundary where more stretching is expected. The initial deformation $\vy_h^0$ is taken to be the continuous Lagrange interpolant of $\vy$ in \eqref{eqn:target_sphere}, see Figure~\ref{fig:sphere_init}, for which $E_h^{\rm pre}=561.443$ and $\mathcal D_h^{\rm aver}=0.169394$. The other numerical parameters are $\tau=0.2$, $tol=10^{-3}$, and $\Delta t=5$ so that $t_m=m\Delta t =5m$. 

\begin{figure}[htbp]
	\begin{center}
		\includegraphics[width=6cm]{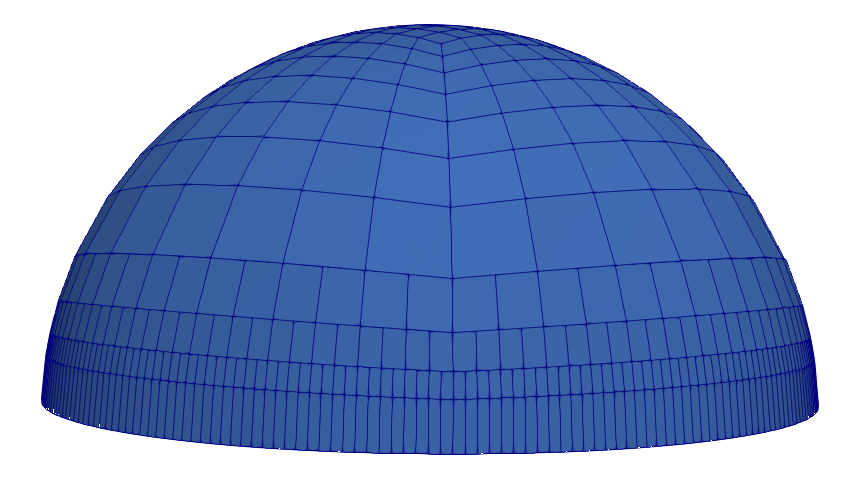}
	\end{center}
	\caption{Initial deformation for the \emph{half-sphere} experiment.} \label{fig:sphere_init}
\end{figure}

For $m=1,\ldots,44$, the equilibrium deformation is similar to the initial one given in Figure~\ref{fig:sphere_init}. The number of gradient flow iterations are 1594 for $m=1$ and between 8 to 101 for $2\le m\le 44$. When $m=45$, $\vf(x_1,x_2,t_m) = -225\vnu(x_1,x_2)$ and the object cannot sustain the corresponding force anymore. It collapses and deforms all the way to a \emph{half-sphere} in the opposite direction. In particular, the range of values for the third component of the deformation is $[0.03094,0.97704]$ at the first iteration and $[-1.06951,0.03254]$ at the last iteration of the gradient flow (81484). We refer to Figure~\ref{fig:sphere_evolution} for the deformation obtained at different steps of the gradient flow and to Table~\ref{table:energy_defect_sphere} for the corresponding energy metric defect. 

\begin{figure}[htbp]
	\begin{center}
		\includegraphics[width=6cm]{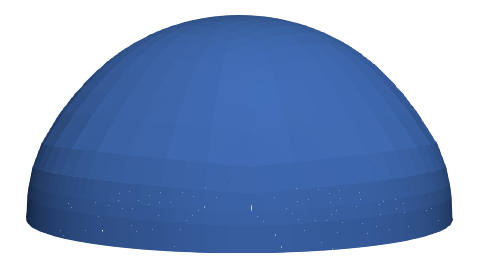}
		\includegraphics[width=6cm]{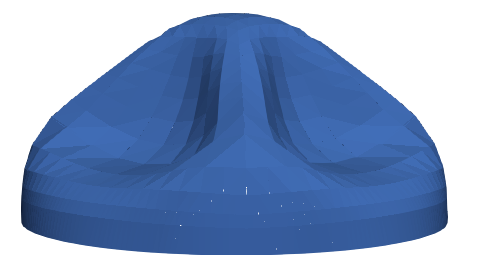}
		\includegraphics[width=6cm]{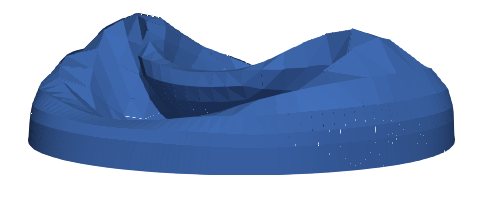}
		\includegraphics[width=6cm]{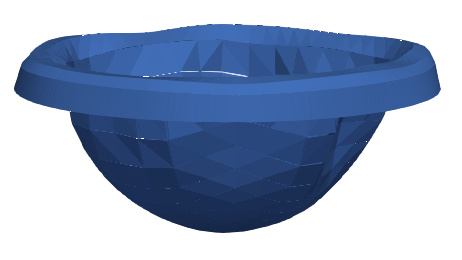}
		\includegraphics[width=6cm]{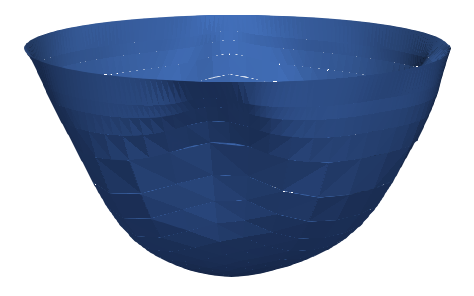}
	\end{center}
	\caption{Deformation for the \emph{half-sphere} experiment when $t=225$ in the forcing term \eqref{def:force_sphere}. Left to right and top to bottom: initial deformation and deformation obtained after 400, 800, 1200, and 81484 steps of the gradient flow.} 
	\label{fig:sphere_evolution}
\end{figure}

\begin{table}[htbp!] \label{table:energy_defect_sphere}
	\begin{center}
		\begin{tabular}{ |r|r|c| } 
			\hline
			\multicolumn{1}{|c|}{$n$} & \multicolumn{1}{c|}{$E_h^{\rm pre}$} & $\mathcal{D}_h^{\rm aver}$ \\
			\hline
			0 & 738.520 & 0.206666 \\
			\hline
			400 & 711.633 & 0.245251 \\
			\hline
			800 & 597.757 & 0.351605 \\
			\hline
			1200 & 313.706 & 0.524551 \\
			\hline 
			81484 & 49.763 & 0.528641 \\
			\hline
		\end{tabular}
		\caption{Prestrain energy $E_h^{\rm pre}$ and metric defect $\mathcal{D}_h^{\rm aver}$ at different steps of the gradient flow for the \emph{half-sphere} experiment when $t=225$, see Figure~\ref{fig:sphere_evolution} for the corresponding deformations.}
	\end{center}
\end{table}

\subsection{Starshade Experiment} \label{sec:starshade}

The starshade technology was developed as part of the NASA exoplanet exploration program. 
Starshades are occulters external to a telescope that shade the light from stars when imaging planets \cite{starshade}. They fold to minimize the space they occupy when transported in rockets and are easily deployable when arrived at destination.

The mathematical model corresponds to the single layer case with isometry constraint and time-dependent boundary conditions mimicking compression (closing) and decompression (opening). The geometry is taken from \cite{NASA} and consists of a dodecagon, see Figure~\ref{fig:geo_nasa}.
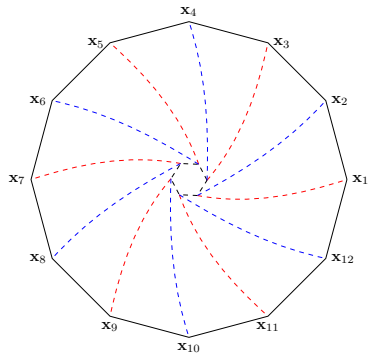
\begin{figure}[htbp]
	\begin{center}
		\scalebox{0.6}{
			\begin{tikzpicture}[x=0.5cm,y=0.5cm]
				
				\draw (7.0,0.0) -- (6.06218,3.5) -- (3.5,6.06218) -- (0.0,7.0) -- (-3.5,6.06218) -- (-6.06218,3.5) -- (-7.0,0.0) -- (-6.06218,-3.5) -- (-3.5,-6.06218) -- (-0.0,-7.0) -- (3.5,-6.06218) -- (6.06218,-3.5) -- (7.0,0.0);
				
				\draw [dashed] (0.799239,-0.034896) -- (0.429840,0.674713) -- (-0.369399,0.709609) -- (-0.799239,0.034896) -- (-0.429840,-0.674713) -- (0.369399,-0.709609) -- (0.799239,-0.034896);
				
				\draw [dashed,red] (0.369399,-0.709609) .. controls (2.7716,-1.1481) and (4.95722,-0.65263) .. (7.0,0.0);
				
				\draw [dashed,blue] (0.369399,-0.709609) .. controls (2.92889,0.64932) and (4.6194,1.9134) .. (6.06218,3.5);	
				
				\draw [dashed,red] (0.799239,-0.034896) .. controls (2.3801,1.8263) and (	3.0438,3.9668) .. (3.5,6.06218);
				
				\draw [dashed,blue] (0.799239,-0.034896) .. controls (0.90212,2.86115) and (0.65263,4.95722) .. (0.0,7.0);	
				
				\draw [dashed,red] (0.429840,0.674713) .. controls (-0.39158,2.97433) and (-1.9134,4.6194) .. (-3.5,6.06218);
				
				\draw [dashed,blue] (0.429840,0.674713) .. controls (-2.0268,2.2118) and (	-3.9668,3.0438) .. (-6.06218,3.5);
				
				\draw [dashed,red] (-0.369399,0.709609) .. controls (-2.7716,1.1481) and (-4.95722,0.65263) .. (-7.0,0.0);
				
				\draw [dashed,blue] (-0.369399,0.709609) .. controls (-2.92889,-0.64932) and (-4.6194,-1.9134) .. (-6.06218,-3.5);
				
				\draw [dashed,red] (-0.799239,0.034896) .. controls (-2.3801,-1.8263) and (-3.0438,-3.9668) .. (-3.5,-6.06218);
				
				\draw [dashed,blue] (-0.799239,0.034896) .. controls (-0.90212,-2.86115) and (-0.65263,-4.95722) .. (-0.0,-7.0);
				
				\draw [dashed,red] (-0.429840,-0.674713) .. controls (0.39158,-2.97433) and (1.9134,-4.6194) .. (3.5,-6.06218);
				
				\draw [dashed,blue] (-0.429840,-0.674713) .. controls (2.0268,-2.2118) and (3.9668,-3.0438) .. (6.06218,-3.5);
				
				\draw (7.0,0.0) node[right]{$\vx_1$};
				\draw (6.06218,3.5) node[right]{$\vx_2$};
				\draw (3.5,6.06218) node[right]{$\vx_3$};
				\draw (0.0,7.0) node[above]{$\vx_4$};
				\draw (-3.5,6.06218) node[left]{$\vx_5$};
				\draw (-6.06218,3.5) node[left]{$\vx_6$};
				\draw (-7.0,0.0) node[left]{$\vx_7$};
				\draw (-6.06218,-3.5) node[left]{$\vx_8$};
				\draw (-3.5,-6.06218) node[below]{$\vx_9$};
				\draw (-0.0,-7.0) node[below]{$\vx_{10}$};
				\draw (3.5,-6.06218) node[below]{$\vx_{11}$};
				\draw (6.06218,-3.5) node[right]{$\vx_{12}$};		
			\end{tikzpicture}
		}
	\end{center}
	\caption{Domain for the \emph{starshade} example.} \label{fig:geo_nasa}
\end{figure}
Inside the computational domain, there are 12 cubic B\'ezier curves as well as 6 straight lines (hexagon) across which folding is possible. 
The red dashed curves correspond to \emph{valley folds} while the blue ones are \emph{mountain folds}. The nodes of the dodecagon are given by
$$\vx_i = \Big(R\cos\Big(\frac{(i-1)\pi}{6}\Big),R\sin\Big(\frac{(i-1)\pi}{6}\Big)\Big), \quad i=1,2,\ldots,12,$$
with $R=7$, while the nodes of the hexagon are given by
$$\vz_i = \Big(r\cos\Big(\frac{(i-1)\pi}{3}-\frac{\pi}{72}\Big),r\sin\Big(\frac{(i-1)\pi}{3}-\frac{\pi}{72}\Big)\Big), \quad i=1,2,\ldots,6,$$
with $r=0.8$. For the curves, we use the following control points ($i=1,2,\ldots,6$): 
$$\vp_i^V = \Big(3\cos\Big(\frac{(i-1)\pi}{3}-\frac{\pi}{8}\Big),3\sin\Big(\frac{(i-1)\pi}{3}-\frac{\pi}{8}\Big)\Big)$$
$$\vq_i^V = \Big(5\cos\Big(\frac{(i-1)\pi}{3}-\frac{\pi}{24}\Big),5\sin\Big(\frac{(i-1)\pi}{3}-\frac{\pi}{24}\Big)\Big)$$
$$\vp_i^M = \Big(3\cos\Big(\frac{(i-1)\pi}{3}+\frac{5\pi}{72}\Big),3\sin\Big(\frac{(i-1)\pi}{3}+\frac{5\pi}{72}\Big)\Big)$$
and
$$\vq_i^M = \Big(5\cos\Big(\frac{(i-1)\pi}{3}+\frac{\pi}{8}\Big),5\sin\Big(\frac{(i-1)\pi}{3}+\frac{\pi}{8}\Big)\Big)$$
The first valley curve is obtained using $\vz_1$, $\vp_1^V$, $\vq_1^V$ and $\vx_1$, the first mountain curve uses $\vz_1$, $\vp_1^M$, $\vq_1^M$ and $\vx_2$, the second valley curve uses $\vz_2$, $\vp_2^V$, $\vq_2^V$ and $\vx_3$, and so on.

For the boundary conditions, we compress the 6 \emph{valley points} $\vx_i$, $i=1,3,5,7,9,11$. More precisely, we consider a dynamical setting by prescribing for $i=1,3,5,7,9,11,$ the following time-dependent pointwise boundary condition for the deformation $\vy_h$
\begin{equation} \label{BC:nasa}
	\vvarphi_i(t) = (1-2c_rt)\Big(R\cos\Big(\frac{(i-1)\pi}{6}\Big),R\sin\Big(\frac{(i-1)\pi}{6}\Big),0\Big)^T,
\end{equation}
where $c_r=0.25$ is a compression ratio. 
Using the strategy described in Section~\ref{subsec:dynamics}, we compute an approximate equilibrium deformation $\vy_h^m$ at time $t=t_m$, $m=0,1,\ldots,M$. Recall that this entails running the BC preprocessing for the variation $\delta \hat \vy_h^{m}= \hat \vy_h^{m}-\hat \vy_h^{m-1}$ imposing the pointwise conditions
\begin{equation*}
	\delta \vvarphi_i = -2c_r\Delta t\Big(R\cos\Big(\frac{(i-1)\pi}{6}\Big),R\sin\Big(\frac{(i-1)\pi}{6}\Big),0\Big)^T, \quad i=1, 3,5,7,9,11.
\end{equation*}
Note that at $t=0$, $\vvarphi_i=\vx_i$ and in particular, the BC preprocessing with the above boundary conditions results in a flat deformation which is conserved by the metric preprocessing and gradient flow. In order to generate an out of plane deformation, the boundary conditions for the initialization step (see case $t=0$ in Section~\ref{subsec:dynamics}) are modified to be \eqref{BC:nasa} for the \emph{valley} points and $\vvarphi_i=(\vx_i,1.5)^T$ for the 6 \emph{mountain} points $\vx_i$, $i=2,4,6,8,10,12$.

The time-step for the time discretization is $\Delta t = 0.05$, and $\mathtt{M}=35$ steps are performed.
The gradient flow parameters are $\tau=0.05$ and $tol=0.1$ while $\widetilde \tau=0.01$, $\widetilde{tol}=0.5$, and $\widetilde\varepsilon_0=0.5$ are chosen for the metric preprocessing.
The computational domain subdivision consists of $620$ quadrilaterals of maximum mesh size $1.20925$ (total of $18600$ degrees of freedom, $16740=620\cdot9\cdot3$ for the deformation $\vy_h$ and $1860=620\cdot 3$ for the Lagrange multipliers used to enforce the linearized isometry constraint), see Figure~\ref{fig:mesh_nasa}. The penalty parameters are $\gamma_0=\gamma_1=\gamma_2=10$. 
\begin{figure}[htbp]
	\begin{center}
		\includegraphics[width=5cm]{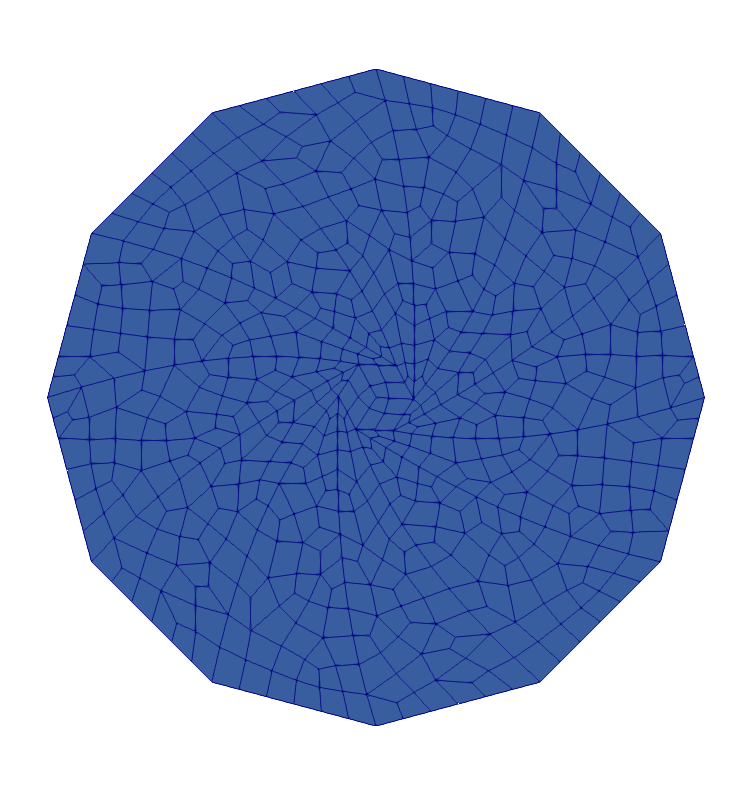}
		\caption{Mesh for the \emph{starshade} example.} \label{fig:mesh_nasa}
	\end{center}
\end{figure}

The equilibrium deformations obtained for several time $t_m=m\Delta t$ are reported in Figure~\ref{fig:iter_NASA}, see Figure~\ref{fig:iter_NASA_side} for corresponding side views.
\begin{figure}[htbp]
	\begin{center}
		\includegraphics[width=6cm]{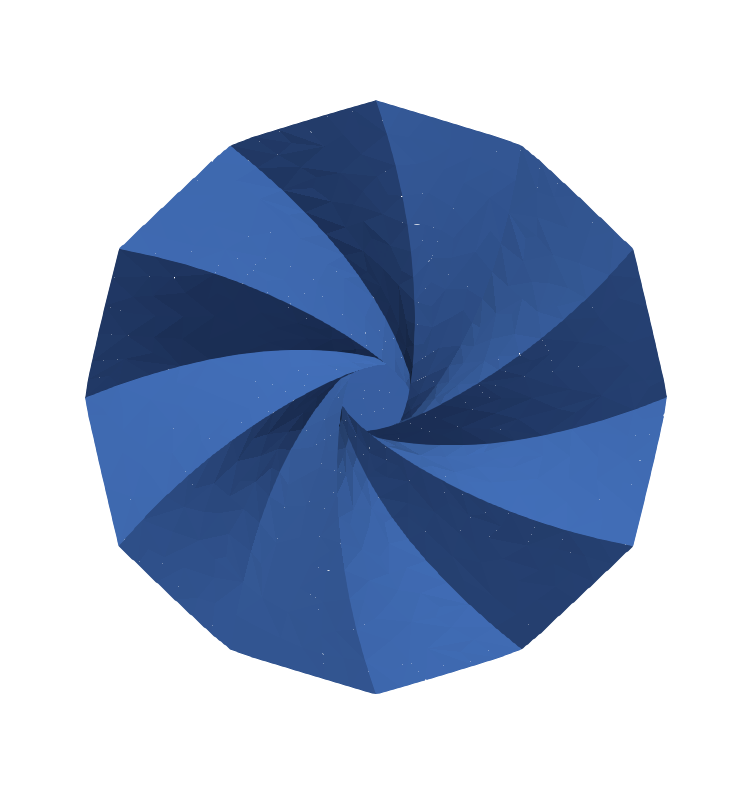}
		\includegraphics[width=6cm]{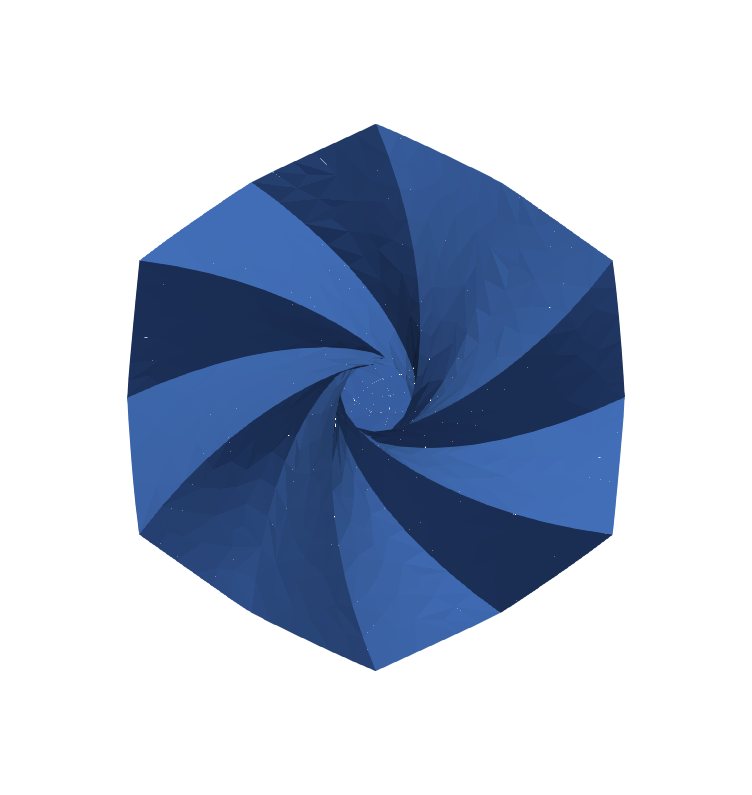} \\		
		\includegraphics[width=6cm]{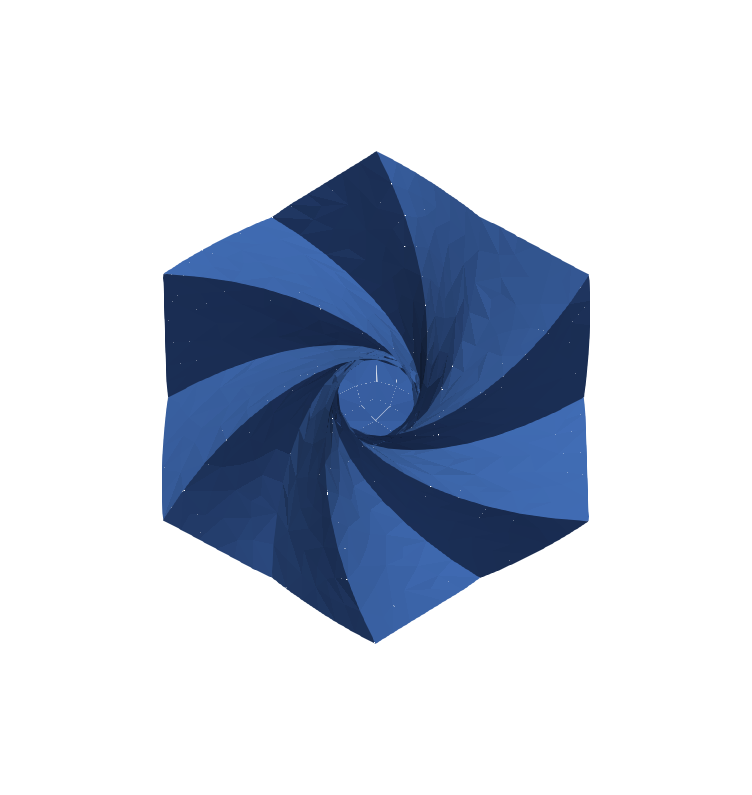}
		\includegraphics[width=6cm]{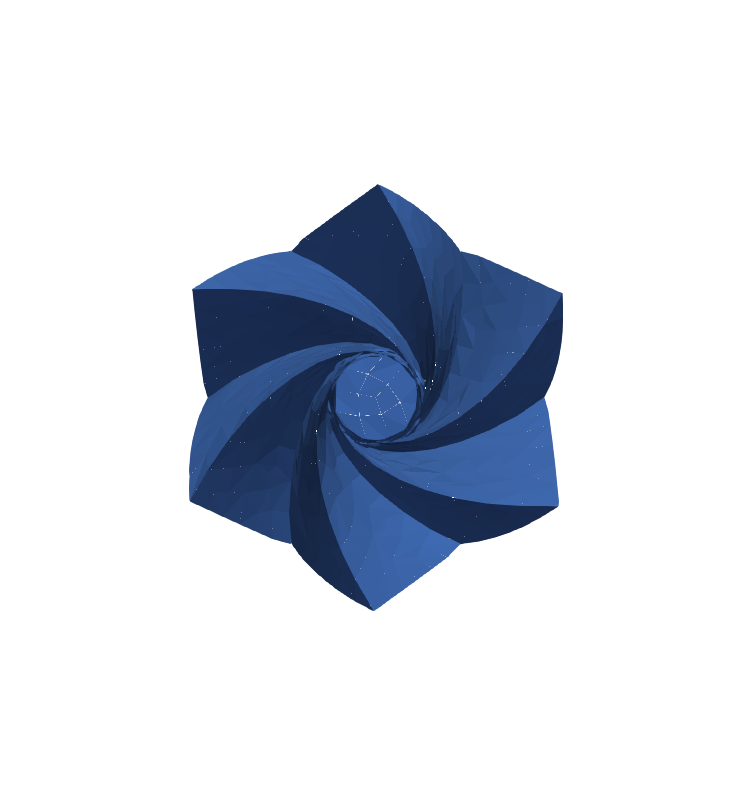} \\
		\includegraphics[width=6cm]{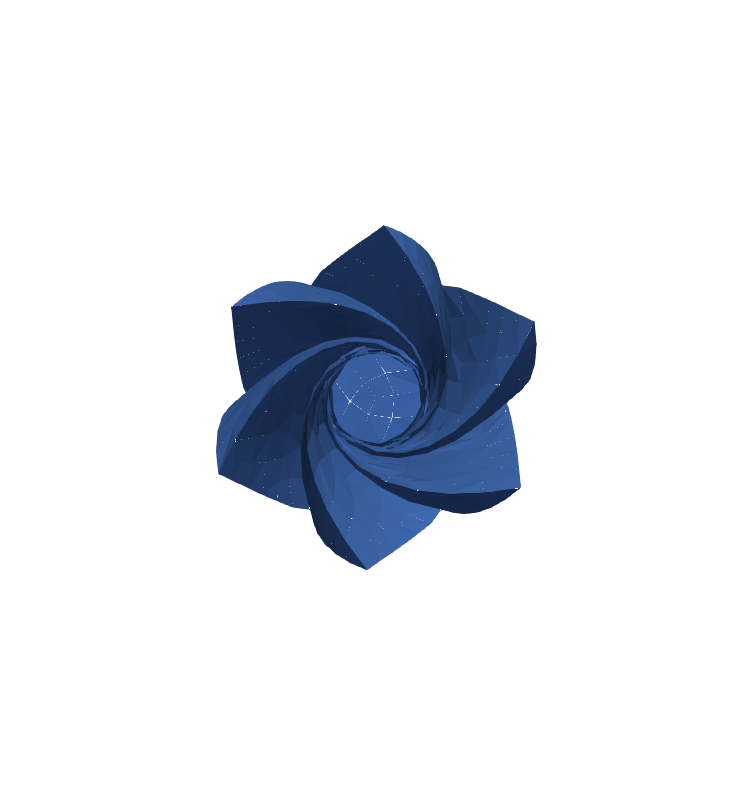}
		\hspace*{-2.8cm}\includegraphics[width=6cm]{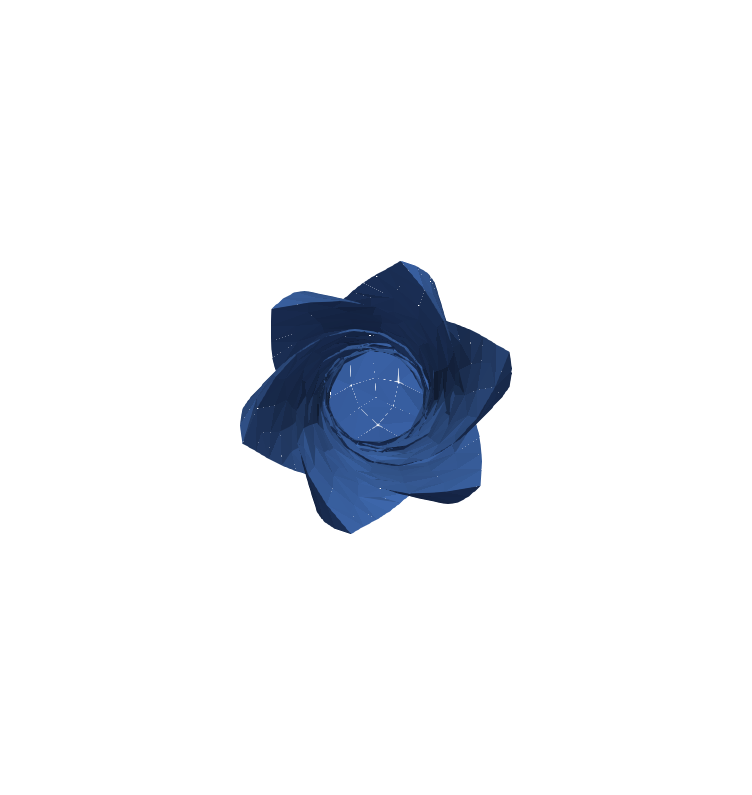}\hspace*{-2.8cm}\includegraphics[width=6cm]{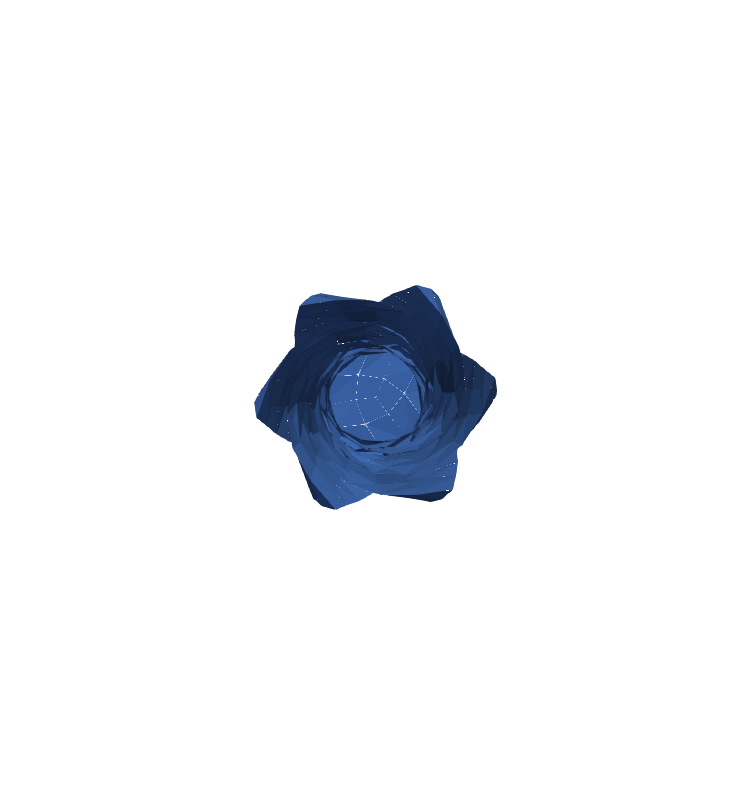}
	\end{center}
	\caption{Left to right and top to bottom: \emph{equilibrium} deformations for the \emph{starshade} experiment at time $t_m$ for $m=5,10,15,20,25,30,35$.} \label{fig:iter_NASA}
\end{figure}

\begin{figure}[htbp]
	\begin{center}
		\includegraphics[width=0.3\textwidth]{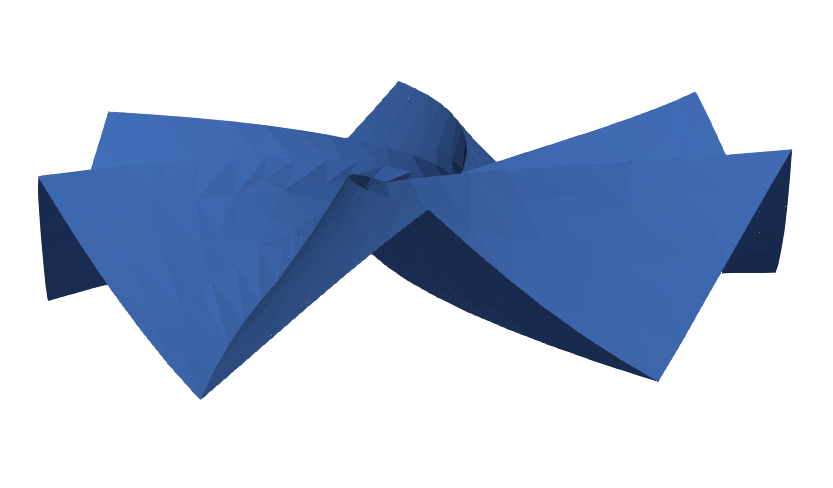}
		\includegraphics[width=0.3\textwidth]{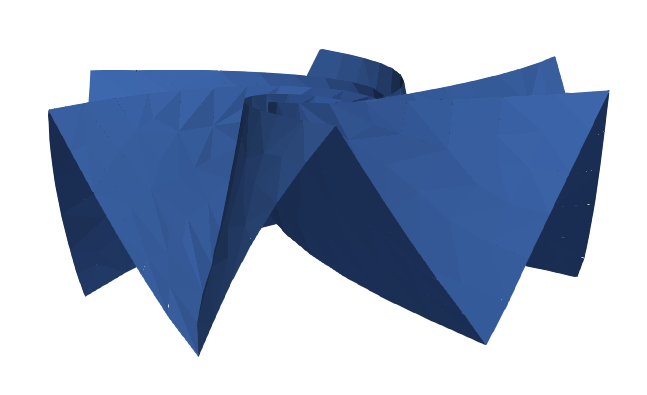} 		
		\includegraphics[width=0.3\textwidth]{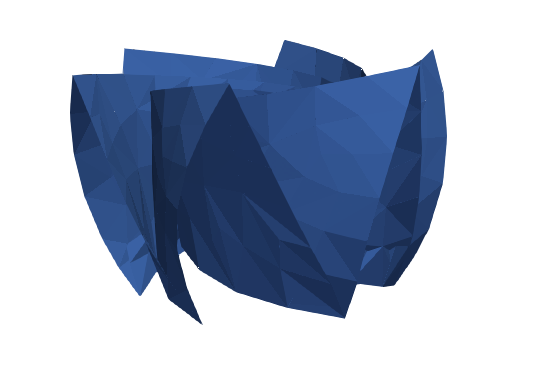}
	\end{center}
	\caption{Left to right: \emph{equilibrium} deformations for the \emph{starshade} experiment at time $t_m$ for $m=10,20,30$ (side view).} \label{fig:iter_NASA_side}
\end{figure}

\bibliographystyle{plain} 
\bibliography{bibliography}

\begin{thebibliography}{10}

\bibitem{ABS2011}
Silas Alben, Bavani Balakrisnan, and Elisabeth Smela.
\newblock Edge effects determine the direction of bilayer bending.
\newblock {\em Nano Lett.}, 11(6):2280--2285, 2011.

\bibitem{arnold2002unified}
Douglas~N. Arnold, Franco Brezzi, Bernardo Cockburn, and L.~Donatella Marini.
\newblock Unified analysis of discontinuous {G}alerkin methods for elliptic
  problems.
\newblock {\em SIAM J. Numer. Anal.}, 39(5):1749--1779, 2002.

\bibitem{BNS2015}
Bavani Balakrisnan, Alek Nacev, and Elisabeth Smela.
\newblock Design of bending multi-layer electroactive polymer actuators.
\newblock {\em Smart Mater. Struct.}, 24(4):045032, 2015.

\bibitem{bartels2013}
S{\"o}ren Bartels.
\newblock Finite element approximation of large bending isometries.
\newblock {\em Numer. Math.}, 124(3):415--440, 2013.

\bibitem{bartels2020finite}
S{\"o}ren Bartels.
\newblock Finite element simulation of nonlinear bending models for thin
  elastic rods and plates.
\newblock In {\em Handbook of Numerical Analysis}, volume~21, pages 221--273.
  Elsevier, 2020.

\bibitem{BBH2021}
S{\"o}ren Bartels, Andrea Bonito, and Peter Hornung.
\newblock Modeling and simulation of thin sheet folding.
\newblock {\em Interfaces Free Bound.}, 24(4):459--485, 2022.

\bibitem{BBN2017}
S{\"o}ren Bartels, Andrea Bonito, and Ricardo~H. Nochetto.
\newblock Bilayer plates: Model reduction, {$\Gamma$}‐convergent finite
  element approximation, and discrete gradient flow.
\newblock {\em Comm. Pure Appl. Math.}, 70(3):547--589, 2017.

\bibitem{bartels2022error}
S{\"o}ren Bartels, Andrea Bonito, and Philipp Tscherner.
\newblock Error estimates for a linear folding model.
\newblock {\em IMA J. Numer. Anal.}, 2023.

\bibitem{BP2021}
S{\"o}ren Bartels and Christian Palus.
\newblock Stable gradient flow discretizations for simulating bilayer plate
  bending with isometry and obstacle constraints.
\newblock {\em IMA J. Numer. Anal.}, 42(3):1903--1928, 2021.

\bibitem{bassi1997high}
Francesco Bassi and Stefano Rebay.
\newblock A high-order accurate discontinuous finite element method for the
  numerical solution of the compressible {N}avier--{S}tokes equations.
\newblock {\em J. Comput. Phys.}, 131(2):267--279, 1997.

\bibitem{BALG2010}
Noy Bassik, Beza Abebe, Kate Laflin, and David Gracias.
\newblock Photolithographically patterned smart hydrogel based bilayer
  actuators.
\newblock {\em Polymer}, 51(26):6093--6098, 2010.

\bibitem{BT09}
Amir Beck and Marc Teboulle.
\newblock A fast iterative shrinkage-thresholding algorithm for linear inverse
  problems.
\newblock {\em SIAM J. Imaging Sci.}, 2(1):183--202, 2009.

\bibitem{BK2014}
Peter Bella and Robert Kohn.
\newblock Metric-induced wrinkling of a thin elastic sheet.
\newblock {\em J. Nonlinear Sci.}, 24:1147--1176, 2014.

\bibitem{BLS2016}
Kaushik Bhattacharya, Marta Lewicka, and Mathia Sch{\"a}ffner.
\newblock Plates with incompatible prestrain.
\newblock {\em Arch. Rational Mech. Anal.}, 221(1):143--181, 2016.

\bibitem{BNPO2022}
Klaus B{\"o}hnlein, Stefan Neukamm, David Padilla-Garza, and Oliver Sander.
\newblock A homogenized bending theory for prestrained plates.
\newblock {\em J. Nonlinear Sci.}, 33(22):1--90, 2022.

\bibitem{BGMpreprint}
Andrea Bonito, Diane Guignard, and Angelique Morvant.
\newblock A note on the numerical approximation of thin structures.
\newblock {\em In preparation}, 2023.

\bibitem{BGNY2022a}
Andrea Bonito, Diane Guignard, Ricardo~H. Nochetto, and Shuo Yang.
\newblock {LDG} approximation of large deformations of prestrained plates.
\newblock {\em J. Comput. Phys.}, 448:110719, 2022.

\bibitem{BGNY2022b}
Andrea Bonito, Diane Guignard, Ricardo~H. Nochetto, and Shuo Yang.
\newblock Numerical analysis of the {LDG} method for large deformations of
  prestrained plates.
\newblock {\em IMA J. Numer. Anal.}, 43:627--662, 2023.

\bibitem{bonito2010quasi}
Andrea Bonito and Ricardo~H. Nochetto.
\newblock Quasi-optimal convergence rate of an adaptive discontinuous
  {G}alerkin method.
\newblock {\em SIAM J. Numer. Anal.}, 48(2):734--771, 2010.

\bibitem{BNN2020}
Andrea Bonito, Ricardo~H. Nochetto, and Dimitris Ntogkas.
\newblock Discontinuous {G}alerkin approach to large bending deformation of a
  bilayer plate with isometry constraint.
\newblock {\em J. Comput. Phys.}, 423:109785, 2020.

\bibitem{BNN2021}
Andrea Bonito, Ricardo~H. Nochetto, and Dimitris Ntogkas.
\newblock {DG} approach to large bending deformations with isometry constraint.
\newblock {\em Math. Models Methods Appl. Sci.}, 31(01):133--175, 2021.

\bibitem{BNY2022}
Andrea Bonito, Ricardo~H. Nochetto, and Shuo Yang.
\newblock {$\Gamma$}-convergent {LDG} method for large bending deformations of
  bilayer plates.
\newblock {\em arXiv preprint arXiv:2301.03151 [math.NA]}, 2023.

\bibitem{brezzi2000discontinuous}
Franco Brezzi, Gianmarco Manzini, Donatella Marini, Paola Pietra, and
  Alessandro Russo.
\newblock Discontinuous {G}alerkin approximations for elliptic problems.
\newblock {\em Numer. Methods Partial Differential Equations}, 16(4):365--378,
  2000.

\bibitem{CD15}
Antoine Chambolle and Charles~H. Dossal.
\newblock On the convergence of the iterates of ``{F}ast {I}terative
  {S}hrinkage/{T}hresholding algorithm''.
\newblock {\em J. Optim. Theory Appl.}, 166(3):25, 2015.

\bibitem{CS17}
Camillo De~Lellis and L{\'a}szl{\'o} Sz{\'e}kelyhidi~Jr.
\newblock High dimensionality and $h$-principle in {PDE}.
\newblock {\em Bull. Am. Math. Soc.}, 54:247--282, 2017.

\bibitem{DPE2010}
Daniele~A. Di~Pietro and Alexandre Ern.
\newblock Discrete functional analysis tools for discontinuous {G}alerkin
  methods with application to the incompressible {N}avier--{S}tokes equations.
\newblock {\em Math. Comp.}, 79(271):1303--1330, 2010.

\bibitem{DPE2011}
Daniele~A. Di~Pietro and Alexandre Ern.
\newblock {\em Mathematical Aspects of Discontinuous {G}alerkin Methods}.
\newblock Math{\'e}matiques et Applications. Springer Berlin Heidelberg, 2011.

\bibitem{ESK2009}
Efi Efrati, Eran Sharon, and Raz Kupferman.
\newblock Elastic theory of unconstrained non-{E}uclidean plates.
\newblock {\em J. Mech. Phys. Solids}, 57(4):762--775, 2009.

\bibitem{FSDM2005}
Yo{\"e}l Forterre, Jan Skotheim, Jacques Dumais, and L.~Mahadevan.
\newblock How the {V}enus flytrap snaps.
\newblock {\em Nature}, 433:421--425, 2005.

\bibitem{FJMM2003}
Gero Friesecke, Richard~D. James, Maria~G. Mora, and Stefan M{\"u}ller.
\newblock Derivation of nonlinear bending theory for shells from three
  dimensional nonlinear elasticity by {G}amma-convergence.
\newblock {\em C.R. Math.}, 336(8):697--702, 2003.

\bibitem{FJM2002c}
Gero Friesecke, Richard~D. James, and Stefan M{\"u}ller.
\newblock The f{\"o}ppl-von k{\'a}rm{\'a}n plate theory as a low energy
  {$\Gamma$}-limit of nonlinear elasticity.
\newblock {\em C.R. Math.}, 335(2):201--206, 2002.

\bibitem{FJM2002b}
Gero Friesecke, Richard~D. James, and Stefan M{\"u}ller.
\newblock A theorem on geometric rigidity and the derivation of nonlinear plate
  theory from three‐dimensional elasticity.
\newblock {\em C.R. Math.}, 55(11):1461--1506, 2002.

\bibitem{FJM2006}
Gero Friesecke, Richard~D. James, and Stefan M{\"u}ller.
\newblock A hierarchy of plate models derived from nonlinear elasticity by
  {G}amma-convergence.
\newblock {\em Arch. Rational Mech. Anal.}, 180(2):183--236, 2006.

\bibitem{FJM2002a}
Gero Friesecke, Stefan M{\"u}ller, and Richard~D. James.
\newblock Rigorous derivation of nonlinear plate theory and geometric rigidity.
\newblock {\em C.R. Math.}, 334(2):173--178, 2002.

\bibitem{goriely2005differential}
Alain Goriely and Martine Ben~Amar.
\newblock Differential growth and instability in elastic shells.
\newblock {\em Phys. Rev. Lett.}, 94(19):198103, 2005.

\bibitem{Gromov1986}
Mikhael Gromov.
\newblock {\em Partial Differential Relations}, volume~13.
\newblock Springer-Verlag, Berlin-Heidelberg, 1986.

\bibitem{HH2006}
Qing Han and Jia-Xing Hong.
\newblock {\em Isometric embedding of {R}iemannian manifolds in {E}uclidean
  spaces}, volume~13.
\newblock American Mathematical Soc., 2006.

\bibitem{JSI2000}
Edwin Jager, Elisabeth Smela, and Olle Ingan{\"a}s.
\newblock Microfabricating conjugated polymer actuators.
\newblock {\em Science}, 290(5496):1540--1545, 2000.

\bibitem{kim2012thermally}
Jungwook Kim, James~A. Hanna, Ryan~C. Hayward, and Christian~D. Santangelo.
\newblock Thermally responsive rolling of thin gel strips with discrete
  variations in swelling.
\newblock {\em Soft Matter}, 8(8):2375--2381, 2012.

\bibitem{Kirchhoff1850}
Gustav Kirchhoff.
\newblock \"uber das {G}leichgewicht und die {B}ewegung einer elastischen
  {S}cheibe.
\newblock {\em J. f{\"u}r Reine Angew. Math.}, 1850(40):51--88, 1850.

\bibitem{klein2007shaping}
Yael Klein, Efi Efrati, and Eran Sharon.
\newblock Shaping of elastic sheets by prescription of non-{E}uclidean metrics.
\newblock {\em Science}, 315(5815):1116--1120, 2007.

\bibitem{KVS2011}
Yael Klein, Shankar Venkataramani, and Eran Sharon.
\newblock An experimental study of shape transitions and energy scaling in thin
  non-{E}uclidean plates.
\newblock {\em Phys. Rev. Lett.}, 106(11):118303--118306, 2011.

\bibitem{Kuiper1955a}
Nicolaas~H. Kuiper.
\newblock On $c^1$-isometric imbeddings. i.
\newblock {\em Indag. Math. (Proceedings)}, 58:545--556, 1955.

\bibitem{Kuiper1955b}
Nicolaas~H. Kuiper.
\newblock On {$C^1$}-isometric imbeddings. ii.
\newblock {\em Indag. Math. (Proceedings)}, 58:683--689, 1955.

\bibitem{KS2014}
Raz Kupferman and Jack~P. Solomon.
\newblock A {R}iemannian approach to reduced plate, shell, and rod theories.
\newblock {\em J. Funct. Anal.}, 266:2989--3039, 2014.

\bibitem{DR1995}
Herv{\'e} Le~Dret and Annie Raoult.
\newblock The nonlinear membrane model as a variational limit of nonlinear
  three-dimensional elasticity.
\newblock {\em J. Math. Pures Appl.}, 73:549--578, 1995.

\bibitem{LM2021}
Marta Lewicka and L.~Mahadevan.
\newblock Geometry, analysis, and morphogenesis: Problems and prospects.
\newblock {\em Bull. Am. Math. Soc.}, 59(3):331--369, 2021.

\bibitem{LP2011}
Marta Lewicka and Mohammad Reza~Pakzad.
\newblock Scaling laws for non-{E}uclidean plates and the {$W^{2,2}$} isometric
  immersions of {R}iemannian metrics.
\newblock {\em ESAIM: Contr. Optim. C.A.}, 17(4):1158--1173, 2011.

\bibitem{LPFJ2021}
Huan Liu, Paul Plucinsky, Fan Feng, and Richard~D. James.
\newblock Origami and materials science.
\newblock {\em Philos. Trans. Royal Soc.}, 379(2201):20200113, 2021.

\bibitem{Love1906}
August~E.H. Love.
\newblock {\em A treatise on the mathematical theory of elasticity (2nd
  edition)}.
\newblock Cambridge University Press, 1906.

\bibitem{MS2019}
Cy~Maor and Asaf Shachar.
\newblock On the role of curvature in the elastic energy of non-{E}uclidean
  thin bodies.
\newblock {\em J. Elast.}, 134:149--173, 2019.

\bibitem{Menges2015}
Achim Menges and Steffen Reichert.
\newblock Performative wood: Physically programming the responsive architecture
  of the {H}ygro{S}cope and {H}ygro{S}kin projects.
\newblock {\em Archit. Des.}, 85(5):66--73, 2015.

\bibitem{modes2010disclination}
Carl~D. Modes, Kaushik Bhattacharya, and Mark Warner.
\newblock Disclination-mediated thermo-optical response in nematic glass
  sheets.
\newblock {\em Phys. Rev. E}, 81(6):060701, 2010.

\bibitem{modes2010gaussian}
Carl~D. Modes, Kaushik Bhattacharya, and Mark Warner.
\newblock Gaussian curvature from flat elastica sheets.
\newblock {\em Proc. Royal Soc.}, 467(2128):1121--1140, 2010.

\bibitem{starshade}
NASA.
\newblock Starshade technology development.
\newblock {\em https://exoplanets.nasa.gov/exep/technology/starshade/},
  accessed December 15, 2022.

\bibitem{Nash1954}
John Nash.
\newblock {$C^1$} isometric imbeddings.
\newblock {\em Ann. Math.}, 60(3):383--396, 1954.

\bibitem{Nesterov1983}
Yurii Nesterov.
\newblock A method for solving the convex programming problem with convergence
  rate $\mathcal{O}(1/k^2)$.
\newblock {\em Dokl. Akad. Nauk SSSR}, 269(9):543--547, 1983.

\bibitem{Nesterov2004}
Yurii Nesterov.
\newblock {\em Introductory Lectures on Convex Optimization: A Basic Course}.
\newblock Applied Optimization. Springer, 2004.

\bibitem{nesterov2018lectures}
Yurii Nesterov.
\newblock {\em Lectures on convex optimization}, volume 137 of {\em Springer
  Optimization and Its Applications}.
\newblock Springer, 2018.

\bibitem{nitsche1971variationsprinzip}
Joachim Nitsche.
\newblock {\"U}ber ein {V}ariationsprinzip zur {L}{\"o}sung von
  {D}irichlet-{P}roblemen bei {V}erwendung von {T}eilr{\"a}umen, die keinen
  {R}andbedingungen unterworfen sind.
\newblock {\em Abh. Math. Semin. Univ. Hambg.}, 36:9--15, 1971.

\bibitem{NASA}
Jet Propulsion Laboratory California~Institute of~Technology.
\newblock Space origami: Make your own starshade.
\newblock {\em
  https://www.jpl.nasa.gov/edu/learn/project/space-origami-make-your-own-starshade},
  accessed September 8, 2021.

\bibitem{Speck2020}
Renate Sachse, Anna Westermeier, Max Mylo, Joey Nadasdi, Manfred Bischoff,
  Thomas Speck, and Simon Poppinga.
\newblock Snapping mechanics of the {V}enus flytrap (\emph{Dionaea muscipula}).
\newblock {\em Proc. Natl. Acad. Sci.}, 117(27):16035--16042, 2020.

\bibitem{SLPSK2015}
Simon Schleicher, Julian Lienhard, Simon Poppinga, Thomas Speck, and Jan
  Knippers.
\newblock A methodology for tranferring principles of plant movements to
  elastic systems in architecture.
\newblock {\em Comput. Aided Des.}, 60:105--117, 2015.

\bibitem{S2007}
Bernd Schmidt.
\newblock Minimal energy configurations of strained multi-layers.
\newblock {\em Calc. Var. Partial Differential Equations}, 30(4):477--497,
  2007.

\bibitem{SIPL1993}
Elisabeth Smela, Olle Ingan{\"a}s, Qibing Pei, and Ingemar Lundstr{\"o}m.
\newblock Electrochemical muscles: Micromachining fingers and corkscrews.
\newblock {\em Adv. Mater. Lett.}, 5(9):630--632, 1993.

\bibitem{TCWSTBM2020}
Yasaman Tahouni, Tiffany Cheng, Dylan Wood, Renate Sachse, Rebecca Thierer,
  Manfred Bischoff, and Achim Menges.
\newblock Self-shaping curved folding: A 4{D}-printing method for fabrication
  of self-folding curved crease structures.
\newblock In {\em Symposium on Computational Fabrication}, number~5 in SCF '20,
  pages 1--11. Association for Computing Machinery, 2020.

\bibitem{WVMR2018}
Dylan Wood, Chiara Vailati, Achim Menges, and Markus R{\"u}ggeberg.
\newblock Hygroscopically actuated wood elements for weather responsive and
  self-forming building parts - facilitating upscaling and complex shape
  changes.
\newblock {\em Constr. Build. Mater.}, 165:782--791, 2018.

\bibitem{wu2013three}
Zi~L. Wu, Michael Moshe, Jesse Greener, Heloise Therien-Aubin, Zhihong Nie,
  Eran Sharon, and Eugenia Kumacheva.
\newblock Three-dimensional shape transformations of hydrogel sheets induced by
  small-scale modulation of internal stresses.
\newblock {\em Nat. Commun.}, 4:1586, 2013.

\bibitem{yavari2010geometric}
Arash Yavari.
\newblock A geometric theory of growth mechanics.
\newblock {\em J. Nonlinear Sci.}, 20(6):781--830, 2010.

\end{thebibliography}

\end{document}